\documentclass[10pt]{amsart}
\usepackage{amsmath}
\usepackage{xypic}
\usepackage{amssymb}
\usepackage{amsfonts}
\usepackage{amsthm}
\usepackage{enumerate}
\usepackage[all]{xy}
\usepackage[latin1]{inputenc}
\usepackage{graphicx}
\usepackage{latexsym}
\usepackage{accents}
\usepackage{mathdots}
\usepackage{mathrsfs}
\input xy

\makeatletter
    \newtheoremstyle{definition}
        {5pt}
        {3pt}
        {}
        {0pt}
        {\scshape}
        {.}
        {5pt}
        {\thmname{#1} \thmnumber{#2} \thmnote{[#3]}} 

\newtheoremstyle{theorems}
        {5pt}
        {3pt}
        {\itshape}
        {0pt}
        {\scshape}
        {.}
        {5pt}
        {\thmname{#1} \thmnumber{#2}\thmnote{[#3]}} 

\swapnumbers \theoremstyle{theorems}

\newtheorem{Theo}{Theorem}[section]
\newtheorem{Prop}[Theo]{Proposition}
\newtheorem{Cor}[Theo]{Corollary}
\newtheorem{Lemma}[Theo]{Lemma}
\newtheorem{Prop(BG)}[Theo]{Proposition (Bongartz-Gabriel)}
\newtheorem{Lemma(Asashiba)}[Theo]{Lemma(Asashiba)}
\newtheorem{Lemma(Gab)}[Theo]{Lemma(Gabriel)}
\newtheorem{Theo(Mil)}[Theo]{Theorem (Milicic)}

\theoremstyle{definition}

\newtheorem{Defn(Asashiba)}[Theo]{Definition (Asashiba)}

\newcommand{\N}{\mathbb{N}}
\newcommand{\Z}{\mathbb{Z}}

\newcommand{\A}{\mathbb{A}}

\def\La{\hbox{$\it\Lambda$}}
\newcommand{\cA}{\mathcal{A}}
\newcommand{\cB}{\mathcal{B}}
\newcommand{\tLa}{\hspace{2.5pt}\tilde{\hspace{-3pt}\La}}
\newcommand{\cF}{\mathcal{F}}

\newcommand{\rep}{{\rm rep}}
\newcommand{\Rep}{{\rm Rep}}
\newcommand{\Hom}{{\rm Hom}}

\def\Mod{\hbox{{\rm Mod}{\hskip 0.3pt}}}

\def\ModbLa{\hbox{{\rm Mod\hspace{0.6pt}}$^b$\hspace{-2pt}\La}}

\def\mf{\mathfrak}
\def\Cdot{\hspace{-1.5pt}\cdot\hspace{-2pt}}

\newcommand{\id}{{\rm 1\hspace*{-0.56ex}\rule{0.1ex}{1.52ex}\hspace*{0.2ex}}}
\def\Ga{\hbox{$\mathit\Gamma$}}
\def\Da{\hbox{${\mathit\Delta}$}}
\def\Sa{\hbox{${\mathit\Sigma}$}}
\def\Oa{\hbox{${\mathit\Omega}$}}
\def\Ta{\hbox{${\mathit\Theta}$}}

\newcommand{\dt}{{\accentset{\hspace{2pt}\mbox{\large\bfseries .}}{}}}

\newcommand{\cdt}{\dt\hspace{2pt}}

\newcommand{\pdt}{{\hspace{1pt}\dt\hspace{1.5pt}}}
\newcommand{\ydt}{{\hspace{0.5pt}\dt\hspace{1.5pt}}}

\newcommand{\mk}{\mathfrak}

\begin{document}

\title[Derived categories]{\sc The bounded derived categories of an algebra \\
with radical squared zero}

\author[Raymundo Bautista]{Raymundo Bautista \vspace{-3pt}}

\author[Shiping Liu]{Shiping Liu \vspace{-3pt}}

\keywords{Representations of a quiver; modules over a linear category; triangulated categories; derived categories; irreducible morphisms; almost split sequences; almost split triangle; Auslander-Reiten quiver; Galois covering.}

\subjclass[2010]{16E35, 16G20, 16G70, 18E30}

\thanks{The first named author is supported in part by Fondo Sectorial SEP-CONACYT via the project 43374F, while the second named author is supported in part by the Natural Sciences and Engineering Research Council of Canada.}

\address{Bautista Raymundo\\ Instituto de Matematicas, UNAM,
Unidad Morelia, Apartado Postal 61-3, 58089 Morelia, Mexico.}
\email{raymundo@matmor.unam.mx}

\address{Shiping Liu\\ D\'epartement de math\'ematiques, Universit\'e de Sherbrooke, Sherbrooke, Qu\'ebec, Canada.}
\email{shiping.liu@usherbrooke.ca}

\maketitle

\begin{abstract}

Let $\La$ be an elementary locally bounded linear category over a field with radical squared zero. We shall show that
the bounded derived category $D^b(\ModbLa)$ of finitely supported left $\La$-modules admits a Galois covering which is the bounded derived category of almost finitely co-presented representations of a gradable quiver. Restricting to the bounded derived category $D^b({\rm mod}^b\hspace{-2pt}\La)$ of finite dimensional left $\La$-modules, we shall be able to describe its indecomposable objects, obtain a complete description of the shapes of its Auslander-Reiten components, and classify those $\La$ such that $D^b({\rm mod}^b\hspace{-2.3pt}\La)$ has only finitely many Auslander-Reiten components. \vspace{-2pt}

\end{abstract}

\medskip

\section*{Introduction}

\medskip

Throughout this paper, $k$ denotes a commutative field. One of the central to\-pics in the representation theory of a finite dimensional $k$-algebra $A$ is to study its bounded derived category $D^b({\rm mod}\hspace{.4pt}A)$ of finitely generated left modules; see, for example, \cite{Ha1,Ha2}. Indeed, the triangulated category $D^b({\rm mod}\hspace{.4pt}A)$ captures all the homological properties of the algebra $A$. Since it is Hom-finite and Krull-Schmidt, we are particularly interested in classifying the indecomposable objects and studying the Auslander-Reiten theory of irreducible morphisms and almost split triangles in $D^b({\rm mod}\hspace{.4pt}A)$. These objectives have been achieved to a certain extent in the heredi\-tary case; see \cite{BLP,Ha1,Ha2}. In case $A$ is a gentle algebra, the indecomposable objects of $D^b({\rm mod}\,A)$ have been explicitly described in \cite{BeM}. Speaking of the Auslander-Reiten theory, it has been shown that an indecomposable complex of $D^b({\rm mod}\,A)$ is the ending (repsectively, starting) term of an almost split triangle if and only if it is isomorphic to a bounded complex of finitely generated projective (respectively, injective) $A$-modules; and consequently, the Auslander-Reiten quiver of $D^b({\rm mod}\hspace{.4pt}A)$ is stable if and only if $A$ is of finite global dimension; see \cite{Ha1, Ha3}. In a general setting, some particular types of stable Auslander-Reiten components of $D^b({\rm mod}\hspace{.4pt}A)$ are investigated in \cite{Sch}. In case $A$ is self-injective with no simple block, the stable Auslander-Reiten components of $D^b({\rm mod}\,A)$ are of shape $\Z \A_\infty$; see \cite[(3.7)]{Whe}, whereas the non-stable ones occur rarely and are explicitly described in \cite[(5.7)]{HKR}.

\medskip

The objective of this paper is to deal with algebras with radical squared zero. Our main tool is the covering technique, which was introduced in \cite{BoG,Gab,Gre} and further developed in \cite{As1, As2, BaL}. This requires us to work in a more general setting. Indeed, let $\La$ be a connected locally bounded $k$-category with radical squared zero. We shall assume that $\La$ is {\it elementary}, that is, all simple left $\La$-modules are one dimensional over $k$. By Gabriel's theorem, $\La\cong kQ/(kQ^+)^2$, where $Q$ is a connected locally finite quiver called the {\it ordinary quiver} of $\La$, and $kQ$ is the path category of $Q$ over $k$ with an ideal $kQ^+$ generated by the arrows; see \cite[(2.2)]{BoG}. Our aim is to study the bounded derived category $D^b(\ModbLa)$ of finitely supported left $\La$-modules and the bounded derived category $D^b({\rm mod}^b\hspace{-2.3pt}\La)$ of finite dimensional left $\La$-modules. Observe that if $\La$ is a finite dimensional $k$-algebra, then $D^b(\ModbLa)$ is the bounded derived category of all left $\La$-modules.

\medskip

In case $Q$ is gradable, using of the Koszul functor considered in \cite{BGS}, we shall obtain a triangle-equivalent
$\mathscr{F}: D^b({\rm Rep}^-(Q^{\rm op})) \to D^b(\ModbLa),$ called {\it Koszul equivalence}; see (\ref{KzEqv}), where ${\rm Rep}^-(Q^{\rm op})$ is the hereditary abelian category of almost finitely co-presented representations of the opposite quiver $Q^{\rm op}$ of $Q$; see \cite[(1.8)]{BLP}. In the general case, we shall choose a minimal gradable covering $\pi: \tilde{Q}\to Q\vspace{1pt}$, which induces a Galois covering $\pi_\lambda^D: D^b({\rm Mod}^b\hspace{-1.5pt}\tLa)\to D^b(\ModbLa)$, where $\tLa=k\tilde{Q}/(k\tilde{Q}^+)^2;$ see \cite[(7.10)]{BaL}. Composing this covering with the Koszul equivalence, we shall obtain a Galois covering $\mf{F}_\pi: D^b({\rm Rep}^-(\tilde{Q}^{\,\rm op}))\to D^b(\ModbLa);$ see (\ref{Main-1}).

\medskip

Restricting to finite dimensional modules, we shall obtain a Galois cove\-ring $\mf{F}_\pi: D^b({\rm rep}^-(\tilde{Q}^{\rm op}))\to D^b({\rm mod}^b\hspace{-2.5pt}\La)$, where ${\rm rep}^-(\tilde{Q}^{\rm op})$ is the hereditary abelian category of finitely co-presented representations of $\tilde{Q}^{\rm op}.\vspace{.5pt}$ As a consequence, the indecomposable complexes and the almost split triangles in $D^b({\rm mod}^b\hspace{-2.5pt}\La)$ can be described in terms of those in ${\rm rep}^-(\tilde{Q}^{\rm op})$; see (\ref{Main-2}). An important feature of the covering $\mf{F}_\pi$ is that it sends an indecomposable injective representation to a shift of a simple module; see (\ref{KPD-image}), called {\it simple complex}. Applying the description of the Auslander-Reiten components of $D^b({\rm rep}^-(\tilde{Q}^{\rm op}))$ obtained in \cite[Section 7]{BLP}, we shall be able to describe the Auslander-Reiten components of $D^b({\rm mod}^b\hspace{-2.5pt}\La)$ containing simple complexes in Theorem \ref{ARC-1} and those containing no simple complex in Theorem \ref{ARC-2}. In particular, if $Q$ has no infinite path, then every Auslander-Reiten component of $D^b({\rm mod}^b\hspace{-2.5pt}\La)$ is a stable tube or of shape $\Z\tilde{Q}$ or $\Z\A_\infty$; see (\ref{ARC-3}). Moreover, $D^b({\rm mod}^b\hspace{-2.5pt}\La)$ has only finitely many Auslander-Reiten components if and only if $Q$ is of Dynkin type or non-gradable of type $\tilde{\A}_n$; and in this case, its Auslander-Reiten components are explicitly described; see (\ref{ARC-4}).
Finally, we should mention that Bekkert and Drozd have described the type, as well as the indecomposable objects, of $D^b({\rm mod}^b\hspace{-2.5pt}\La)$ by using a completely different approach, that is the representation theory of bocs; see \cite{BeD}.

\smallskip

\section{Preliminaries}

\medskip

\noindent The objective of this section is to collect some basic notion and terminology and to fix some notation which will be used throughout this paper.

\vspace{2pt}

\subsection{\sc Quivers.} Let $Q=(Q_0, Q_1)$ be a locally finite quiver, where $Q_0$ is the set of vertices and $Q_1$ is the set of arrows. Given an arrow $\alpha: a\to b$, write $a=s(\alpha)$ and $b=e(\alpha)$ and introduce a {\it formal inverse} $\alpha^{-1}$ with $s(\alpha^{-1})=e(\alpha)$ and $e(\alpha^{-1})=s(\alpha)$. With each $a\in Q_0$, we associate a {\it trivial path} $\varepsilon_a$ with $s(\varepsilon_a)=e(\varepsilon_a)=a$. For $a, b\in Q_0$, we shall denote by $Q_1(a, b)$ the set of arrows from $a$ to $b$, by $Q_{\le 1}(a, b)$ the set of paths of length $\le 1$ from $a$ to $b$, and by $Q(a, b)$ the set of all paths from $a$ to $b$. Moreover, $a^+$ denotes the set of arrows in $Q$ starting at $a$, and $a^-$ is the set of arrows ending at $a$. Finally, one says that $Q$ is {\it strongly locally finite} if the number of paths from any given vertex to another one is finite; see \cite[Section 1]{BLP}.

\medskip

Given a walk $w=\alpha_r^{e_r} \cdots \alpha_2^{e_2} \alpha_1^{e_1}$ in $Q$, where $\alpha_i\in Q_1$ and $e_i=\pm 1$ such that $e(\alpha_i^{e_i})=s(\alpha_{i+1}^{e_{i+1}})$, we define its {\it degree} to be $e_1+\cdots+e_r.$ By convention, a trivial path is of degree $0$. We say that $Q$ is {\it gradable} if the walks from $a$ to $b$ have the same degree for any vertices $a, b$, or equivalently, every closed walk is of degree 0; see \cite[(7.1)]{BaL}. In general, one associates a gradable quiver $Q^{\hspace{0.5pt}\mathbb{Z}}$ with $Q$, whose vertices are the pairs $(a, i)\in Q_0\times \mathbb{Z}$, and whose arrows are the pairs $(\alpha ,i): (a, i)\to (b,i+1)$, with $\alpha: a\to b\in Q_1$ and $i\in \Z$; see \cite[(7.2)]{BaL}.
%
%
%
%
%

\medskip

Let $G$ be a group acting on $Q$. The $G$-action is called {\it free} if $g\cdot a\ne a$ for all $a\in Q_0$ and
non-identity $g\in G$. In this case, a quiver morphism $\varphi: Q\to Q'$ is called a {\it Galois $G$-covering}\hspace{.4pt}; see \cite[(4.1)]{BaL} if the following conditions are satisfied:

\begin{enumerate}[(1)]

\item The action of $\varphi$ on the vertices is surjective.

\item If $g\in G$, then $\varphi \circ g=\varphi.$

\item If $a, b\in Q_0$ with $\varphi(a)=\varphi(b)$, then $b=g\cdot a$ for some $g\in G.$

\item If $x=\varphi(a)$ with $a\in Q_0$, then $\varphi$ induces bijections $a^+\to x^+$ and $a^-\to x^-$.

\end{enumerate}


\subsection{\sc Translation quivers.} Let $(\Ga, \tau)$ be a translation quiver; see, for definition, \cite[Page 47]{Rin2}, where $\Ga$ is a quiver and $\tau$ is the translation. A vertex $a$ of $\Ga$ is called {\it projective} (respectively, {\it injective}) if $\tau a$ (respectively, $\tau^-a$) is not defined, and
{\it left stable} (respectively, {\it right stable}, {\it stable}) if $\tau^na$ is defined for all $n\ge 0$ (respectively, for all $n\le 0$,  for all $n\in \Z$). Moreover, $\Ga$ is called {\it left stable} (respectively, {\it right stable}, {\it stable}) if every vertex of $\Ga$ is left stable (respectively, right stable, stable).

\medskip

Given a quiver $\Da$ without oriented cycles, one can construct a stable translation quiver $\Z\Da$ in a canonical way; see, for example, \cite[Section 2]{Liu}. We shall denote by $\N^-\hspace{-2pt}\Da$ the full subquiver of $\mathbb{Z}\Da$ generated by the vertices $(n, x)$ with $x\in \Da_0$ and $n\le 0$; and by $\N\Da$ the one generated by the vertices $(n, x)$ with $x\in \Da_0$ and $n\ge 0$. If $\Da$ is a tree, then $\mathbb{Z}\Da$ does not depend on the orientation of $\Da$; see \cite{Ha2}. For instance, if $\Da$ is a quiver of type $\A_\infty$, then $\Z\Da$ will be written as $\Z\A_\infty$.

\medskip

A connected subquiver $\Da$ of $\Ga$ is called a {\it section} if it is convex, contains no oriented cycle, and meets every $\tau$-orbit in $\Ga$ exactly once; see \cite[(2.1)]{Liu}. In this case, $\Ga$ embeds in $\Z \Da$; see \cite[(2.3)]{Liu}. Moreover, a section $\Da$ of $\Ga$ is called {\it left-most} (respectively, {\it right-most}) if all its vertices are projective (respectively, injective); and in this case, $\Ga$ embeds in $\N\Da$ (respectively, $\N^-\hspace{-1pt}\Da$).

\medskip

Assume that $(\Ga, \tau)$ is equipped with a free action of a group $G$. Let $(\Ga', \tau')$ be another translation quiver. A {\it Galois $G$-covering} $\varphi: (\Ga, \tau) \to (\Ga', \tau')$ is a Galois $G$-covering $\varphi: \Ga\to \Da$ such, for any $a\in \Ga_0$, that $a$ is non-projective if and only if $\varphi(a)$ is non-projective; and in this case, $\tau'(\varphi(a))=\varphi(\tau a)$; compare \cite[(4.6)]{BaL}.

\vspace{2pt}

\subsection{\sc Linear categories.} Let $\cA$ be a $k$-category; see, for definition, \cite[(2.1)]{BoG}, which is equipped with an action of a group $G$. The $G$-action is called

\begin{enumerate}[$(1)$]

\item {\it free} if $g \Cdot X\not\cong X$, for any indecomposable object $X$ of $\cA$ and any non-identity element $g$ of $G\hspace{0.2pt};$

\vspace{1pt}

\item {\it locally bounded} if, for any objects $X, Y \hspace{-1pt} \in \hspace{-1pt} \cA$, $\Hom_{\mathcal{A}}(X, g\cdot Y)=0$ for all but finitely many $g\in G\hspace{0.2pt};$

\vspace{1pt}

\item {\it directed} if, for any indecomposable objects $X, Y$ of $\cA$, there exists at most one element $g$ of $G$ such that $\Hom_{\mathcal{A}}(X, \; g\cdot Y)\ne 0$ and $\Hom_{\mathcal{A}}(g\cdot Y, \; X)\ne 0.$

\end{enumerate}

\medskip

\noindent{\sc Remark.} In Definition 2.1 in \cite{BaL}, the local boundedness of the $G$-action only requires that the condition that $\Hom_{\mathcal{A}}(X, g\cdot Y)=0$ for all but finitely many $g\in G$ be satisfied for {\it indecomposable} objects $X, Y$ of $\cA$. In order for Theorem 2.12 in \cite{BaL} to hold, however, one needs to require this condition be satisfied for {\it all} objects $X, Y$ of $\cA$. Observe that all relevant statements in \cite{BaL}, such as Lemmas 6.2 and 6.6, still hold under this new definition with the same proof.

\medskip

Let $E: \cA\to \cB$ be a $k$-linear functor between $k$-categories. One says that $E$ is $G$-{\it stable} if it admits a $G$-{\it stabilizer} $\delta=(\delta_g)_{g\in G}$, where $\delta_g: E\circ g\to E$ are functorial isomorphisms such that $\delta_{h, X}\circ \delta_{g, \,h \hspace{0.3pt}\cdot X} = \delta_{gh\hskip 0.2pt, \hskip 0.5pt X},$ for $g, h\in G$ and $X\in \cA$. In this case, $E$ is called a {\it Galois $G$-precovering}; see \cite[(2.5)]{BaL} if it induces a $k$-linear isomorphism \vspace{-8pt}
$$\hspace{10pt} {\oplus}_{g\in G}\, \Hom_{\mathcal{A}}(X, g\Cdot Y) \to \Hom_{\mathcal{B}}(E(X), E(Y)): \; (u_g)_{g\in G} \mapsto {\textstyle\sum}_{g\in G}\, \delta_{g, Y}\circ E(u_g),$$ for each pair $X, Y$ of objects of $\cA$. Recall that a $G$-stabilizer $\delta=(\delta_g)_{g\in G}$ is called {\it trivial} if $\delta_{g, X}=\id_{X}$ for all $g\in G$ and $X\in \cA$; see \cite[(2.3)]{BaL}.

\medskip

We shall modify slightly the notion of a Galois covering as defined in \cite[(2.8)]{BaL}.

\medskip

\noindent{\sc Definition.} Let $\cA$ be a $k$-category with a free and locally bounded action of a group $G$. A Galois $G$-precovering $E: \cA\to \cB$ is called a {\it Galois $G$-covering} if the following conditions are satisfied.

\begin{enumerate}[$(1)$]

\item The functor $E$ is dense.

\item If $X\in \cA$ is indecomposable, then $E(X)$ is indecomposable.

\vspace{1.5pt}

\item If $X, Y\in \cA$ are indecomposable with $E(X)\cong E(Y)$, then $Y=g\Cdot X$ with $g\in G$.

\end{enumerate}

\medskip

\noindent{\sc Remark.} The definition of a Galois $G$-covering in \cite[(2.8)]{BaL} only requires $E$ be {\it almost dense}, that is, every indecomposable object of $\cB$ is isomorphic to the image under $E$ of an object of $\cA$. However, the denseness of $E$ will enables one to identify $\cB$ with the orbit category $\cA/G$ as defined in \cite[(3.1)]{Gab}. Observe that Theorem 7.10 in \cite{BaL} still holds under this new definition with the same proof; see also the proof of Theorem \ref{Main-1}.

\vspace{2pt}

\subsection{\sc Auslander-Reiten Theory.} Let $\cA$ be a Hom-finite Krull-Schmidt additive $k$-category, that is, morphism spaces are finite dimensional over $k$ and idempotent endomorphisms split; compare \cite[(1.1)]{LNP}. If $X, Y\in \cA$ are indecomposable, then a morphism $f: X\to Y$ is irreducible; see, for definition, \cite[Section 2]{AuR} if and only if it has a non-zero image in ${\rm irr}(X, Y)={\rm rad}(X, Y)/{\rm rad}^2(X, Y)$, where ${\rm rad}(X, Y)$ denotes the $k$-space of morphisms in the Jacobson radical of $\cA$; see \cite{Bau}. We shall say that $\cA$ has {\it symmetric irr-spaces} if, for any indecomposable objects $X, Y$ of $\cA$, the dimension $d_{XY}$ of ${\rm irr}(X, Y)$ over ${\rm End}(X)/{\rm rad}({\rm End}(X))$ is equal to its dimension over ${\rm End}(Y)/{\rm rad}({\rm End}(Y))$.

\medskip

A sequence $\xymatrixcolsep{16pt}\xymatrix{X\ar[r]^f & Y \ar[r]^g & Z}$ in $\cA$ is called {\it pseudo-exact} if $f$ is a pseudo-kernel of $g$, while $g$ is a pseudo-cokernel of $f$. Such a pseudo-exact sequence with $Y\ne 0$ is called {\it almost split}; see \cite[(2.1)]{Liu2}, provided that $f$ is minimal left almost and $g$ is minimal right almost; see, for definition, \cite[Section 2]{AuR}. In case $\cA$ is abelian, this notion of an almost split sequence coincides with the classical one defined in \cite[Section 2]{AuR}. In case $\cA$ is triangulated with shift functor $[1]$, a sequence $\xymatrixcolsep{16pt}\xymatrix{\hspace{-3pt}X\ar[r]^f & Y \ar[r]^g & Z}$ is almost split if and only if it embeds in an almost split triangle $\xymatrixcolsep{16pt}\xymatrix{\hspace{-2pt}X\ar[r]^f & Y \ar[r]^g & Z \ar[r]^h & X[1]}$ with $Y\ne 0$; see \cite[(6.1)]{Liu2}, which is called an {\it Auslander-Reiten triangle} in \cite[(4.1)]{Ha2}.

\medskip

Assume that $\cA$ has symmetric irr-spaces. The {\it Auslander-Reiten quiver} $\Ga_{\hspace{-2pt}\mathcal{A}}$ of $\cA$ is a translation quiver defined as follows. The vertex set is a complete set of isomorphism class representatives of the indecomposable objects of $\cA$. Given two vertices $X, Y$, the number of arrows from $X$ to $Y$ is equal to the dimension $d_{X,Y}$. The translation is the {\it Auslander-Reiten translation} $\tau_{\hspace{-1pt}_\mathcal{A}}$ which is defined so that $\tau_{\hspace{-1pt}_\mathcal{A}} Z=X$ if and only if $\cA$ has an almost split
sequence $\xymatrixcolsep{18pt}\xymatrix{X\ar[r] & Y \ar[r] & Z.}$ 
The connected components of $\Ga_{_{\mathcal{A}}}$ are called the {\it Auslander-Reiten components} of $\mathcal A$.

\vspace{2pt}

\subsection{\sc Categories of complexes.} Let $\mathcal{A}$ be a full additive subcategory of an abelian $k$-category. We shall denote by $C(\cA)$ the {\it complex category} of $\cA$, that is, the category of all complexes over $\cA$; and by $K(\cA)$ the {\it homotopy category} of $\cA$, that is, the quotient category of $C(\cA)$ modulo the null-homotopy morphisms. We shall denote by $C^{-,b}(\cA)$ and $K^{-,b}(\cA)$ the full subcategories of $C(\cA)$ and $K(\cA)$, respectively, of bounded-above complexes with bounded homology; and by $C^{b}(\cA)$ and $K^{b}(\cA)$ those of bounded complexes. Finally, $D(\cA)$ and $D^b(\cA)$ denote the {\it derived category} and the {\it bounded derived category} of $\cA$ respectively, which are localizations of $K(\cA)$ and $K^b(\cA)$, respectively, with respect to the quasi-isomorphisms. The shift functor for complexes is written as $[1]$. As usual, we shall identify an object $X$ of $\cA$ with the stalk complex $X[0]$ concentrated in degree $0$. In this way, $\cA$ becomes a full subcategory of each of $C^b(\cA)$, $K^b(\cA)$ and $D^b(\cA)$.

\vspace{2pt}

\subsection{\sc Representations of quivers.} Let $Q$ be a strongly locally finite quiver. The {\it path category} $kQ$ of $Q$ over $k$ is defined as follows; see \cite[(2.1)]{BoG}. The objects are the vertices of $Q$, and $\Hom_{kQ}(a, b)$ with $a, b\in Q_0$ is the $k$-space spanned by $Q(a,b)$. The composition of morphisms is induced from the concatenation of paths.

\medskip

A {\it $k$-representation} $M$ of $Q$ consists of a family of $k$-spaces $M(a)$ with $a\in Q_0,$ and a family of $k$-linear maps $M(\alpha): M(a)\to M(b)$ with $\alpha: a\to b\in Q_1$.
One says that $M$ is {\it locally finite dimensional} if ${\rm dim}_kM(a)$ is finite for all $a\in Q_0$; and
{\it finite dimensional} if $\sum_{a\in Q_0}{\rm dim}_kM(a)$ is finite for all $a\in Q_0$. The category ${\rm Rep}(Q)$ of all $k$-representations of $Q$ is a hereditary abelian $k$-category; see \cite{GaR}. We shall denote by ${\rm rep}^b(Q)$ the full subcategory of ${\rm Rep}(Q)$ of finite dimensional representations.

\medskip

Throughout this paper, all tensor products are over the base field $k$. With a vertex $a$ of $Q$, one associates an indecomposable injective representation $I_a$ such that $I_a(x)$, with $x\in Q_0$, is spanned by $Q(x, a)$, and $I_a(\alpha): I_a(x)\to I_a(y)$, with $\alpha: x\to y\in Q_1$, sends $\rho \alpha$ to $\rho$ and vanishes on the paths not factoring through $\alpha$. In a dual fashion, one associates with $a$ an indecomposable projective representation $P_a$ of $Q$. 
Let ${\rm Inj}\,(Q)$ stand for the full additive subcategory of ${\rm Rep}(Q)$ generated by the injective representations $I_a\otimes V_a$ with $a\in Q_0$ and $V_a$ some $k$-space; and ${\rm Proj}\,(Q)$ for the one generated by the projective representations $P_a\otimes U_a$ with $a\in Q_0$ and $U_a$ some $k$-space; see \cite[(1.3)]{BLP}. Moreover, we shall denote by ${\rm inj}\,(Q)$ the full additive subcategory of ${\rm Inj}\,(Q)$ generated by the representations $I_a$ with $a\in Q_0$, and by ${\rm proj}\,(Q)$ the full additive subcategory of ${\rm Proj}\,(Q)$ generated by the representations $P_a$ with $a\in Q_0$.

\medskip

A $k$-representation $M$ of $Q$ is called {\it almost finitely co-presented} if it admits an
injective co-resolution \vspace{-1pt}
$\xymatrixcolsep{20pt}\xymatrix{0 \ar[r] &  M \ar[r] & I_0 \ar[r] & I_1 \ar[r] & 0}\vspace{-0pt}$ with $I_0, I_1\in {\rm Inj}(Q)$, and {\it finitely co-presented} if, in addition, $I_0, I_1\in {\rm inj}(Q)$. In the dual situations, one says that $M$ is {\it almost finitely presented} and {\it finitely presented}, respectively.
The full subcategories $\Rep^-(Q)$ and $\Rep^+(Q)$ of $\Rep(Q)$ generated respectively by the almost finitely co-presented representations and by the almost finitely presented representations are hereditary abelian $k$-categories, while the full subcategories $\rep^-(Q)$ and $\rep^+(Q)$ of $\Rep(Q)$ generated by the finitely co-presented representations and by the finitely presented representations, respectively, are Hom-finite Krull-Schmidt hereditary abelian $k$-categories, whose intersection is $\rep^b(Q)$; see \cite[(1.8),(1.15)]{BLP}. Since $\rep^-(Q)$ has enough projective objects, $D^b(\rep^-(Q))$ is a Hom-finite Krull-Schmidt $k$-category; see \cite{BaS} and also \cite[(1.9)]{BaL}.

\smallskip

\section{Complexes over a locally bounded category}

\medskip

\noindent Throughout this section, $\La$ stands for a locally bounded $k$-category; see, for defi\-nition, \cite[(2.1)]{BoG}, whose object set is written as $\La_0$ and its Jacobson radical is written as ${\rm rad}\La$. A {\it left $\La$-module} is a $k$-linear functor $M: \La\to {\rm Mod}\hskip 0.4pt k$, where ${\rm Mod}\hskip 0.4pt k$ is the category of $k$-spaces; see \cite[(2.2)]{BoG}. The {\it support} of such a module $M$ is the set ${\rm supp}M$ of objects $x\in \La_0$ for which $M(x)\ne 0$. We shall say that $M$ is {\it finitely supported} if ${\rm supp}M$ is finite, and {\it finite dimensional} if $\sum_{x\in \it\Lambda_0}{\rm dim}_kM(x)$ is finite. Let $\Mod\La$ denote the category of all left $\Lambda$-modules. The full subcategory of $\Mod\La$ of finitely supported modules and that of finite dimensional ones will be denoted by $\ModbLa$ and ${\rm mod}^{\hspace{0.5pt}b}\hspace{-2.5pt}\La$, respectively. A morphism $f: M\to N$ in $\ModbLa$ is called {\it radical} if its image is contained in the radical of $N$.

\medskip

Given $x\in \La_0$, it is known that  $P[x]=\La(x, -)$ is an indecomposable projective finite dimensional $\A$-module. By Yoneda's Lemma, $\Hom_{{\rm Mod}\it\Lambda}(P[x], P[y])\cong \La(y, x)$ in such a way that a $\La$-linear morphism $f$ corresponds to $f(\id_x)\in \Lambda(y,x)$. In particular, a radical $\La$-linear morphism $P[x]\to P[y]$ corresponds to an element in $({\rm rad}\hspace{.4pt}\La)(y,x)$. Observe, for any $k$-space $V$, that $P[x]\otimes V$ is a finitely supported projective $\La$-module. We shall denote by ${\rm Proj}\hspace{0.5pt}\La$ the full additive subcategory of $\ModbLa$ generated by the modules isomorphic to some $P[x]\otimes V$ with $x\in \La_0$ and $V\in {\rm Mod}k$, and by ${\rm proj}\hspace{0.5pt}\La$ the full additive subcategory of ${\rm mod}^b\hspace{-3pt}\La$ generated by the modules isomorphic to some $P[x]$ with $x\in \La_0$.

\medskip

Every module $M$ in ${\rm Mod}^{\hspace{0.5pt}b}\hspace{-2.5pt}\La$ admits a minimal projective cover $P\to M$ with $P\in {\rm Proj} \, \La$ such that $P\in {\rm proj} \, \La$ if and only if $M\in {\rm mod}^{\hspace{0.5pt}b}\hspace{-2.5pt}\La$; see \cite [(6.1)]{BaL}. Thus, sending a bounded complex over $\ModbLa$ to its projective resolution yields a triangle-equivalence from $D^b(\ModbLa)$ to $K^{-,b}(\mbox{Proj}\,\La)$. We fix a quasi-inverse
$$\mathcal{E}_{_{\hspace{-1pt}\it\Lambda}}: K^{-,b}(\mbox{Proj}\hspace{0.5pt}\La)\to D^b(\ModbLa)$$ of this equivalence.
We shall replace $K^{-,b}(\mbox{Proj}\,\La)$ by another better behaved category. Indeed, a complex over ${\rm Proj}\hspace{0.5pt}\La$ is called {\it radical} if all its differentials are radical. We denote by $RC({\rm Proj}\hspace{0.5pt}\La)$ and $RC^{-,b}({\rm Proj}\hspace{0.5pt}\La)$ the full subcategories of $C({\rm Proj}\hspace{0.5pt}\La)$ and $C^{-,b}({\rm Proj}\hspace{0.5pt}\La)$, respectively, generated by the radical complexes. Restricting the canonical projection functor $\mathcal{P}_{_{\hspace{-1.5pt}\it\Lambda}}: C^{-,b}({\rm Proj}\hspace{0.5pt}\La)\to K^{-,b}({\rm Proj}\hspace{0.5pt}\La)$ to the radical complexes, we obtain a projection functor $\mathcal{P}_{_{\hspace{-1.5pt}\it\Lambda}}: RC^{-,b}({\rm Proj}\hspace{0.5pt}\La)\to K^{-,b}({\rm Proj}\hspace{0.5pt}\La).$

\medskip

\begin{Prop}\label{Der-Hmtp} Let $\La$ be a locally bounded $k$-category with a nilpotent radical.

\vspace{-.5pt}

\begin{enumerate}[$(1)$]

\item The projection functor $\mathcal{P}_{_{\hspace{-1.5pt}\it\Lambda}}: RC^{-,b}({\rm Proj}\hspace{0.5pt}\La)\to K^{-,b}({\rm Proj}\hspace{0.5pt}\La)$ is dense and reflects isomorphisms.

%

\item A complex $X^\cdt$ in $RC^{-, b}({\rm Proj}\hspace{0.5pt}\La)$ is indecomposable if and only if $\mathcal{P}_{_{\hspace{-1.5pt}\it\Lambda}}(X^\cdt)$ is indecomposable in $K^{-, b}({\rm Proj}\hspace{0.5pt}\La)$.

\vspace{.5pt}

\end{enumerate}\end{Prop}

\noindent{\it Proof.} Since ${\rm rad}\La$ is nilpotent, all null-homotopic endomorphisms in $RC^{-, b}({\rm Proj}\hspace{0.5pt}\La)$ are nilpotent. Let $X^\cdt$ be a complex in $RC^{-,b}({\rm Proj}\hspace{0.5pt}\La)$. Since the null-homotopic endomorphisms of $X^\cdt$ are nilpotent, every idempotent endomorphism of $\mathcal{P}_{_{\hspace{-1.5pt}\it\Lambda}}(X^\cdt)$ can be lifted to an idempotent endomorphism of $X^\cdt$ in $RC^{-, b}({\rm Proj}\hspace{0.5pt}\La)$; see \cite[(27.1)]{AnF}. As a consequence, $\mathcal{P}_{_{\hspace{-1.5pt}\it\Lambda}}(X^\cdt)$ is indecomposable in case $X^\cdt$ is indecomposable. On the other hand, if $\mathcal{P}_{_{\hspace{-1.5pt}\it\Lambda}}(X^\cdt)=0$, that is, $\id_{X^\ydt}$ is null-homotopic, then $\id_{X^\ydt}$ is nilpotent. Thus, $\id_{X^\ydt}=0$, that is, $X^\cdt=0$ in $RC^{-, b}({\rm Proj}\hspace{0.5pt}\La)$. As a consequence, $\mathcal{P}_{_{\hspace{-1.5pt}\it\Lambda}}(X^\cdt)$ is decomposable in case $X^\cdt$ is decomposable. This establishes Statement (2).

Let $f^\ydt: X^\cdt\to Y^\ydt$ be a morphism in $RC^{-,b}({\rm Proj}\hspace{0.5pt}\La)$ such that $\mathcal{P}_{_{\hspace{-1.5pt}\it\Lambda}}(f^\ydt)$ is an isomorphism. Then, $\mathcal{P}_{_{\hspace{-1.5pt}\it\Lambda}}(f^\ydt) \, \mathcal{P}_{_{\hspace{-1.5pt}\it\Lambda}}(g^\cdt)= \mathcal{P}_{_{\hspace{-1.5pt}\it\Lambda}}(g^\cdt) \, \mathcal{P}_{_{\hspace{-1.5pt}\it\Lambda}}(f^\ydt) =\id_{\,\mathcal{P}_{_{\hspace{-1.5pt}\it\Lambda}}(X^\cdt)}\vspace{1pt}$, for some morphism $g^\cdt: Y^\ydt\to X^\cdt$ in $RC^{-,b}({\rm Proj}\hspace{0.5pt}\La)$. That is, $f^\ydt g^\cdt -\id_{X^\cdt}$ and $g^\cdt f^\ydt  -\id_{X^\cdt}$ are null-homotopic, and hence, they are nilpotent. As a consequence, $f^\ydt g^\cdt$ and $g^\cdt f^\ydt$ are automorphisms of $X^\cdt$. Therefore, $f^\ydt$ is an automorphism. This shows that $\mathcal{P}_{_{\hspace{-1.5pt}\it\Lambda}}$ reflects isomorphisms.

\vspace{1pt}

It remains to prove that $\mathcal{P}_{_{\hspace{-1.5pt}\it\Lambda}}$ is dense. Consider a complex $(X^\cdt, d_X^\pdt)\vspace{.5pt}\in K^{-,b}({\rm Proj}\hspace{0.5pt}\La)$. We may assume that ${\rm H}^n(X^\cdt)=0$ for $n\le 0$. Write $d_X^{\hspace{1pt}0}=pj$, where $p: X^0\to L$ is an epimorphism and $j: L\to X^1$ is an monomorphism. Let\vspace{-2pt}
$$\xymatrix{
\cdots \ar[r] & P^n \ar[r]^{d^n}& P^{n+1}\ar[r] & \cdots \ar[r] & P^{-1} \ar[r]^-{d^{-1}} & P\hspace{1pt}^0 \ar[r]^-{d\hspace{1pt}^0} &
L\ar[r] & 0
}$$ be a minimal projective resolution of $L$ in $\ModbLa$, where the $d^n$ with $n<0$ are radical.
Set $Z^n=P^n$ for $n\le 0$ and $Z^n=X^n$ for $n>0$. Put $d_Z^n=d^n$ for $n<0$, and $d_Z^{\hspace{1pt}0}=d^{\hspace{1pt}0}j$, and $d_Z^n=d_X^n$ for $n>0$. Then, $Z^\ydt\in C^{-,b}({\rm Proj}\hspace{0.5pt}\La)$ such that $X^\cdt\cong Z^\ydt$ in $K^{-,b}({\rm Proj}\hspace{0.5pt}\La)$. Therefore, we may assume that $d_X^n$ is radical for every $n<0$.

\vspace{1pt}

Suppose that $d_X^{\hspace{0.5pt}s}$ is not radical for some $s\ge 0$. Since $X^\cdt$ is bounded-above, we may assume that $d_X^{\hspace{0.5pt}n}$ is radical for all $n>s$. Since $X^{s+1}$ is projective, we may write 
$$
d_X^{\hspace{0.5pt}s}=\left(\begin{array}{cc} \id_M & h \\ 0 & g^s \end{array}\right): X^s=M\oplus N^s\longrightarrow M\oplus N^{s+1}=X^{s+1},
$$
where $g^s$ is a radical morphism. Since $d_X^{s+1}$ is radical such that $d_X^{s+1} d^s_X=0$, we see that $d_X^{s+1}=(0, g^{s+1})$, where $g^{s+1}$ is radical such that $g^{s+1}g^s=0$. Writing $d_X^{s-1}=(f^{s-1}, g^{s-1})^T$,
we obtain $f^{s-1}+h g^{s-1}=0$ and $g^sg^{s-1}=0.$

Set $Y^n=X^n$ for $n\not\in \{s, s+1\}$ and $Y^n=N^n$ for $n\in \{s, s+1\}$. Put $d_Y^n=d_X^n$ for $n\not\in \{s-1, s, s+1\}$ and $d_Y^{\hspace{0.5pt}n}=g^n$ if $n\in \{s-1, s, s+1\}.$ This yields a complex $Y^\ydt$ in $C^{-,b}({\rm Proj}\hspace{0.5pt}\La).$ We claim that $X^\cdt\cong Y^\ydt$ in $K^{-,b}({\rm Proj}\hspace{0.5pt}\La).$ Indeed,
we define $\psi^n=\varphi^n=\id_{X^n}$ for $n\notin \{s, s+1\}$ and $\psi^n=(0, \id_{N^n}): M\oplus N^n\to N^n$
for $n\in \{s, s+1\}$. Moreover, we set

$$\varphi^s=\left(\begin{array}{l} -h \\ \id_{N^s}\end{array}\right): N^s \to M\oplus N^s; \quad
\varphi^{s+1}=\left(\begin{array}{c} 0\\ \id_{N^{s+1}} \end{array}\right): N^{s+1}\to M\oplus N^{s+1}.$$

Then $\psi^\pdt=(\psi^n)_{n\in \mathbb{Z}}: X^\cdt \to Y^\pdt$ and $\varphi^\pdt=(\varphi^n)_{n\in \mathbb{Z}}: Y^\pdt\to  X^\cdt$ are morphisms such that $\psi^\pdt \phi^\pdt =\id_{Y^\ydt}$ and $\phi^\pdt \psi^\pdt $ is homotopic to $\id_{X^\cdt}$. This establishes our claim. By definition, $d_Y^n$ is radical for all $n>s-1$. If $s=0$ then, since $d_X^{-1}$ is radical, $g^{-1}$ is radical; and in this case, $Y^\pdt\in RC^{-, b}({\rm Proj}\hspace{0.5pt}\La)$. Otherwise, $d_Y^n$ is radical for all $n<0$. By induction, we may find a complex $M^\pdt\in RC^{-, b}({\rm Proj}\hspace{0.5pt}\La)$ such that $X^\cdt\cong M^\pdt$ in $K^{-, b}({\rm Proj}\hspace{0.5pt}\La).$ The proof of the proposition is completed.

\medskip

For the rest of this section, suppose that $\La$ is elementary (that is, all simple $\La$-modules are one-dimensional over $k$) with radical squared zero. By Gabriel's Theorem, we may assume that $\La=kQ/(kQ^+)^2,\vspace{.5pt}$ where $Q$ is a locally finite quiver and $kQ^+$ is the ideal in $kQ$ generated by the arrows; see \cite[(2.1)]{BoG}.
%
%
For $u\in kQ$, we shall write $\bar{u}=u+(kQ^+)^2\in \La$. For an arrow $\alpha: x\to y$ in $Q$, let $P[\alpha]: P[y]\to P[x]$ be the $\La$-linear morphism given by the left multiplication by $\bar{\alpha};$ and for a trivial path $\varepsilon_x$ with $x\in Q_0$, we shall write $P[\varepsilon_x]=\id_{P[x]}$. We quote from \cite[(7.6)]{BaL} the following description of the morphisms in ${\rm Proj}\hspace{0.5pt}\La$.

\medskip

\begin{Lemma}\label{rqz-pm}

Let $\La=kQ/(kQ^+)^2,\vspace{1pt}$ where $Q$ is a locally finite quiver. If $x, y\in Q_0$ and $U, V$ are $k$-vector spaces, then every $\La$-linear morphism $f: P[x]\otimes U\to P[y]\otimes V$ can be uniquely written as
$f=\textstyle{\sum}_{\,\gamma\in Q_{\le 1}(y, x)}\, P[\gamma]\otimes f_\gamma,\vspace{1pt}$ where $f_\gamma\in {\rm Hom}_k(U, V),$
such that $f$ is radical if and only if $f={\sum}_{\alpha\in Q_1(y, x)}\, P[\alpha]\otimes f_\alpha$, where $f_\alpha\in {\rm Hom}_k(U, V).$

\end{Lemma}

\medskip

We shall give a decomposition for each complex $M^\ydt\in RC({\rm Proj}\hspace{0.5pt}\La).$ Indeed, for each $n\in \Z$, we can write $M^n=\oplus_{x\in Q_0}\, P[x]\otimes M_x^n,$ where $M_x^n \in \mbox{Mod}\,k,$ and $d^n_M=(d^n_M(y, x))_{(y, x)\in Q_0\times Q_0},\vspace{1pt}$ where $d^n_M(y, x): P[x]\otimes M_x^n \to P[y]\otimes M_y^{n+1}$ is a $\La$-linear morphism. In view of Lemma \ref{rqz-pm}, we obtain
$$d^n_M(y, x)={\sum}_{\alpha\in Q_1(y,x)}\, P[\alpha]\otimes M_\alpha^n\,, \; \mbox{ where }\;  M_\alpha^n \in {\rm Hom}\,_k(M_x^n, M_y^{n+1}).$$

We shall need to consider the opposite quiver $Q^{\rm \hspace{.5pt}op}$ of $Q$ with $(Q^{\rm \hspace{.5pt}op})_0=Q_0$ and $(Q^{\rm \hspace{.5pt}op})_1=\{ \alpha^{\rm o}: y\to x \mid \alpha: x\to y\in Q_1\}.$ We define the {\it support-quiver} $\Oa(M^\ydt)$ of $M^\ydt$ to be a subquiver of the gradable quiver $(Q^{\rm \hspace{.5pt}op})^{\Z},$ whose vertices are $(x, n)\in Q_0\times \Z$ with $M_x^n\ne 0;$ and whose arrows are $(\alpha^{\rm o}, n): (x, n)\to (y, n+1)$, where $\alpha: y\to x \in Q_1$ and $n\in \Z$ such that $M_\alpha^n\ne 0.$ Given an integer $n$, since $M^n$ is a finite direct sum of modules of form $P[x]\otimes U$, the set $\Oa^{\hspace{.5pt}n}\hspace{-0.5pt}(M^\pdt)$ of vertices $(x, n)\in \Oa(M^\pdt)$ with $x\in Q_0$ is finite such that
$M^n=\oplus_{(x, n)\in {\it\Omega}^{\hspace{.5pt}n}\hspace{-0.5pt}(M^\pdt)}\, P[x]\otimes M_x^n.$

\medskip

Given a connected component $\mathcal{C}$ of $\Oa(M^\cdt)$, the {\it restriction} $M^\pdt_\mathcal{C}$ of $M^\ydt$ to $\mathcal{C}$ is a radical complex over ${\rm Proj}\hspace{0.5pt}\La$ defined as follows. For each integer $n$, we write $\mathcal{C}^n=\mathcal{C}\cap \Oa^{\hspace{.5pt}n}\hspace{-0.5pt}(M^\pdt)$, and we put
$M^n_\mathcal{C}=\oplus_{(x,n)\in \mathcal{C}^n} P[x]\otimes M_x^n,$ a summand of $M^n.$ Moreover, we define
$d^n_{M_\mathcal{C}}$ to be the composite of the following morphisms
$$\xymatrix{M^n_\mathcal{C} \ar[r]^{q_{_\mathcal{C}}^n} &  M^n \ar[r]^{d_M^n} & M^{n+1} \ar[r]^{p_{_\mathcal{C}}^{n+1}} & M^{n+1}_\mathcal{C},
}$$ where $q_{_\mathcal{C}}^n$ is the canonical injection and $p_{_\mathcal{C}}^{n+1}$ is the canonical projection.

\medskip

\begin{Prop}\label{RC-decom}

Let $\La=kQ/(kQ^+)^2$ with $Q$ a locally finite quiver, and let $M^\ydt$ be a non-zero complex in $RC({\rm Proj}\hspace{.5pt}\La)$. If $\mathcal{C}_i$ with $i\in I$ are the connected components of the support-quiver of $M^\ydt$, then
$M^\ydt=\oplus_{i\in I}\, M_{\mathcal{C}_i}^\ydt.$

\end{Prop}

\noindent{\it Proof.} We shall keep the notation introduced above. Let $\mathcal{C}_i$ with $i\in I$ be the connected components of $\Oa(M^\ydt)$. For each $n\in \Z$, since
$\Oa^{\hspace{.5pt}n}\hspace{-0.5pt}(M^\pdt)$ is a disjoint union of $\mathcal{C}_i^n$ with $i\in I$, we obtain
$$M^n=\oplus_{(x, n)\in {\it\Omega}^{\hspace{.5pt}n}\hspace{-0.5pt}(M^\pdt)}\, P[x]\otimes M_x^n=\oplus_{i\in I}\left(\oplus_{(x, n)\in \mathcal{C}_i^n}\, P[x]\otimes M_x^n\right)=\oplus_{i\in I}\, M_{\mathcal{C}_i}^n.$$

Consider a pair $(i, j)\in I\times I$ with $i\ne j$. Let $(x, n)\in \mathcal{C}_i^n$ and $(y, n+1)\in \mathcal{C}_j^{n+1}$. Since $\Oa(M^\cdt)$ has no arrow $(x, n)\to (y, n+1)$, we have $M_\alpha^n=0$ for all possible arrows $\alpha: y\to x$ in $Q$. Thus, $d^n_M(y, x)={\sum}_{\alpha\in Q_1(y,x)}\, P[\alpha]\otimes M_\alpha^n=0$. As a consequence, the sequence \vspace{-2pt}
$$\xymatrix{M^n_{\mathcal{C}_i} \ar[r]^{q_{_{\mathcal{C}_i}}^n} &  M^n \ar[r]^{d_M^n} & M^{n+1} \ar[r]^{p_{_{\mathcal{C}_j}}^{n+1}} & M^{n+1}_{\mathcal{C}_j},
}$$ has null composite. This shows that $d_{M}^n=\oplus_{i\in I}\,d_{M_{\mathcal{C}_i}}^n.$ The proof of the proposition is completed.

\medskip

A radical complex over ${\rm Proj}\hspace{.5pt}\La$ is said to be {\it support-connected} if its support-quiver is connected. In view of Proposition \ref{RC-decom}, we see that an indecomposable radical complex over ${\rm Proj}\hspace{.5pt}\La$ is support-connected.

\medskip

\begin{Cor}\label{RC-fin-decom}

Let $\La=kQ/(kQ^+)^2,\vspace{1pt}$ where $Q$ is a locally finite quiver. If $M^\ydt$ is a non-zero complex in $RC^{-, b}({\rm Proj}\hspace{.5pt}\La)$, then $M^\ydt=M_1^\ydt\oplus \cdots \oplus M_m^\ydt,$
where the $M_i^\ydt$ are support-connected complexes in $RC^{-, b}({\rm Proj}\hspace{.5pt}\La).$

\end{Cor}

\noindent{\it Proof.} Let $M^\ydt$ be a non-zero complex in $RC^{-, b}({\rm Proj}\hspace{.5pt}\La)$. There exist integers $s, t$ with $s\le t$ such that $M^i=0$ for all $i>t$ and $H^i(M^\ydt)\ne 0$ for all $i<s.$ Let $\mathcal{C}_i$ with $i\in I$ be the connected components of $\Oa(M^\ydt)$. Letting $M_i^\ydt$ be the restriction of $M^\ydt$ to
$\mathcal{C}_i,$ by Proposition \ref{RC-decom}, we obtain $M^\ydt=\oplus_{i\in I}\, M_i^\ydt\hspace{.4pt}.$

For each $i\in I$, there exists a maximal integer $n_i\le t$ such that $M_i^{n_i}\ne 0.$ Thus there exists some $x_i\in Q_0$ such that $(x_i, n_i)\in \mathcal{C}_i$. In particular, $(x_i, n_i)\in \Oa^{n_i}(M^\ydt)$ with $n_i\le t$. Since $M_i^\ydt$ is a radical complex, by the maximality of $n_i$, we see that ${\rm H}^{n_i}(M_i^\ydt)\ne 0$, and consequently, $n_i\ge s$. Since each $\Oa^n(M^\ydt)$ with $n\in \Z$ is finite, the set $\{(x_i, n_i) \mid i\in I\}$ is finite. Since the $\mathcal{C}_i$ with $i\in I$ are pairwise disjoint, we see that $I$ is finite. The proof of the corollary is completed.

\smallskip

\section{\sc Koszul equivalence in the graded case}

\medskip

\noindent Throughout this section, let $\La=kQ/(kQ^+)^2,$ where $Q$ is a gradable locally finite quiver. Our objective is to show that the bounded derived categories $D^b(\ModbLa)$ and $D^b({\rm mod}\hspace{.4pt}^b\hspace{-2pt}\La)$ are equivalent to the bounded derive categories $D^b({\rm Rep}^-(Q^{\rm op}))$ and $D^b({\rm rep}^-(Q^{\rm op}))$, respectively. These equivalences are given by the Koszul dua\-lity considered in \cite{BGS}. Note, however, that our results in this section cannot be derived directly from theirs.

\medskip

We start with a grading $Q_0=\cup_{n\in \Z}\, Q^n$ for $Q$, where the $Q^n$ are pairwise disjoint such that every arrow of $Q$ is of form $x\to y$ with $x\in Q^n$ and $y\in Q^{n+1}$ for some $n;$ see the remark following \cite[(7.1)]{BaL}. In particular, if $x\in Q^m$ and $y\in Q^n$, then every possible path from $x$ to $y$ is of length $n-m$. Being locally finite, $Q$ is strongly locally finite, and so is $Q^{\,\rm op}$. More importantly, $Q^{\rm \hspace{.5pt}op}$ is also gradable with a grading $(Q^{\rm \hspace{.5pt}op})^n=Q^{-n}$, for $n\in \Z$.

\medskip

A $k$-representation $M$ of $Q^{\rm \hspace{.5pt}op}$ is called {\it locally support-finite} if $(Q^{\rm \hspace{.5pt}op})^n \cap {\rm supp} M$ is finite, for all $n\in \Z$. We denote by ${\rm Rep}^*(Q^{\, \rm op})$ the full subcategory of ${\rm Rep}(Q^{\rm \hspace{.5pt}op})$ generated by the locally support-finite representations. With each representation $M\in {\rm Rep}^*(Q^{\, \rm op})$, we shall associate a radical complex over $\mbox{Proj} \hspace{0.5pt}\La$ as follows. Indeed, for each $n\in \mathbb{Z}$, we set
$$F(M)^n=\oplus_{x\in Q^{-n}}\, P[x]\otimes M(x)\in {\rm Proj}\,\La,$$
and define $
d_{F(M)}^n: \oplus_{x\in Q^{-n}} \, P[x]\otimes M(x) \to \oplus_{y\in Q^{-n-1}}\, P[y]\otimes M(y)
\vspace{2pt} $
to be the $\La$-linear morphism given by the matrix$
(d_{F(M)}^n(y,x))_{(y,x) \in Q^{-n-1}\times Q^{-n}},
$
where
$$
d^n_{F(M)}(y,x)={\sum}_{\alpha\in Q_1(y,x)}\,P[\alpha]\otimes M(\alpha^{\rm o}): P[x]\otimes M(x)\to P[y]\otimes M(y).
$$

Since ${\rm rad}^{\hspace{0.4pt}2}\hspace{-2.4pt}\La=0$, this yields a radical complex $F(M)^\cdt$ over $\mbox{Proj} \hspace{0.5pt}\La$. Let now $f: M\to N$ be a morphism in ${\rm Rep}^*(Q^{\rm \hspace{.5pt}op})$. For $n\in \Z$, we define
$$
F(f)^n=\oplus_{x\in Q^{-n}}\, \id_{P[x]}\otimes f(x)\,: \oplus_{x\in Q^{-n}} \, P[x]\otimes M(x) \to \oplus_{x\in Q^{-n}} \, P[x]\otimes N(x).
$$

It is easy to verify that $F(f)^\cdt=(F(f)^n)_{n\in \mathbb{Z}}$ is a morphism from $F(M)^\cdt$ to $F(N)^\cdt$. In this way, we obtain a $k$-linear functor $F: {\rm Rep}^*(Q^{\rm \hspace{.5pt}op})\to RC({\rm Proj}\hspace{0.5pt}\La),$
called {\it Koszul functor}.

\medskip

\begin{Lemma}\label{F-property}

Let $\La=kQ/(kQ^+)^2,$ where $Q$ is a gradable locally finite quiver.
The Koszul functor $F: {\rm Rep}^*(Q^{\rm \hspace{.5pt}op})\to RC({\rm Proj}\hspace{0.6pt}\La)$ is full, faithful and exact.

\end{Lemma}

\noindent{\it Proof.} First of all, in view of the property of the tensor product over $k$, we see easily that $F$ is exact and faithful. Consider a morphism $f^\ydt: F(M)^\cdt\to F(N)^\cdt$ in $RC({\rm Proj}\hspace{0.6pt}\La)$ with
$M, N\in {\rm Rep}^*(Q^{\rm \hspace{.5pt}op}).$ For $n\in \Z$, write
$$f^n=(f^n(y,x))_{(y,x)\in Q^{-n}\times Q^{-n}}: \oplus_{x\in Q^{-n}}\, P[x]\otimes M(x)\to \oplus_{y\in Q^{-n}}\, P[y]\otimes N(y),$$ where $f^n(y,x): P[x]\otimes M(x)\to P[y]\otimes N(y)$ is $\La$-linear and can be written as
$f^n(y,x)=\textstyle{\sum}_{\gamma\in Q_{\le 1}(y, x)}\, P[\gamma]\otimes f_\gamma,$ with $f_\gamma\in \Hom_k(M(x), N(y));$ see (\ref{rqz-pm}). Since $x, y\in Q^{-n}$, we see that $Q_{\le 1}(y, x)\ne \emptyset$ if and only if $y= x$. Thus, $f^n(y,x)=0$ if $y\ne x,$ and $f^n(x, x)=\id_{P[x]}\otimes f_{\varepsilon_x}$, and consequently, we obtain
$$f^n={\rm diag}\{\id_{P[x]}\otimes f_{\varepsilon_x}\}_{x\in Q^{-n}}: \oplus_{x\in Q^{-n}}\, P[x]\otimes M(x)\to \oplus_{x\in Q^{-n}}\, P[x]\otimes N(x).$$

Now, let $(y, x)\in Q^{-n-1}\times Q^{-n}\vspace{1pt}$ with $Q_1(y,x)\ne \emptyset$. In view of the equation $f^{n+1} \circ d_{F(M)}^n=d_{F(N)}^n \circ f^n$, we see that
$$\textstyle{\sum}_{\alpha\in Q_1(y, x)}\, P[\alpha]\otimes (f_{\varepsilon_y}\circ M(\alpha^{\rm o}))
=\textstyle{\sum}_{\alpha\in Q_1(y, x)}\, P[\alpha]\otimes (N(\alpha^{\rm o}) \circ f_{\varepsilon_x}).$$
By the uniqueness stated in Lemma \ref{rqz-pm}, $f_{\varepsilon_y}\circ M(\alpha^{\rm o})=N(\alpha^{\rm o}) \circ f_{\varepsilon_x}$, for every arrow $\alpha: y\to x$ in $Q$. This shows that $g=(f_{\varepsilon_x})_{x\in Q_0}: M\to N$ is a morphism such that $F(g)=f$. The proof of the lemma is completed.

\medskip

We shall describe the image of the Koszul functor. The following statement describes the objects up to shift, which generalizes slightly Proposition 7.7 in \cite{BaL}.

\medskip

\begin{Prop}\label{Ind-complex}

Let $\La=kQ/(kQ^+)^2$ with $Q$ a gradable locally finite quiver, and let $F: {\rm Rep}^*(Q^{\rm \hspace{.5pt}op})\to RC({\rm Proj}\hspace{0.6pt}\La)$ be the Koszul functor. If $M^\ydt\in RC({\rm Proj}\hspace{.6pt}\La)\vspace{1pt}$ is support-connected, then $M^\ydt\cong F(N)^\cdt[s]$ for some $N\in {\rm Rep}^*(Q^{\rm \hspace{.5pt}op})$ and $s\in \Z$.

\end{Prop}

\noindent{\it Proof.} Let $M^\ydt$ be a non-zero complex in $RC({\rm Proj}\hspace{.6pt}\La)\vspace{1pt}$ with a connected support-quiver $\Oa(M^\ydt)$. For $n\in \Z$, write $M^n=\oplus_{x\in Q_0}\, P[x]\otimes M_x^n$ with $M_x^n \in \mbox{Mod}\,k,$ and $d^n_M=(d^n_M(y, x))_{(y, x)\in Q_0\times Q_0},\vspace{2pt}$ where
$$d^n_M(y, x)={\sum}_{\alpha\in Q_1(y,x)}\, P[\alpha]\otimes M_\alpha^n\,, \; \mbox{ with }\;  M_\alpha^n \in {\rm Hom}\,_k(M_x^n, M_y^{n+1}).$$

We claim that there exists some integer $s$ such that for all $n\in \Z$, we have
$$M^n=\oplus_{x\in Q^{-n-s}}\, P[x]\otimes M_x^n,$$ where $M_x^n=0$ for all but finitely many $x\in  Q^{-n-s}.$ Indeed, let $M^{n_0}$ be non-zero for some $n_0$. Then $M_a^{n_0}\ne 0$ for some $a\in Q_0$, that is, $(a, n_0)\in \Oa(M^\ydt).$ Assume that $a\in Q^t$ for some integer $t$. Set $s=-t-n_0$. Let $n$ be an integer such that $M^n$ is non-zero. Then, $M_x^n\ne 0$ for some
$x\in Q_0$, and hence, $(x, n)\in \Oa(M^\ydt).$ Being connected, $\Oa(M^\pdt)$ contains a walk from $(a, n_0)$ to $(x, n)$, which is a walk in $(Q^{\rm \hspace{.5pt}op})^{\mathbb{Z}}$ of degree $n-n_0$. As a consequence, $Q^{\rm \hspace{.5pt}op}$ contains a walk from $a$ to $x$ of degree $n-n_0$; see \cite[(7.2)]{BaL}, that is, $Q$ contains a walk from $a$ to $x$ of degree $n_0-n$. As a consequence, $x\in Q^{t+n_0-n}=Q^{-n-s}$. The establishes our claim. As a consequence, for each $n\in \Z$, we obtain $d^n_M=(d^n_M(y, x))_{(y, x)\in Q^{-s-n-1}\times Q^{-s-n}},\vspace{1pt}$ where
$d^n_M(y, x)={\sum}_{\alpha\in Q_1(y,x)}\, P[\alpha]\otimes M_\alpha^n,$ with $M_\alpha^n \in {\rm Hom}\,_k(M_x^n, M_y^{n+1}).\vspace{1.5pt}$

Now, we shall define an object $N\in {\rm Rep}^*(Q^{\rm \hspace{.5pt}op}).$ For a vertex $x\in Q^n$, we set $N(x)=M^{-n-s}_x$; and for an arrow $\alpha: y\to x$ with $(y, x)\in Q^{n-1}\times Q^n$, set
$N(\alpha^{\rm o})=(-1)^s M_\alpha^{-n-s}$. For each $n\in \Z$, we have
$$F(N)^n[s]=\oplus_{x\in Q^{-n-s}}\, P[x]\otimes N(x) =\oplus_{x\in Q^{-n-s}}\, P[x]\otimes M_x^{(n+s)-s}=M^n$$
and $d_{F(N)[s]}^{\hspace{.5pt}n}=(-1)^sd_{F(N)}^{\hspace{.5pt}n+s}
=\,((-1)^sd_{F(N)}^{\hspace{.5pt}n+s}(y,x))_{(y,x) \in Q^{-n-s-1} \times Q^{-n-s}},
\vspace{2pt}$
where
$$
\begin{array}{rcl}
(-1)^sd_{F(N)}^{\hspace{.5pt}n+s}(y,x)
& = & (-1)^s{\sum}_{\alpha\in Q_1(y,x)}\,P[\alpha]\otimes N(\alpha^{\rm o})\\
\vspace{-9pt}\\
& = & {\sum}_{\alpha\in Q_1(y,x)}\,P[\alpha]\otimes M_\alpha^{(n+s)-s}
\\ \vspace{-9pt}\\
& = & d_M^{\hspace{.5pt}n}(y,x). \vspace{3pt}
\end{array}$$
This shows that $d_{F(N)[s]}^{\hspace{.5pt}n} = d_M^{\hspace{.5pt}n}\hspace{.5pt}.\vspace{1pt}$ The proof of the proposition is completed.

\medskip

In order to describe morphisms in the image of the Koszul functor, we shall say that a morphism $f^\ydt: M^\ydt\to N^\ydt$ in $RC({\rm Proj}\hspace{.6pt}\La)$ is {\it radical} if $f^n: M^n\to N^n$ is a radical morphism for every $n\in \Z$.

\medskip

\begin{Lemma}\label{Morph-F-images}

Let $\La=kQ/(kQ^+)^2$ with $Q$ a gradable locally finite quiver, and let $F: {\rm Rep}^*(Q^{\rm \hspace{.5pt}op})\to RC({\rm Proj}\hspace{0.6pt}\La)$ be the Koszul functor. Consider a non-zero mor\-phism $f^\pdt: F(M)^\cdt\to F(N)^\cdt[s]$ in $RC({\rm Proj}\hspace{0.5pt}\La)$ with $M, N\in {\rm Rep}^*(Q^{\rm \hspace{.5pt}op})$ and $s\in \Z$.

\begin{enumerate}[$(1)$]

\item If $f^\pdt$ is non-radical, then $s=0$.

\vspace{.5pt}

\item If $f^\pdt$ is radical, then $s=1.$

\end{enumerate} \end{Lemma}

\noindent{\it Proof.} For $n\in \Z$, we write the $\La$-linear morphism $f^n: F(M)^n\to F(N)^{n+s}$ as
$$f^n=(f^n(y, x))_{(y,x)\in Q^{-n-s} \times Q^{-n}}: \oplus_{x\in Q^{-n}}\, P[x]\otimes M(x) \to \oplus_{y\in Q^{-n-s}}\, P[y]\otimes M(y),$$
where $f^n(y,x): P[x]\otimes M(x) \to P[y]\otimes N(y)$
is a $\La$-linear morphism. By Lemma \ref{rqz-pm}, we have
$f^n(y,x)={\textstyle\sum}_{\gamma\in Q_{\le 1}(y, x)}\, P[\gamma]\otimes f_\gamma,$ with $f_\gamma\in \Hom_k(M(x), N(y)).$

\vspace{1.5pt}

Suppose that $f^m\ne 0$ for some integer $m$. Then, $f^m(y,x)\ne 0$ for some pair $(y, x)\in Q^{-m-s} \times Q^{-m}.$ In case $f^\cdt$ is non-radical, we may assume that $f^m(y,x)$ is non-radical. Then, $Q_{\le 1}(y, x)$ contains a trivial path, that is, $y=x$. As a consequence, $s=0$. If $f^\cdt$ is radical, then $f^m(y,x)$ is radical, and thus, $Q_1(y, x)\ne \emptyset$. As a consequence, $s=1$. The proof of the lemma is completed.

\medskip

Let $x\in Q_0$. We denote by $I_{x^{\rm o}}$ the indecomposable injective representation of $Q^{\hspace{.5pt}\rm op}$ associated with $x$. Since $Q^{\rm op}$ is gradable locally finite, $I_{x^{\rm o}}$ is locally support-finite, and so is every representation in ${\rm Rep}^{-}(Q^{\hspace{.5pt}\rm op})$. 
Moreover, we denote by $S[x]$ the simple $\La$-module supported by $x$, and by $P^\pdt_{S[x]}$ its deleted minimal projective resolution in ${\rm mod}\La$, which is an object of $RC^{-,b}(\mbox{Proj}\La)$ such that $P_{S[x]}^n=0$ for $n>0$, and $P_{S[x]}^0=P[x]$, and ${\rm H}^n(P_{S[x]}^\pdt)=0$ for $n<0$.

\medskip

\begin{Lemma}\label{F-inj-image}

Let $\La=kQ/(kQ^+)^2$ with $Q$ a gradable locally finite quiver, and let $F: {\rm Rep}^*(Q^{\rm \hspace{.5pt}op})\to RC({\rm Proj}\hspace{0.6pt}\La)$ be the Koszul functor. If $a\in Q^t$ for some $t\in \Z$, then $F(I_{a^{\rm o}})^\cdt\cong P^\pdt_{S[a]}[t].$

\end{Lemma}

\noindent{\it Proof.} We shall regard $\La$ as a $k$-algebra and identify a module $M\in \ModbLa$ with $\oplus_{x\in Q_0}M(x)$. Let $a\in Q^t$. For $n\in \Z$, we have $F(I_{a^{\rm o}})^{\hspace{0.4pt}n}=\oplus_{x\in Q^{-n}}\, P[x] \otimes I_{a^{\rm o}}(x)$ and
$d_{F(I_{a^{\rm o}})}^{\hspace{0.4pt}n}=(d_{F(I_{a^{\rm o}})}^{\hspace{0.4pt}n}(y, x))_{(y, x)\in Q^{-n-1}\times Q^{-n}},$
where $$d_{F(I_{a^{\rm o}})}^{\hspace{0.4pt}n}(y, x)={\sum}_{\alpha\in Q_1(y, x)}\, P[\alpha]\otimes I_{a^{\rm o}}(\alpha^{\rm o}).\vspace{0pt}$$

Let $a\in Q^t$. Write $s=-t$. Every path in $Q^{\rm \hspace{.5pt}op}$ ending at $a$ is the opposite path $p^{\rm o}$ of a path $p$ in $Q$ starting at $a$, and hence, $s(p^{\rm o})\in Q^{-n}$ for some $n\le s$. In particular, $F(I_{a^{\rm o}})^n=0$ for all $n>s$. For each $n\le s$, the projective $\La$-module $F(I_{a^{\rm o}})^n=\oplus_{x\in Q^{-n}}\, P[x] \otimes I_{a^{\rm o}}(x)$ has a $k$-basis
$$\mathit\Omega_n=\{\bar{\varepsilon}_x \otimes p^{\rm o}, \bar{\beta}\otimes p^{\rm o} \mid x\in Q^{-n}, \beta\in Q_1(x, -), p\in Q(a, x)\},$$ where $Q_1(x, -)$ denotes the set of arrows of $Q$ starting at $x,$ while ${\rm rad}\hspace{0.4pt} F(I_{a^{\rm o}})^n$ has a $k$-basis $\mathit\Theta_n=\{\bar{\beta}\otimes p^{\rm o} \mid x\in Q^{-n}, \beta\in Q_1(x, -), p\in Q(a, x)\}.$ In particular, $F(I_{a^{\rm o}})^s\cong P[a]$.

\vspace{2pt}

Fix an integer $n<s$. Since ${\rm rad}^{\hspace{.6pt}2\hspace{-2pt}}\La=0$, we see that $d_{F(I_{a^{\rm o}})}^{\hspace{0.4pt}n}\vspace{1pt}$ vanishes on $\mathit\Theta_n$. Consider a top element $\bar{\varepsilon}_x \otimes p^{\rm o}$ of $F(I_{a^{\rm o}})^n$, where $x\in Q^{-n}$ and $p\in Q(a, x)$.
Since $p$ is of length $s-n>0,$ we have $p=\gamma \circ p_*$, where $\gamma \in Q_1(y, x)$ with $y\in Q^{-n-1}$
and $p_*\in Q(a, y)$. By definition, $(P[\gamma]\otimes I_{a^{\rm o}}(\gamma^{\rm o}))(\bar{\varepsilon}_x\otimes p^{\rm o})=\bar{\gamma}\otimes p_*^{\rm o}.$ Moreover, if $\alpha\in Q_1(x, b)$ with $b\in Q^{-n-1}$ is different from $\gamma$, then $p$ does not factor through $\alpha$, and hence, $p^{\rm o}$ does not factor through $\alpha^{\rm o}$. By definition, $I_{a^{\rm o}}(\alpha)(p^{\rm o})=0$, and $(P[\alpha]\otimes I_{a^{\rm o}}(\alpha^{\rm o}))(\bar{\varepsilon}_x\otimes p^{\rm o})=0$. Thus, $d_{F(I_{a^{\rm o}})}^n\vspace{1pt}$ sends $\bar{\varepsilon}_x\otimes p^{\rm o}$ to
$\bar{\gamma}\otimes p_*^{\rm o}\in \mathit\Theta_{n+1}.\vspace{1pt}$

Conversely, let $\bar{\beta}\otimes q^{\rm o}\in \mathit\Theta_{n+1}$, where $\beta\in Q_1(y, x)$ with $(y, x)\in Q^{-n-1}\times Q^{-n}$ and $q\in Q(a, z)$. Then
$\varepsilon_x\otimes q^{\rm o}\beta^{\rm o}\vspace{1pt}\in \mathit\Omega_n$ which, as we have seen, is sent by $d_{F(I_{a^{\rm o}})}^n\vspace{1pt}$ to $\bar{\beta}\otimes q^{\rm o}$. As a consequence, ${\rm Im} \hspace{0.6pt} d_{F(I_{a^{\rm o}})}^n\vspace{1pt}={\rm rad} \hspace{0.6pt} F(I_{a^{\rm o}})^{n+1}.$ This in turn implies that
${\rm Ker}\hspace{0.6pt} d_{F(I_{a^{\rm o}})}^n ={\rm rad}\hspace{0.6pt}F(I_{a^{\rm o}})^n={\rm Im}\hspace{0.5pt}d_{F(I_{a^{\rm o}})}^{\hspace{0.6pt}n-1}.\vspace{1.5pt}$

To sum up, it is shown that $F(I_{a^{\rm o}})^n=0$ for all $n>s,$ and ${\rm H}^s(F(I_{a^{\rm o}})^\cdt)\cong S[a]$ and ${\rm H}^n(F(I_{a^{\rm o}})^\cdt)=0$ for all $n<s$. That is, $F(I_{a^{\rm o}})^\cdt\cong P^\pdt_{S[a]}[-s].$ The proof of the lemma is completed.

\medskip

The following statement determines those representations whose image under the Koszul functor is bounded-above with bounded homology.

\medskip

\begin{Prop}\label{F-bch}

Let $\La=kQ/(kQ^+)^2$ with $Q$ a gradable locally finite quiver, and let $F: {\rm Rep}^*(Q^{\rm \hspace{.5pt}op})\to RC({\rm Proj}\hspace{0.6pt}\La)$ be the Koszul functor. If $M\in {\rm Rep}^*(Q^{\rm \hspace{0.4pt}op})$, then

\vspace{-1.5pt}

\begin{enumerate}[$(1)$]

\item $F(M)^\pdt\in RC^{-, b}({\rm Proj}\hspace{.4pt}\La)$ if and only if $M\in {\rm Rep}^-(Q^{\rm \hspace{0.4pt} op});$

\vspace{1pt}

\item $F(M)^\pdt\in RC^{-, b}({\rm proj}\hspace{.4pt}\La)$ if and only if $M\in {\rm rep}^-(Q^{\rm \hspace{0.4pt} op});$

\vspace{1pt}

\item $F(M)^\pdt\in RC^{b}({\rm proj}\hspace{.5pt}\La)$ if and only if $M\in {\rm rep}^b(Q^{\rm \hspace{0.4pt} op}).$

\end{enumerate} \end{Prop}

\noindent{\it Proof.} Let $M$ be a representation in ${\rm Rep}^*(Q^{\rm \hspace{0.4pt}op})$ with an injective co-presentation $\xymatrixcolsep{16pt}\xymatrix{0 \ar[r] & M \ar[r] & I^{\hspace{.4pt}0} \ar[r] & I^1 \ar[r] & 0,}$ where
$I^{\hspace{.4pt}0}, I^1\in {\rm Inj}(Q^{\rm \hspace{.5pt}op}).$ Applying the exact functor $F$, we obtain a short exact sequence
$\xymatrixcolsep{16pt}\xymatrix{0 \ar[r] & F(M)^\pdt \ar[r] & F(I^{\hspace{.4pt}0})^\pdt \ar[r] & F(I^1)^\pdt \ar[r] & 0}$ in $C(\ModbLa)$. By Lemma \ref{F-inj-image}, $F(I^{\hspace{.4pt}0})^\pdt$ and $F(I^1)^\pdt$ belong to $RC^{-, b}({\rm Proj}\hspace{.4pt}\La)$. As a consequence, $F(M)^\cdt\in RC^{-, b}({\rm Proj}\hspace{.4pt}\La).$ Moreover, since there exists an exact sequence
$\xymatrixcolsep{16pt}\xymatrix{{\rm H}^{n-1}(F(I^1)^\pdt) \ar[r] & {\rm H}^n(F(M)^\pdt) \ar[r] & {\rm H}^n(F(I^{\hspace{.4pt}0})^\pdt)}$ for all $n\in \Z$; see \cite[(1.3.1)]{Wei}, we conclude that $F(M)^\pdt\in RC^{-, b}({\rm Proj}\hspace{.4pt}\La)$.

Conversely, assume that there exist $s, t$ with $s<t$ such that $F(M)^n=0$ for all $n>t$, and ${\rm H}^n(F(M)^\cdt)=0$ for all $n<s.$ In particular, the full subquiver $\Sa(M)$ of $Q$ generated by the vertices $x$ with $M(x)\ne 0$ is contained in $\cup_{n\le t}\hspace{.4pt} (Q^{\rm op})^n$. As a consequence, the full subquiver $(Q^{\rm \hspace{.5pt}op})^{\le s}$ of $Q^{\rm \hspace{.5pt}op}$ generated by the vertices in $\cup_{n\le s}\hspace{0.4pt}Q^{-n}$ is co-finite in $\mathit\Sigma(M)$, that is, its complement in $\mathit\Sigma(M)$ is finite.

Let $N$ be the restriction of $M$ to $(Q^{\rm \hspace{.5pt}op})^{\le s}$; see \cite[Section 1]{BLP}, and denote by $S$ the image of $d_{F(M)}^{\,s}$. Being a radical complex, $F(N)^\cdt$ is isomorphic to the shift by $[-s]$ of the deleted minimal projective resolution of $S$. If $S=0$, then $F(N)^\pdt=0.$ Since $F$ is faithful, $N=0$. In particular, $N\in {\rm Inj}(Q^{\rm \hspace{.5pt}op})$. Otherwise, since ${\rm rad}^{\hspace{.5pt}2}\hspace{-2pt}\La=0$, we see that $S$ is semi-simple, say $S= \left(S[x_1]\otimes U_1\right) \oplus \cdots \oplus \left(S[x_t]\otimes U_t\right),$ where $x_1, \ldots, x_t\in Q_0$ and $U_1, \ldots, U_t$ are non-zero $k$-spaces. By Lemma \ref{F-inj-image}, we have $$F(N)^\pdt\cong \left(F(I_{x_1^{\rm o}})^\pdt\otimes U_1\right)[s_1] \oplus \cdots \oplus \left(F(I_{x_t^{\rm o}})^\pdt\otimes U_t\right)[s_t],$$ where $s_1, \ldots, s_t$ are some integers. Since $F$ is fully faithful,
we have a decomposition $N=N_1 \oplus \cdots \oplus N_t$ such that
$F(N_i)^\pdt \cong \left(F(I_{x_i})^\pdt\otimes U_i\right)[s_i] = F(I_{x_i}\otimes U_i)^\pdt [s_i],$ for $i=1, \ldots, t.$ By Lemma \ref{Morph-F-images}(1), $s_i=0$, and consequently, $N_i\cong I_{x_i}\otimes U_i$, for
$i=1, \ldots, t.$ That is, $N\in {\rm Inj}(Q^{\rm \hspace{.5pt}op})$. In view of Theorem 1.12(2) in \cite{BLP}, we conclude that $M$ is almost finitely co-presented. This establishes Statement (1).

For Statement (2), it suffices to prove the necessity. Indeed, suppose that $F(M)^\pdt$ lies in $RC^{-, b}({\rm proj}\hspace{.4pt}\La)$. By Statement (1), $M$ is almost finitely co-presented, and hence, ${\rm soc}\,M$ is finitely supported and essential in $M$; see \cite[(1.6)]{BLP}. Moreover, by the definition of $F$, we see that $M$ is locally finite dimensional, and so is ${\rm soc}\,M$. Therefore, ${\rm soc}\,M$ is finite dimensional. As a consequence, $M$ is finitely co-presented. This establishes Statement (2), from which Statement (3) follows easily. The proof of the proposition is completed.

\medskip

We shall extend the Koszul functor to bounded complexes. Indeed, in view of Proposition \ref{F-bch}, we see that $F: {\rm Rep}^*(Q^{\hspace{0.5pt}\rm op})\to RC({\rm Proj}\hspace{.5pt}\La)$ restricts to a functor $F: {\rm Rep}^-(Q^{\hspace{0.5pt}\rm op})\to RC^{-,b}({\rm Proj}\hspace{.5pt}\La)$, called again {\it Koszul functor}.  Given a complex $M^\ydt\in C^b({\rm Rep}^-(Q^{\rm \hspace{.5pt}op}))$, applying $F$ to the components of $M^\ydt$, we obtain a double complex $F(M^\pdt)^\cdt$ over $\mbox{Proj}\hspace{.4pt}\La$ as follows: \vspace{-2pt}
$$\xymatrix{
             &  \vdots                                                    &&             \vdots\\
\cdots\ar[r] & F(M^i)^{j+1} \ar[u] \ar[rr]^{F(d_M^i)^{j+1}}  && F(M^{i+1})^{j+1} \ar[u] \ar[r]
& \cdots\\
\cdots\ar[r] & F(M^i)^j \ar[rr]^{F(d_M^i)^j} \ar[u]^{(-1)^i d_{F(M^i)}^j} && F(M^{i+1})^j \ar[u]_{(-1)^{i+1}d_{F(M^{i+1})}^j} \ar[r] & \cdots\\
&\vdots\ar[u] && \vdots\ar[u] \\
}$$


Let $\mathcal{F}(M^\pdt)$ be the sum total complex of $F(M^\pdt)^{\cdt}$, which lies in $C({\rm Proj}\hspace{.6pt}\La)$. More explicitly, for each $n\in \Z$, we have $\mathcal{F}(M^\pdt)^n=\oplus_{i\in \Z} \, F(M^i)^{n-i}\in {\rm Proj}\hspace{.6pt}\La;$ and
$$d^n_{\mathcal{F}(M^\pdt)}=(d^n_{\mathcal{F}(M^\pdt)}(j,i))_{(j,i)\in \Z\times \Z}:
\oplus_{i\in \mathbb{Z}} \, F(M^i)^{n-i} \to   \oplus_{j\in \mathbb{Z}} \, F(M^j)^{n+1-j},$$
where $d^n_{\mathcal{F}(M^\pdt)}(j,i):  F(M^i)^{n-i} \to F(M^j)^{n+1-j}$ is given by
$$d^n_{\mathcal{F}(M^\ydt)}(j,i) = \left\{\begin{array}{ll}
(-1)^i d_{F(M^i)}^{\hspace{.4pt}n-i}, & j=i;   \\ \vspace{-9pt} \\
F(d_M^i)^{n-i},             & j=i+1; \\ \vspace{-9pt} \\
0,                          & j\ne i, i+1.
\end{array}\right.
\vspace{5pt}$$

Since $M^\cdt$ is bounded with $F(M^i)^\cdt$ bounded above, we see that $\mathcal{F}(M^\pdt)\in C^-({\rm Proj}\hspace{.6pt}\La).$ Next, let $f^\ydt: M^\ydt\to N^\ydt$ be a morphism in $C^b({\rm Rep}(Q^{\rm \hspace{.5pt}op}))$. For each $n\in \Z$, we set
$$\mathcal{F}(f^\ydt)^n=\oplus_{i\in \Z } \{F(f^i)^{n-i}\}: \oplus_{i\in \Z} \, F(M^i)^{n-i} \to \oplus_{i\in \Z} \, F(N^i)^{n-i}.$$

For each pair $(i, j)\in \Z\times \Z$, we have
$$F(f^j)^{n+1-j}\circ d^n_{\mathcal{F}(M^\ydt)}(j,i)= \left\{\begin{array}{ll}
(-1)^i F(f^j)^{n+1-j} \circ d_{F(M^i)}^{\hspace{.4pt}n-i}, & j=i;   \\ \vspace{-9pt} \\
F(f^j)^{n+1-j} \circ F(d_M^i)^{n-i},                       & j=i+1; \\ \vspace{-9pt} \\
0,                                                         & j\ne i, i+1,
\end{array}\right.$$
and

$$d^n_{\mathcal{F}(N^\ydt)}(j,i)\, F(f^i)^{n-i} = \left\{\begin{array}{ll}
(-1)^i d_{F(N^i)}^{\hspace{.4pt}n-i} \circ F(f^i)^{n-i}, & j=i;   \\ \vspace{-9pt} \\
F(d_N^i)^{n-i} \circ F(f^i)^{n-i},                       & j=i+1; \\ \vspace{-9pt} \\
0,                                                       & j\ne i, i+1.
\end{array}\right.$$

Since $F(f^i)^\cdt: F(M^i)^\cdt \to F(N^i)^\cdt$ is a morphism of complexes, we obtain
$$F(f^i)^{n+1-i} d_{F(M^i)}^{n-i}=d_{F(N^i)}^{n-i} \, F(f^i)^{n-i};$$
and since $f^{i+1} d_M^i=d_N^i f^i,\vspace{1pt}$ we have $F(f^{i+1})^\cdt\, F(d_M^i)^\cdt=F(d_N^i)^\cdt\, F(f^i)^\cdt$. In particular, we have $F(f^{i+1})^{n-i}\, F(d_M^i)^{n-i}=F(d_N^i)^{n-i}\, F(f^i)^{n-i}.$
This gives rise to an equation
$\mathcal{F}(f^\pdt)^{n+1} \circ d^n_{\mathcal{F}(M^\ydt)}=d^n_{\mathcal{F}(N^\ydt)} \circ \mathcal{F}(f^\pdt)^n.\vspace{1.5pt}$ Therefore, the construction yields a functor
$\mathcal{F}: C^b({\rm Rep}^-(Q^{\rm \hspace{.5pt}op}))\to C^{-}({\rm Proj}\hspace{.6pt}\La),$ called {\it complex Koszul functor}.

\medskip

\begin{Lemma}\label{hat-F-functor}
Let $\La=kQ/(kQ^+)^2,\vspace{-2pt}$ where $Q$ is a gradable locally finite quiver. The complex Koszul functor is a $k$-linear functor
$\xymatrixcolsep{16pt}\xymatrix{\mathcal{F}: C^b({\rm Rep}^-(Q^{\rm \hspace{.5pt}op}))\ar[r] & C^{-,b}({\rm Proj}\hspace{.6pt}\La),
}\vspace{-1pt}$ which sends acyclic complexes to acyclic ones and restricts to the Koszul functor \vspace{.8pt}
$F: {\rm Rep}^-(Q^{\hspace{0.5pt}\rm op})\to RC^{-,b}({\rm Proj}\hspace{.5pt}\La)$.

\end{Lemma}

\smallskip

\noindent{\it Proof.} First of all, for any $M\in {\rm Rep}^-(Q^{\hspace{.5pt}\rm op})$, it follows easily from the construction that $\mathcal{F}(M)=F(M)^\cdt$. Let $M^\ydt\in C^b({\rm Rep}^-(Q^{\rm \hspace{.5pt}op}))$. With no loss of generality, we may assume that there exists an integer $s<0$ such that $M^n\ne 0$ only if $s\le i\le 0$.


If $M^\cdt$ is acyclic then, since $F$ is exact, the double complex $F(M^\ydt)^\ydt$ has exact rows. Since the complexes $F(M^i)^\cdt$ with $s\le i\le 0$ are bounded above, we deduce from the Acyclic Assembly Lemma stated in \cite[(2.7.3)]{Wei} that $\mathcal{F}(M^\ydt)$ is acyclic.

It remains to show that $\mathcal{F}(M^\ydt)\in C^{-, b}({\rm Proj}\hspace{.6pt}\La)$. Indeed, we have seen that $\mathcal{F}(M^\ydt)\in C^-({\rm Proj}\hspace{.6pt}\La)$. Since the $F(M^i)^\cdt$ with $s\le i\le 0$ have bounded homology, there exists an integer $m<0$ such that ${\rm H}^n(F(M^i)^\ydt)=0$ for $i\in \Z$ and $n\le m$. We define a double complex $(L^\cdt\hspace{.6pt}^\cdt, v^\cdt\hspace{.6pt}^\cdt, h^\cdt\hspace{.6pt}^\cdt)$ as follows. The objects are given by
$$L^{ij}=\left\{\begin{array}{lll} 0\hspace{.5pt}, & i>m+1;\\
{\rm Im}\, d_{F(M^i)}^m\hspace{.5pt}, & i=m+1; \\ \vspace{-10pt}\\
F(M^i)^j\hspace{.5pt}, & i\le m.
\end{array}\right.$$
The vertical morphisms are given by
$$v^{ij}=\left\{\begin{array}{lll} 0\hspace{.5pt}, & i\ge m+1;\\
v^{i,m}, & i=m; \\
(-1)^id_{F(M^i)}^j\hspace{.5pt}, & i < m,
\end{array}\right.$$ where $v^{i, m}$ is the co-restriction of $(-1)^id_{F(M^i)}^m$ to its image. The horizontal morphisms are given by
$$h^{ij}=\left\{\begin{array}{lll} 0\hspace{.5pt}, & j > m+1;\\
h^{i,m+1}, & i=m+1; \\
F(d_M^i)^j\hspace{.5pt}, & j \le m,
\end{array}\right.
\vspace{1pt}$$ where $h^{i, m+1}$ is induced from $F(d_M^i)^m$. This yields a lower half-plane double complex $L^\cdt\hspace{.6pt}^\cdt$ with exact columns. By the Acyclic Assembly Lemma in \cite[(2.7.3)]{Wei}, the sum total complex $(N^\cdt, d_N^\ydt)$ of $L^\cdt\hspace{.6pt}^\cdt$ is acyclic. On the other hand, by the construction, we see that ${\rm H}^n(N^\cdt)={\rm H}^n(\mathcal{F}(M^\ydt))$ for all $n<s+m$. Thus, $\mathcal{F}(M^\ydt)\in C^{-, b}({\rm Proj}\hspace{.6pt}\La)$. The proof of the lemma is completed.

\medskip

We shall collect more properties of the complex Koszul functor in the following statement. For this purpose, we denote by $C_{f^\ydt}$ the mapping cone of a morphism $f^\cdt$ of complexes over an additive category.

\medskip

\begin{Lemma}\label{cal-F-property}

Let $\La=kQ/(kQ^+)^2\vspace{-2pt}$ with $Q$ a gradable locally finite quiver, and let $\xymatrixcolsep{16pt}\xymatrix{\mathcal{F}: C^b({\rm Rep}^-(Q^{\rm \hspace{.5pt}op}))\ar[r] & C^{-,b}({\rm Proj}\hspace{.6pt}\La)}$ be the complex Koszul functor.

\vspace{-3pt}

\begin{enumerate}[$(1)$]

\item If $M^\ydt$ is an object of $C^b({\rm Rep}^-(Q^{\rm \hspace{.5pt}op}))$, then $\mathcal{F}(M^\ydt[1]) = \mathcal{F}(M^\ydt)[1].$

\vspace{2pt}

\item If $f^\ydt\in C^b({\rm Rep}^-(Q^{\rm \hspace{.5pt}op}))$ is a null-homotopic, then $\mathcal{F}(f^\cdt)$ is null-homotopic.

\vspace{2pt}

\item If $f^\ydt$ is a morphism in $C^b({\rm Rep}^-(Q^{\rm \hspace{.5pt}op}))$, then $C_{\mathcal{F}(f^\cdt)}=\mathcal{F}(C_{\hspace{-1pt}f^\ydt})$.

\end{enumerate}\end{Lemma}

\noindent{\it Proof.} Statement (1) can be shown by a routine verification. Consider a morphism $f^\ydt: M^\ydt\to N^\cdt$ in $C^b({\rm Rep}^-(Q^{\rm \hspace{.5pt}op}))$. Suppose that $h^i: M^i\to N^{i-1}$, $i\in \Z$, are morphisms such that $f^i=h^{i+1}\circ d_M^i+d_N^{i-1}\circ h^i$, for every $i\in \Z$. Fix $n\in \Z$. We set
$$\varphi^n=(\varphi^n(j,i))_{(j, i)\in \mathbb{Z}\times \mathbb{Z}}: \oplus_{i\in \mathbb{Z}}F(M^i)^{n-i}\to
\oplus_{j\in \mathbb{Z}}F(N^j)^{n-1-j},$$
where $\varphi^n(j,i): F(M^i)^{n-i}\to F(N^j)^{n-1-j}$ is given by
$$\varphi^n(j,i)=\left\{\begin{array}{ll}
F(h^i)^{n-i}, & j=i-1;\\
0, & j\ne i-1.
\end{array}\right.$$

Given a pair $(j, i)\in \Z\times \Z$, we have
$$\begin{array}{rcl}
\left(\varphi^{n+1}\, d_{\mathcal{F}(M^\ydt)}^n\right)(j, i) &=& \varphi^{n+1}(j, j+1)\, d_{\mathcal{F}(M^\ydt)}^n(j+1, i)\\ \vspace{-8pt}\\
&=&\left\{\begin{array}{ll}
(-1)^i F(h^i)^{n+1-i}\, d_{F(M^i)}^{\hspace{.4pt}n-i}, & j=i-1;   \\ \vspace{-9pt} \\
F(h^{i+1})^{n-i}\,F(d_M^i)^{n-i},             & j=i; \\ \vspace{-9pt} \\
0,                          & j\ne i-1, i.
\end{array}\right.
\end{array}$$
and
$$\begin{array}{rcl}
\left(d_{\mathcal{F}(N^\ydt)}^{n-1} \, \varphi^n \right)(j, i) &=& d_{\mathcal{F}(N^\ydt)}^{n-1}(j, i-1) \, \varphi^n(i-1, i) \\ \vspace{-8pt}\\
&=&\left\{\begin{array}{ll}
(-1)^{i-1} d_{F(N^{i-1})}^{\hspace{.4pt}n-i} \, F(h^i)^{n-i}, & j=i-1;   \\ \vspace{-9pt} \\
F(d_N^{i-1})^{n-i} \, F(h^i)^{n-i},             & j=i; \\ \vspace{-9pt} \\
0,                          & j\ne i-1, i.
\end{array}\right.
\end{array}$$

Since $F(h^i)^\cdt: F(M^i)^\cdt\to F(N^{i-1})^\cdt$ is a complex morphism, we conclude that $d_{F(N^{i-1})}^{\hspace{.4pt}n-i} \, F(h^i)^{n-i}=F(h^i)^{n+1-i}\, d_{F(M^i)}^{\hspace{.4pt}n-i}.$ In view of this equation, we deduce that
$(\varphi^{n+1}\, d_{\mathcal{F}(M^\ydt)}^n + d_{\mathcal{F}(N^\ydt)}^{n-1} \, \varphi^n)(j, i)=0\vspace{1pt}$ for any $j\ne i$. On the other hand, since $d_M^i \, h^{i+1}+d_N^{i-1} h^i=f^i$, we have
$$F(h^{i+1})^{n-i}\,F(d_M^i)^{n-i}+F(d_N^{i-1})^{n-i} \, F(h^i)^{n-i}=F(f^i)^{n-i}.$$
Therefore, $\varphi^{n+1}\circ d_{\mathcal{F}(M^\ydt)}^n+d_{\mathcal{F}(N^\ydt)}^{n-1}\circ \varphi^n=\mathcal{F}(f^\pdt)^n.
\vspace{2.5pt}$ This establishes Statement (2).

\vspace{2pt}

In order to prove Statement (3), we write $(C^\pdt, d_C^\ydt)$ for the mapping cone of $f^\ydt: M^\ydt\to N^\cdt$. We shall show that $\mathcal{F}(C^\pdt)=C_{\mathcal{F}(f^\ydt)}$. Fix an integer $n$.
By definition, we have $F(C^n)^\cdt=F(M^{n+1})^\cdt \oplus F(N^n)^\cdt$ and \vspace{-3pt}
$$
d_{F(C^n)}^{\pdt}=\left(\begin{array}{cc}
d_{F(M^{n+1})}^\pdt & \hspace{-8pt} 0                            \\ \vspace{-8pt}\\
                  0 & \hspace{-8pt} d_{F(N^n)}^{\pdt}
\end{array}\right); \; F(d_C^{\hspace{.6pt}n})^\cdt=\left(\begin{array}{cc}
- F(d_M^{n+1})^\cdt & \hspace{-8pt} 0\\ \vspace{-8pt}\\
F(f^{n+1})^\cdt & \hspace{-8pt} F(d_N^n)^{\hspace{.6pt}\cdt}
\end{array}\right).
$$

First, we have
$$\begin{array}{rcl}
C_{\mathcal{F}(f^\ydt)}^n
&=& \mathcal{F}(M^\ydt)^{n+1}\oplus \mathcal{F}(N^\cdt)^n\\
&=& \left(\oplus_{j\in \Z}\, F(M^j)^{n+1-j} \right) \oplus \left(\oplus_{i\in \Z}\, F(N^i)^{n-i} \right) \\ \vspace{-9pt}\\
&=& \oplus_{i\in \Z}\, F(M^{i+1} \oplus  N^i)^{n-i} \\ \vspace{-9pt}\\
&=& \oplus_{i\in \Z}\, F(C^i)^{n-i} 
= \mathcal{F}(C^\ydt)^n.
\end{array}$$

Next, we shall identify their $n$-th differentials
$$d_{C_{\hspace{-1pt}\mathcal{F}(f^\ydt)}}^n \hspace{-2pt} = \hspace{-2pt}
\left(\hspace{-5pt}\begin{array}{cc}
- d_{\mathcal{F}(M^\ydt)}^{n+1} & \hspace{-3pt} 0\\ \vspace{-8pt}\\
\mathcal{F}(f^\ydt)^{n+1} & \hspace{-3pt} d_{\mathcal{F}(N^\cdt)}^n
\end{array}\hspace{-5pt}\right)\hspace{-3pt}: \mathcal{F}(M^\ydt)^{n+1}\oplus \mathcal{F}(N^\cdt)^n\to \mathcal{F}(M^\ydt)^{n+2}\oplus \mathcal{F}(N^\cdt)^{n+1}$$

and $d_{\mathcal{F}(C^\ydt)}^n=(d_{\mathcal{F}(C^\ydt)}^n(j,i))_{(j,i)\in \Z\times \Z}: \;
\oplus_{i\in \Z} F(C^i)^{n-i} \to \oplus_{j\in \Z} F(C^j)^{n+1-j},$
where $$d^n_{\mathcal{F}(C^\ydt)}(j,i)=\left\{\begin{array}{ll}
(-1)^i d_{F(C^i)}^{\hspace{.4pt}n-i}, & j=i;   \\ \vspace{-9pt} \\
F(d_C^i)^{n-i},             & j=i+1; \\ \vspace{-9pt} \\
0,                          & j\ne i, i+1.
\end{array}\right.$$

\medskip

It suffices to identify, for any $(j, i)\in \Z\times \Z$ and $(X^\cdt, Y^\cdt)$ with $X^\cdt, Y^\cdt\in \{M^\cdt, N^\cdt\}$, the composite
$$g_{j,i}(X^\cdt, Y^\cdt): \xymatrix{F(X^{i+1})^{n-i}\ar[r]^-{q_1} & C_{\mathcal{F}(f^\ydt)}^n \ar[r]^-{d^n_{C_{\hspace{-1.5pt}\mathcal{F}(f^\ydt)}}} &
C_{\mathcal{F}(f^\ydt)}^{n+1} \ar[r]^-{p_1} & F(Y^{j+1})^{n+1-j}}$$ with the composite
$$h_{j,i}(X^\cdt, Y^\cdt): \xymatrix{F(X^{i+1})^{n-i}\ar[r]^-{q_2} & \mathcal{F}(C^\ydt)^n\ar[r]^-{d^n_{\mathcal{F}(C^\ydt)}} & \mathcal{F}(C^\ydt)^{n+1} \ar[r]^-{p_2} & F(Y^{j+1})^{n+1-j},}$$ where $q_1, q_2$ are the canonical injections, and $p_1, p_2$ are the canonical projections. We start with the pair $(M^\cdt, N^\cdt)$. It is easy to see that both composites
$$g_{j,i}(M^\cdt, N^\cdt): \xymatrix{F(N^i)^{n-i}\ar[r]^-{q_1} & C_{\hspace{-1pt}\mathcal{F}(f^\ydt)}^n \ar[r]^-{d^n_{C_{\mathcal{F}(f^\ydt)}}} &
C_{\mathcal{F}(f^\ydt)}^{n+1} \ar[r]^-{p_1} & F(M^{j+1})^{n+1-j}}$$
and
$$h_{j,i}(M^\cdt, N^\cdt): \xymatrix{F(N^i)^{n-i}\ar[r]^-{q_2} & \mathcal{F}(C^\pdt)^n\ar[r]^-{d^n_{\mathcal{F}(C^\pdt)}} &\mathcal{F}(C^\pdt)^{n+1} \ar[r]^-{p_2} & F(M^{j+1})^{n+1-j}}$$
are null. Next, we consider the pair $(M^\cdt, M^\cdt)$. Observe that the composite
$$g_{j,i}(M^\cdt, M^\cdt): \xymatrix{ F(M^{i+1})^{n-i}\ar[r]^-{q_1} & C_{\mathcal{F}(f^\ydt)}^n \ar[r]^-{d^n_{C_{\hspace{-1.5pt}\mathcal{F}(f^\ydt)}}} &
C_{\mathcal{F}(f^\ydt)}^{n+1} \ar[r]^-{p_1} & F(M^{j+1})^{n+1-j}}$$
is the composite of the following morphisms:
$$\xymatrix{F(M^{i+1})^{n-i}\ar[r]^-{q} & \mathcal{F}(M^{\ydt})^{n+1} \ar[rr]^-{-d_{\mathcal{F}(M^\ydt)}^{n+1}} &&
\mathcal{F}(M^\ydt)^{n+2} \ar[r]^-{p} & F(M^{j+1})^{n+1-j}.}$$
That is,
$$g_{j,i}(M^\cdt, M^\cdt)=-d_{\mathcal{F}(M^\pdt)}^{n+1}(j+1, i+1)=\left\{\begin{array}{ll}
(-1)^{i+2} d_{F(M^{i+1})}^{\hspace{.7pt}n-i}, & j=i;\\ \vspace{-10pt}\\
-F(d_M^{\hspace{.7pt}i+1})^{n-i},            & j=i+1; \\ \vspace{-10pt}\\
0, & j\ne i, i+1.
\end{array}\right.$$

\medskip

On the other hand, the composite
$$h_{j,i}(M^\cdt, M^\cdt): \xymatrix{F(M^{i+1})^{n-i}\ar[r]^-{q_2} & \mathcal{F}(C^\ydt)^n\ar[r]^-{d^n_{\mathcal{F}(C^\ydt)}} & \mathcal{F}(C^\ydt)^{n+1} \ar[r]^-{p_2} & F(M^{j+1})^{n+1-j}}$$
is the composite of the following morphisms
$$\xymatrix{F(M^{i+1})^{n-i} \ar[r]^{q} & F(C^i)^{n-i} \ar[rr]^{d^n_{\mathcal{F}(C^\ydt)}(j,i)}
&& F(C^j)^{n+1-j} \ar[r]^{p} & F(M^{j+1})^{n+1-j}.
}$$

If $j\ne i, i+1$, since $d^n_{\mathcal{F}(C^\ydt)}(j,i)=0$, we obtain $h_{j,i}(M^\cdt, M^\cdt)=0.$ If $j=i$, then $h_{i,i}(M^\cdt, M^\cdt)$ is the composite of the following morphisms
$$\xymatrix{F(M^{i+1})^{n-i} \ar[r]^{q} & F(C^i)^{n-i} \ar[rr]^-{(-1)^i d^{n-i}_{F(C^i)}} && F(C^i)^{n+1-i} \ar[r]^p &
F(M^{i+1})^{n+1-i},}$$
that is, $h_{i,i}(M^\cdt, M^\cdt)=(-1)^id_{F(M^{i+1})}^{n-i}.$ If $j=i+1$, then $h_{i+1,i}(M^\cdt, M^\cdt)$ is the composite of the morphisms
$$\xymatrix{F(M^{i+1})^{n-i} \ar[r]^{q} & F(C^i)^{n-i} \ar[rr]^{F(d_C^i)^{n-i}} && F(C^{i+1})^{n-i} \ar[r]^p &
F(M^{i+2})^{n-i}.}$$
Thus, $h_{i+1, i}(M^\cdt, M^\cdt)=-F(d_M^{i+1})^{n-i}$. This shows that $g_{j,i}(M^\cdt, M^\cdt)=h_{j,i}(M^\cdt, M^\cdt).\vspace{1pt}$
Similarly, $g_{j,i}(N^\cdt, M^\cdt)=h_{j,i}(N^\cdt, M^\cdt)$ and $g_{j,i}(N^\cdt, N^\cdt)=h_{j,i}(N^\cdt, N^\cdt)$. The proof of the lemma is completed.

\medskip

In view of Lemmas \ref{hat-F-functor} and \ref{cal-F-property}, we deduce that the complex Koszul functor \vspace{-1pt} $\xymatrixcolsep{16pt}\xymatrix{\mathcal{F}: C^b({\rm Rep}^-(Q^{\rm \hspace{.5pt}op}))\ar[r] & C^{-,b}({\rm Proj}\hspace{.6pt}\La)}$ induces
a commutative diagram \vspace{-8pt}

$\begin{array}{lllll}
\hspace{-20pt} \begin{array}{c}  \\ \\ \\
 (3.1) \end{array}  &&&
\xymatrix{C^b({\rm Rep}^-(Q^{\rm \hspace{.5pt}op})) \ar[r]^{\mathcal{P}_{Q^{\rm \hspace{.5pt}op}}}\ar[d]^-{\mathcal{F}} & K^b({\rm Rep}^-(Q^{\rm \hspace{.5pt}op})) \ar[r]^{\mathcal{L}_{Q^{\rm \hspace{.5pt}op}}} \ar[d]^-{\mathscr{F}} &  D^b({\rm Rep}^{-}(Q^{\rm \hspace{.5pt}op})) \ar[d]^-{\mf{F}} \\
C^{-, b}({\rm Proj}\hspace{.6pt}\La) \ar[r]^{\mathcal{P}_{_{\hspace{-1pt}\it\Lambda}}} & K^{-, b}({\rm Proj}\hspace{.6pt}\La) \ar[r]^{\mathcal{E}_{_{\hspace{-1pt}\it\Lambda}}} & D^b(\ModbLa),}
\end{array}$

\medskip

\noindent where $\mathscr{F}, \mf{F}$ are triangle-exact. We shall call $\mf{F}$ {\it derived Koszul functor} and show that it is actually a triangle-equivalence.

\smallskip

\begin{Lemma}\label{Equiv-lemma}

Let $\La=kQ/(kQ^+)^2\vspace{-1pt}$ with $Q$ a gradable locally finite quiver, and let $\xymatrixcolsep{16pt}\xymatrix{\mf{F}: D^b({\rm Rep}^-(Q^{\rm \hspace{.5pt}op}))\ar[r] & D^b(\ModbLa)}\hspace{-2pt}$ be the derived Koszul functor.

\vspace{-2pt}

\begin{enumerate}[$(1)$]

\item If $X^\pdt$ is a complex in $D^b(\ModbLa),$ then
$X^\ydt\cong \mf{F}(M_1)[n_1]\oplus \cdots \oplus \mf{F}(M_s)[n_s],$
where $M_1, \ldots, M_s\in {\rm Rep}^-(Q^{\rm \hspace{.5pt}op})$ and $n_1,$ $\ldots,$ $n_s\in \Z$.

\vspace{2pt}

\item If $M, N\in {\rm Rep}^-(Q^{\rm \hspace{.5pt}op})$ and $s\in \Z$ such that $\Hom_{D^b({\rm Mod}^b\hspace{-2pt}\it\Lambda)}(\mf{F}(M), \mf{F}(N)[s])\ne 0,\vspace{.5pt}$ then $s=0, 1$.

\end{enumerate}

\end{Lemma}

\noindent{\it Proof.} Recall that $\mathcal{E}_{_{\hspace{-1pt}\it\Lambda}}: K^{-,b}(\mbox{Proj}\hspace{0.5pt}\La)\to D^b(\ModbLa)$ is a triangle-equivalence. Let $X^\pdt\in D^b(\ModbLa).$ By Proposition \ref{Der-Hmtp}(3), $X^\pdt\cong \mathcal{P}_{\it\Lambda}(\mathcal{E}_{\it\Lambda}(Y^\ydt))$ for a radical complex $Y^\cdt\in RC^{-,b}({\rm Proj}\hspace{.6pt}\La).$ By Corollary \ref{RC-fin-decom}, $Y^\cdt=Y_1^\cdt\oplus \cdots \oplus Y_s^\cdt$, where the $Y_i^\cdt$ are support-connected complexes in $RC^{-,b}({\rm Proj}\hspace{.6pt}\La).$ By Propositions \ref{Ind-complex} and \ref{F-bch}, we have $Y_i^\cdt\cong F(M_i)^\ydt[n_i]$ for some $M_i\in {\rm Rep}^-(Q^{\rm \hspace{.5pt}op})$ and $n_i\in \Z.$ Since $\mathcal{F}$ restricts to $F$; see (\ref{hat-F-functor}), we obtain $\mf{F}(M_i)[n_i]=\mathcal{E}_{\it\Lambda}(\mathcal{P}_{\it\Lambda}(F(M_i)^\cdt[n_i]))\cong \mathcal{E}_{\it\Lambda}(\mathcal{P}_{\it\Lambda}(Y_i^\cdt)),$ for $i=1, \ldots, s$. Thus,
$$X^\ydt\cong \mathcal{E}_{\it\Lambda}(\mathcal{P}_{\it\Lambda}(Y_1^\cdt)) \oplus \cdots \oplus
\mathcal{E}_{\it\Lambda}(\mathcal{P}_{\it\Lambda}(Y_s^\cdt))
\cong  \mf{F}(M_1)[n_1]\oplus \cdots \oplus \mf{F}(M_s)[n_s].$$

To prove Statement (2), let $M, N\in {\rm Rep}^-(Q^{\rm \hspace{.5pt}op})$ and $s\ne 0, 1$. By Lemma \ref{Morph-F-images}, $\Hom_{C({\rm Proj}\hspace{.6pt}\it\Lambda)}(F(M)^\cdt, F(N)^\cdt[s])=0,\vspace{1pt}$
that is, $\Hom_{C({\rm Proj}\hspace{.6pt}\it\Lambda)}(\mathcal{F}(M), \mathcal{F}(N)[s])=0.\vspace{1pt}$ In particular,
$\Hom_{K^{-, b}({\rm Proj}\hspace{.6pt}\it\Lambda)}(\mathscr{F}(M), \mathscr{F}(N)[s])=0.\vspace{1pt}$ Since $\mathcal{E}_{\it\Lambda}$ is an equivalence,
$\Hom_{D^b({\rm Mod}^b\hspace{-2pt}\it\Lambda)}(\mf{F}(M), \mf{F}(N)[s])=0.$ The proof of the lemma is completed.

\medskip

We are ready to obtain the main result of this section.

\medskip

\begin{Theo}\label{KzEqv}

Let $\La=kQ/(kQ^+)^2\vspace{1pt}$ with $Q$ a gradable locally finite quiver. The derived Koszul functor
$\mf{F}: D^b({\rm Rep}^-(Q^{\rm \hspace{.5pt}op}))\to D^b(\ModbLa)$ is a triangle-equivalence, which restricts to a triangle-equivalence $\mf{F}: D^b({\rm rep}^-(Q^{\rm \hspace{.5pt}op}))\to D^b({\rm mod}\hspace{.6pt}^b\hspace{-2pt}\La).$

\end{Theo}

\noindent{\it Proof.} We shall make use of the commutative diagram (3.1). In view of Proposition \ref{F-bch}(2), we need only to prove the first part of the theorem. By Lemma \ref{Equiv-lemma}(1), $\mf{F}$ is dense. To show that $\mf{F}$ is fully faithful, by Lemma \ref{Equiv-lemma}, it suffices to show that
$$\xymatrix{
\mf{F}_{M, N[s]}: \Hom_{D^b({\rm Rep}^-(Q^{\rm \hspace{.5pt}op}))}(M, N[s]) \ar[r] & \Hom_{D^b({\rm Mod}\hspace{.6pt}^b\hspace{-1.2pt}\it\Lambda))}(\mf{F}(M), \mf{F}(N)[s])}$$
is an isomorphism, for all $M, N\in {\rm Rep}^-(Q^{\rm \hspace{.5pt}op})$ and $s\in \Z$. If $s\ne 0, 1$, then $\Hom_{D^b({\rm Rep}^-(Q^{\rm \hspace{.5pt}op}))}(M, N[s])=0\vspace{2pt}$ since ${\rm Rep}^-(Q^{\rm \hspace{.5pt}op})$ is hereditary, and on the other hand, $\Hom_{D^b({\rm Mod}\hspace{.6pt}^b\hspace{-2pt}\it\Lambda))}(\mf{F}(M), \mf{F}(N)[s])=0\vspace{1pt}$ by Lemma \ref{Equiv-lemma}(2).
In case $s=0$, we have a commutative diagram
$$\xymatrix{
\Hom_{C^b({\rm Rep}^-(Q^{\rm \hspace{.5pt}op}))}(M, N) \ar[d]^{\cong} \ar[r]^-{\mathcal{F}_{M, N}} &  \Hom_{\,RC^{-,b}({\rm Proj}\hspace{.6pt}\it\Lambda)}(\mathcal{F}(M), \mathcal{F}(N)) \ar[d]^{\mathcal{P}_{_{M, N}}}\\
\Hom_{K^b({\rm Rep}^-(Q^{\rm \hspace{.5pt}op}))}(M, N) \ar[d]^{\cong} \ar[r]^-{\mathscr{F}_{M, N}} & \Hom_{K^{-,b}({\rm Proj}\hspace{.6pt}\it\Lambda)}(\mathscr{F}(M), \mathscr{F}(N))) \ar[d]^{\cong}\\
\Hom_{D^b({\rm Rep}^{-}(Q^{\rm \hspace{.5pt}op}))}(M, N) \ar[r]^-{\mf{F}_{M, N}} & \Hom_{D^b({\rm Mod}\hspace{.6pt}^b\hspace{-1.2pt}\it\Lambda))}(\mf{F}(M), \mf{F}(N)).
}$$

If $f^\cdt: \mathcal{F}(M)\to \mathcal{F}(N)$ is a null-homotopic in $RC^{-,b}({\rm Proj}\hspace{.6pt}\it\Lambda)$, then it is radical, and by Lemma \ref{Morph-F-images}(2), $f^\cdt=0$. This shows that $\mathcal{P}_{_{M, N}}\vspace{.5pt}$ is an isomorphism. Since $\mathcal{F}$ restricts to $F$, we deduce from Lemma \ref{F-property} that $\mathcal{F}_{M, N}$ is an isomorphism. As a consequence, $\mf{F}_{M, N}$ is an isomorphism.

\smallskip

It remains to consider the case where $s=1.$ Let $\theta^\ydt: \mf{F}(M)\to \mf{F}(N)[1]$ be a morphism in $D^b(\ModbLa)$, which embeds in an exact triangle
$$
\xymatrix{\mf{F}(N)\ar[r]^-{\phi^\pdt}& Y^\ydt\ar[r]^-{\psi^\pdt} & \mf{F}(M)\ar[r]^-{\theta^\ydt} &
\mf{F}(N)[1].
}$$

If $U\in {\rm Rep}^-(Q^{\rm \hspace{.5pt}op})$ and $t\in \Z$ are such that $\mf{F}(U)[t]\vspace{.5pt}$ is a non-zero summand of $Y^\cdt$, we claim that $t=0$. Indeed, writing $Y^\ydt=Z^\ydt \oplus \mf{F}(U)[t]$, we obtain an exact triangle
$$\xymatrix{
\mf{F}(N)\ar[r]^-{\zeta^\ydt \choose \mu^\cdt} & Z^\ydt \oplus \mf{F}(U)[t] \ar[r]^-{(\xi^\cdt, \hspace{.6pt}\nu^\cdt)} &\mf{F}(M)\ar[r]^{\theta^\ydt} & \mf{F}(N)[1]}
$$
in $D^b(\ModbLa)$. If $t> 0$, then $\nu^\ydt=0$ by Lemma \ref{Equiv-lemma}(2). Since ${\zeta^\ydt \choose \mu^\cdt}\vspace{1pt}$ is a pseudo-kernel of $(\xi^\cdt, 0)$, there exists some $u^\cdt: \mf{F}(U)[t] \to \mf{F}(N)$ such that $u^\cdt\circ \zeta^\ydt=0$ and $u^\cdt\circ \mu^\cdt=\id_{\mf{F}(U)[t]}.$ However, $u^\cdt=0$ by Lemma \ref{Equiv-lemma}(2), a contradiction. Thus $t\le 0$. A dual argument shows that $t\ge 0$. This establishes our claim.

Now, we deduce from Lemma \ref{Equiv-lemma}(1) that $Y^\ydt\cong \mf{F}(L)$ for some $L\in {\rm Rep}^-(Q^{\rm \hspace{.5pt}op})$. Since $\mf{F}_{N, L}$ is surjective, $\phi^\pdt=\mf{F}(f)$ for some $f: N\to L$ in ${\rm Rep}^-(Q^{\rm \hspace{.5pt}op})$. Consider an exact
triangle
$\xymatrix{N\ar[r]^{f} & L\ar[r]^{g} & V \ar[r]^h & N[1]}$ in $D^b({\rm Rep}^-(Q^{\rm \hspace{.5pt}op}))$. Applying $\mf{F}$ to this triangle yields a commutative diagram in
$D^b(\ModbLa)$ as follows:

$$
\xymatrix{
\mf{F}(N)\ar[r]^-{\phi^\pdt}\ar@{=}[d] & Y^\pdt \ar[d]_\cong \ar[r]^-{\psi^\pdt} & \mf{F}(M)\ar[d]^{\eta^\cdt}_\cong \ar[r]^-{\theta^\cdt} & \mf{F}(N)[1]\ar@{=}[d]\\
\mf{F}(N)\ar[r]^{\mf{F}(f)} & \mf{F}(L)\ar[r]^-{\mf{F}(g)} &
\mf{F}(V)\ar[r]^-{\mf{F}(h)} & \mf{F}(N)[1].}$$

Since $\mf{F}_{M, V}$ is surjective, $\eta^\ydt=\mf{F}(v)$ for some $v: M\to V$ in ${\rm Rep}^-(Q^{\rm \hspace{.5pt}op})$, and consequently, $\theta^\pdt=\mf{F}(hv)$. This shows that $\mf{F}_{M, N[1]}$ is surjective. If $\mf{F}(h)=0$, then $\mf{F}(g)$ is a retraction. Since $\mf{F}_{L, V}$ and $\mf{F}_{V, V}$ are isomorphisms, $g$ is a retraction, and hence, $h=0$. That is, $\mf{F}_{V, N[1]}$ is injective, and so is $\mf{F}_{M, N[1]}$. The proof of the theorem is completed.

\medskip

In the sequel, the triangle-equivalences $\mf{F}: D^b({\rm Rep}^-(Q^{\rm \hspace{.5pt}op}))\to D^b(\ModbLa)$ and $\mf{F}: D^b({\rm rep}^-(Q^{\rm \hspace{.5pt}op}))\to D^b({\rm mod}\hspace{.6pt}^b\hspace{-2pt}\La)$ will be called {\it Koszul equivalences}.

\smallskip

\section{\sc Galois Covering in the general case}

\medskip

\noindent Throughout this section, let $\La$ stand for a connected elementary locally bounded $k$-category with radical squared zero. By Gabriel's theorem, we may assume that $\La=kQ/(kQ^+)^2,$ where $Q$ is a connected locally finite quiver. The main objective of this section is to construct Galois coverings $\mk{F}_\pi: D^b({\rm Rep}^-(\tilde{Q}^{\rm \hspace{.5pt}op}))\to D^b(\ModbLa)\vspace{.6pt}$ and $\mk{F}_\pi: D^b({\rm rep}^-(\tilde{Q}^{\rm \hspace{.5pt}op}))\to D^b({\rm mod}\hspace{.6pt}^b\hspace{-2pt}\La)$, where $\tilde{Q}$ is a gradable covering
of $Q$. In particular, the indecomposable objects of $D^b(\ModbLa)$ and those of $D^b({\rm mod}\hspace{.6pt}^b\hspace{-2pt}\La)$ can be described in terms of the indecomposable representations of ${\rm Rep}^-(\tilde{Q}^{\rm \hspace{.5pt}op})$ and those of ${\rm rep}^-(\tilde{Q}^{\rm \hspace{.5pt}op}),$ respectively.

\medskip

To start with, we fix a connected component $\tilde{Q}$ of the gradable quiver $Q^{\mathbb{Z}}$ containing a vertex of form $(a_0, 0)$ with $a_0\in Q_0$. Observe that $\tilde{Q}$ admits a natural grading $\tilde{Q}_0=\cup_{n\in \mathbb{Z}}\hspace{.5pt} \tilde{Q}^n,$ where $\tilde{Q}^n$ consists of the vertices of form $(a, n)$ with $a\in Q_0$. We shall denote by $r_{\hspace{-1.5pt}_Q}$ the grading period of $Q$, which is defined by $r_{\hspace{-1.5pt}_Q}=0$ if $Q$ is gradable; and otherwise, $r_{\hspace{-1.5pt}_Q}$ is the minimum of the positive degrees of the closed walks in $Q$; see \cite[(7.3)]{BaL}. The {\it translation} $\rho$ of $\tilde{Q}$ is the automorphism sending a vertex $(a, n)$ to $(a, n+r_{\hspace{-1.5pt}_Q})$, and the {\it translation group} $G$ of $\tilde{Q}$ is the torsion-free group generated by $\rho$; see \cite[(7.4)]{BaL}. It is evident that $\rho$ is trivial if and only if $Q$ is gradable. It is known that there exists a Galois $G$-covering $\pi: \tilde{Q}\to Q$, sending a vertex $(a, n)$ to $a$, which we shall call a {\it minimal gradable covering} of $Q$. Observe that $Q$ is gradable if and only if $\pi$ is an isomorphism; see \cite[(7.5)]{BaL}. For convenience, we shall also call $\tilde{Q}$ a {\it minimal gradable covering} of $Q$.

\medskip

\begin{Lemma}\label{Quiver-cov}

Let $Q$ be a connected locally finite quiver with  a minimal gradable covering $\pi: \tilde{Q}\to Q$.

\vspace{-1pt}

\begin{enumerate}[$(1)$]

\item If $x\in \tilde{Q}^m$ and $y\in \tilde{Q}^n$ with $\pi(x)=\pi(y)$, then $n\equiv m \hspace{.7pt} (\hspace{.4pt}{\rm mod}\, r_{\hspace{-1.5pt}_Q})$.

\vspace{1pt}

\item If $x, y\in \tilde{Q}^n$ for some integer $n$, then $\pi(x)=\pi(y)$ if and only if $x=y$.

\vspace{1.4pt}

\item If $(a, 0)\in \tilde{Q}$ with $a\in Q_0$, then $(b, n)\in \tilde{Q}$ with $b\in Q_0$ and $n\in \Z$ if and only if $n\equiv d \hspace{.7pt} (\hspace{.4pt}{\rm mod}\, r_{\hspace{-1.5pt}_Q})$, where $d$ is the degree of some walk in $Q$ from $a$ to $b$.

\vspace{-2pt}

\end{enumerate}

\end{Lemma}

\smallskip

\noindent{\it Proof.} Assume that $\pi(x)=\pi(y)$, for some $x\in\tilde{Q}^m$ and $y\in \tilde{Q}^n$. Since $\pi$ is a Galois $G$-covering, $y=g \cdot x\vspace{1pt}$ for some $g\in G$. Since $g=\rho^s$ for some integer $s$, we see that $y\in \tilde{Q}^n\cap \tilde{Q}^{m+sr_{\hspace{-1.5pt}_Q}}.$ Hence, $n=m+sr_{\hspace{-1.5pt}_Q}$. This establishes Statement (1), and Statement (2) follows trivially from the definition of $\pi$.

Let $a, b\in Q_0$ with $(a, 0)\in \tilde{Q}.$ Being connected, $Q$ has a walk from $a$ to $b$, say, of degree $d$. Then, $(b, d)\in \tilde{Q}$; see \cite[(7.2)]{BaL}, that is $(b, d)\in \tilde{Q}^d$. If $(b, n)\in \tilde{Q}$, then it belongs to $\tilde{Q}^n$, and hence, $n\equiv d \hspace{.7pt} (\hspace{.4pt}{\rm mod}\, r_{\hspace{-1.5pt}_Q})$ by Statement (1). Conversely, assume that $n=d+s r_{\hspace{-1pt}_Q}$ for some $s\in \Z$, then $(b, n)=\rho^s\cdot (a, d)\in \tilde{Q}$. The proof of the lemma is the lemma is completed.

\medskip

Since $Q\vspace{.6pt}$ is locally finite, so is $\tilde{Q}$. In particular, we have a locally bounded $k$-category $\tLa=k\tilde{Q}/(k\tilde{Q}^+)^2\vspace{1pt}$ with radical squared zero. Since $\tilde{Q}$ is gradable, we have a Koszul functor $F: {\rm Rep}^-(\tilde{Q}^{\rm \hspace{.5pt}op})\to RC^{-,b}({\rm Proj}\,\tLa)$. In view of the commutative diagram (3.1), $F$ induces a commutative diagram

\vspace{-5pt}

$\begin{array}{lllll}
\hspace{-20pt}\begin{array}{c}  \\ \\ \\
 (4.1) \end{array}  &&
\xymatrix{C^b({\rm Rep}^-(\tilde{Q}^{\rm \hspace{.5pt}op})) \ar[r]^{\mathcal{P}_{\hspace{-.4pt}\tilde{Q}^{\rm \hspace{.5pt}op}}} \ar[d]^-{\mathcal{F}} &  K^b({\rm Rep}^-(\tilde{Q}^{\rm \hspace{.5pt}op})) \ar[r]^{\mathcal{L}_{\tilde{Q}^{\rm \hspace{.5pt}op}}} \ar[d]^-{\mathscr{F}} & D^b({\rm Rep}^{-}(\tilde{Q}^{\rm \hspace{.5pt}op})) \ar[d]^-{\mf{F}} \ar[d]\\
C^{-, b}({\rm Proj}\,\tLa) \ar[r]^{\mathcal{P}_{\hspace{.5pt}\tilde{\hspace{-1pt}\it\Lambda}}} & K^{-, b}({\rm Proj}\,\tLa) \ar[r]^{\mathcal{E}_{\tilde{\hspace{-1pt}\it\Lambda}}} & D^b(\Mod^{\hspace{.5pt}b\hspace{-1.5pt}}\tLa),}\end{array}$

\medskip

\noindent where $\mathscr{F}$ is triangle-exact and $\mf{F}$ is a triangle-equivalences; see (\ref{KzEqv}).

\medskip

Now, the $G$-action on $\tilde{Q}$ induces a $G$-action on $\tilde{Q}^{\rm \hspace{.5pt}op}$ such that $g\cdot \alpha^{\rm o}=(g\cdot \alpha)^{\rm o}$, for $\alpha\in \tilde{Q}_1$. In particular, $\rho\cdot (\tilde{Q}^{\rm \hspace{.5pt}op})^n=(\tilde{Q}^{\rm \hspace{.5pt}op})^{n-r_{\hspace{-1.2pt}_Q}},\vspace{1pt}$ for all $n\in \Z.$ Moreover, the $G$-action on $\tilde{Q}^{\rm \hspace{.5pt}op}$ induces a $G$-action on ${\rm Rep}(\tilde{Q}^{\rm \hspace{.5pt}op})$ as follows; see \cite[(3.2)]{Gab}. 

\vspace{-1pt}

\begin{enumerate}[$(1)$]

\item For an object $M\in {\rm Rep}(\tilde{Q}^{\rm \hspace{.5pt}op}),\vspace{.6pt}$ one defines $\rho\cdot M$ by $(\rho\cdot M)(x)=M(\rho^{-1}\cdot x),$ for $x\in \tilde{Q}_0;$ and $(\rho\cdot M)(\alpha^{\rm o})= M(\rho^{-1}\cdot \alpha^{\rm o}),$ for $\alpha\in \tilde{Q}_1.$

\vspace{1pt}

\item For a morphism $f: M\to N$ in ${\rm Rep}(\tilde{Q}^{\hspace{.5pt}\rm op}),\vspace{1pt}$ one defines $\rho\cdot f: \rho \cdot M \to \rho \cdot N$ by setting $(\rho\cdot f)(x)=f(\rho^{-1} \cdot x),$ for $x\in \tilde{Q}_0.$

\end{enumerate}

\medskip

Clearly, ${\rm Rep}^-(\tilde{Q}^{\hspace{.5pt}\rm op})$ is stable under this $G$-action on ${\rm Rep}(\tilde{Q}^{\hspace{.5pt}\rm op}).\vspace{1pt}$ Thus, we have a $G$-action on ${\rm Rep}^-(\tilde{Q}^{\hspace{.5pt}\rm op})$, which induces a $G$-action on each of $C^b({\rm Rep}^-(\tilde{Q}^{\hspace{.5pt}\rm op}))$, $K^b({\rm Rep}^-(\tilde{Q}^{\hspace{.5pt}\rm op}))$ and $D^b({\rm Rep}^-(\tilde{Q}^{\hspace{.5pt}\rm op}))$; see \cite[(5.4)]{BaL}.
%

\medskip

We shall need another group. Indeed, regarding $\rho$ as an automorphism of $C^b({\rm Rep}^{-}(\tilde{Q}^{\rm \hspace{.5pt}op}))$,
we obtain an automorphism $\vartheta=\rho\circ [-r_{\hspace{-1.5pt}_Q}]$, called {\it shifted translation}, of $C^b({\rm Rep}^{-}(\tilde{Q}^{\rm \hspace{.5pt}op}))$. Observe that $\vartheta$ is trivial if $r_{\hspace{-1.5pt}_Q}=0$ and of infinite order otherwise. Thus,  $\vartheta$ generates a torsion-free group $\mathfrak{G}$ of automorphisms, called {\it shifted translation group}, of $C^b({\rm Rep}^{-}(\tilde{Q}^{\rm \hspace{.5pt}op}))$. As usual, the $\mf{G}$-action on $C^b({\rm Rep}^{-}(\tilde{Q}^{\rm \hspace{.5pt}op}))$ induces a $\mf{G}$-action on each of $K^b({\rm Rep}^{-}(\tilde{Q}^{\rm \hspace{.5pt}op}))$ and $D^b({\rm Rep}^{-}(\tilde{Q}^{\rm \hspace{.5pt}op}))$; see \cite[(5.4)]{BaL}

\medskip

\begin{Lemma}\label{Der-gr}


Keep the notation introduced above.

\begin{enumerate}[$(1)$]

\item Let $M, N\in {\rm Rep}^{-}(\tilde{Q}^{\rm \hspace{.5pt}op})$ be such that $\Hom_{D^b({\rm Rep}^{-}(\tilde{Q}^{\rm \hspace{.5pt}op}))}(M,\, \vartheta^s \hspace{-1.5pt} \cdot \hspace{-1.5pt} N)\ne 0\vspace{1.5pt}$ for some $s\in \Z$. If
    $sr_{\hspace{-1.5pt}_Q}\ne 0$, then $s r_{\hspace{-1.5pt}_Q}=-1$ and $\Hom_{D^b({\rm Rep}^{-}(\tilde{Q}^{\rm \hspace{.5pt}op}))}(\vartheta^s \hspace{-1.5pt} \cdot \hspace{-1.5pt} N, M)=0.$

\vspace{1pt}

\item The $\mf{G}$-action on $D^b({\rm Rep}^{-}(\tilde{Q}^{\rm \hspace{.5pt}op}))\vspace{1.5pt}$ is free, locally bounded, and directed.

\end{enumerate} \end{Lemma}

\noindent{\it Proof.} (1) By hypothesis, $\Hom_{D^b({\rm Rep}^{-}(\tilde{Q}^{\rm \hspace{.5pt}op}))}(M, (\rho^s \hspace{-1.5pt} \cdot \hspace{-1.5pt} N)[-s\hspace{.55pt}r_{\hspace{-1.5pt}_Q}])\ne 0.$ If $sr_{\hspace{-1.5pt}_Q}\ne 0$, since ${\rm Rep}^{-}(\tilde{Q}^{\rm \hspace{.5pt}op})$ is hereditary, $s\hspace{.55pt}r_{\hspace{-1.5pt}_Q}=-1$, that is, $r_{\hspace{-1.5pt}_Q}=1$ and $s=-1$. Then,
$$\Hom_{D^b({\rm Rep}^{-}(\tilde{Q}^{\rm \hspace{.5pt}op}))}(\vartheta^s \hspace{-1.5pt} \cdot \hspace{-1.5pt} N, M)\cong
\Hom_{D^b({\rm Rep}^{-}(\tilde{Q}^{\rm \hspace{.5pt}op}))}(\rho^{-1} \hspace{-1.5pt} \cdot \hspace{-1.5pt} N,  M[-1])=0.$$

(2) We need only to consider the case where $r_{\hspace{-1.5pt}_Q}>0$. Let $M^\cdt, N^\ydt\in D^b({\rm Rep}^{-}(\tilde{Q}^{\rm \hspace{.5pt}op})).$ Since every object of $D^b({\rm Rep}^{-}(\tilde{Q}^{\rm \hspace{.5pt}op}))$ is a finite sum of stalk complexes; see, for example, \cite[(3.1)]{Len}, we may assume that $M^\cdt=M$ and $N^\cdt=N[s]$, where $M, N\in {\rm Rep}^{-}(\tilde{Q}^{\rm \hspace{.5pt}op})$ and $s\in \Z$. For each $t\in \Z$, we have
$$\Hom_{D^b({\rm Rep}^{-}(\tilde{Q}^{\rm \hspace{.5pt}op}))}(M^\cdt, \vartheta^t\cdot N^\cdt)=
\Hom_{D^b({\rm Rep}^{-}(\tilde{Q}^{\rm \hspace{.5pt}op}))}(M, (\rho^t \cdot N)[s-t \hspace{.5pt} r_{\hspace{-1.5pt}_Q}]).$$

If $\Hom_{D^b({\rm Rep}^{-}(\tilde{Q}^{\rm \hspace{.5pt}op}))}(M^\cdt, \vartheta^t\cdot N^\cdt)\ne 0,$ then
$s-t \hspace{.5pt} r_{\hspace{-1.5pt}_Q} \in \{ 0, 1\}$. Since $r_{\hspace{-1.5pt}_Q}\ne 0$, the number of such integers $t$ is at most two. Thus, the $\mf{G}$-action on
$D^b({\rm Rep}^{-}(\tilde{Q}^{\rm \hspace{.5pt}op}))\vspace{1.5pt}$ is locally bounded; and since $\mk{G}$ is torsion-free, it is also free; \vspace{1pt} see \cite[(2.2)]{BaL}.

If both $\Hom_{D^b({\rm Rep}^{-}(\tilde{Q}^{\rm \hspace{.5pt}op}))}(M^\cdt, \vartheta^t \hspace{-2pt}\cdot \hspace{-2pt} N^\cdt)\vspace{1.5pt}$ and
$\Hom_{D^b({\rm Rep}^{-}(\tilde{Q}^{\rm \hspace{.5pt}op}))}(\vartheta^t \hspace{-2pt} \cdot \hspace{-2pt} N^\cdt, M^\cdt)$ are non-zero, then $s-t\hspace{.5pt} r_{\hspace{-1.5pt}_Q}, t\hspace{.5pt} r_{\hspace{-1.5pt}_Q}-s\in \{0, 1\}$. Since $r_{\hspace{-1.5pt}_Q}\ne 0$, the number of such integers $t$ is at most one. That is, the $\mf{G}$-action on $D^b({\rm Rep}^{-}(\tilde{Q}^{\rm \hspace{.5pt}op}))\vspace{1.5pt}$ is directed. The proof of the lemma is completed.

\medskip

On the other hand, the $G$-action on $\tilde{Q}$ induces naturally a $G$-action on $\tLa,$ which in turn
induces a $G$-action on $\Mod\tLa\vspace{.5pt}$ as follows; see \cite[(3.2)]{Gab}. Fix an element $g\in G$, which will be regarded as an automorphism of $\tLa$.

\begin{enumerate}[$(1)$]

\item For a $\tLa$-module $M: \tLa\to \Mod\hspace{.5pt} k,$ one puts $g\cdot M=M\circ g^{-1}: \tLa\to \Mod\hspace{.5pt} k.$

\vspace{1.5pt}

\item For a $\tLa$-linear morphism $u: M\to N$, a $\tLa$-linear morphism $g \cdot u: g\cdot M\to g\cdot N$ is defined by $(g\cdot u)(x)=u(g^{-1}\cdot x),\vspace{1pt}$ for $x\in \tilde{Q}_0$. \vspace{-2pt}

\end{enumerate}

%

\medskip

Furthermore, the $G$-action on $\Mod\tLa$ induces a $G$-action on each of $C(\Mod\hspace{.6pt}\tLa)$, $K(\Mod\hspace{.6pt}\tLa)$ and $D(\Mod\hspace{.6pt}\tLa);$
see \cite[(5.4)]{BaL}. Clearly, $RC^{-,b}({\rm Proj}\,\tLa)$ is stable under the $G$-action on $C(\Mod\hspace{.4pt}\tLa).$ To study the behavior of $F: {\rm Rep}^-(\tilde{Q}^{\rm \hspace{.5pt}op})\to RC^{-,b}({\rm Proj}\,\tLa)$ with respect to these $G$-actions, we shall need some general considerations. Given a complex $X^\cdt$ over an additive category $\cA,$ we shall denote by $\mk{t}(X^\cdt)$ the complex whose $n$-th component is $X^n$ and whose $n$-th differential is $-d_X^n$, for $n\in \Z$; and for a morphism of complexes $f^\cdt=(f^n)_{n\in \mathbb{Z}}: X^\cdt\to Y^\ydt$, we have a complex morphism $\mk{t}(f^\cdt)=(f^n)_{n\in \mathbb{Z}}: \mk{t}(X^\cdt)\to \mk{t}(Y^\cdt)$. This yields an automorphism $\mk{t}$, called {\it twist functor}, of $C(\cA)$. Observe,
for each $p\in \Z$, that there exists a functorial isomorphism $\kappa_p: \mk{t}^{\hspace{.4pt}p} \to \id_{\,C(\mathcal{A})}$ defined by $\kappa_{\hspace{-1pt}_{p, X^\cdt}}=((-\id_{X^n})^{pn})_{n\in \mathbb{Z}}: \mk{t}^p(X^\cdt)\to X^\cdt,\vspace{1pt}$ for all $X^\cdt\in C(\cA)$. Recall that a linear functor $E: \cA\to \cB$ induces a linear functor $E^C: C(\cA)\to C(\cB)$ by component-wise action. The following statement is evident.

\medskip

\begin{Lemma}\label{twist}

Let $E: \cA\to \cB$ be a linear functor between additive categories, which induces a functor $E^C: C(\cA)\to C(\cB)$. If $X^\cdt\in C(\cA)$ and $p, q\in \Z$, then

\vspace{-1.5pt}

\begin{enumerate}[$(1)$]

\item $\kappa_{\hspace{-1.5pt}_{p, X^\cdt}} \circ \kappa_{\hspace{-1.5pt}_{q, \mf{t}^{^p\hspace{-1pt}}(X^\cdt)}}=\kappa_{\hspace{-1.5pt}_{p+q, X^\cdt}};$

\vspace{3pt}

\item $E^C(\kappa_{{\hspace{-1.5pt}_{p, X^\cdt}}})=\kappa_{\hspace{-1.5pt}_{p, E^C(X^\cdt)}}\,;$

\vspace{3pt}

\item $E^C(\mf{t}^p(X^\cdt))=\mf{t}^p(E^C(X^\cdt)).$

\end{enumerate}

\end{Lemma}

%

\medskip

\begin{Lemma}\label{FG-action}


Let $F: {\rm Rep}^-(\tilde{Q}^{\rm \hspace{.5pt}op})\to RC^{-,b}({\rm Proj}\,\tLa)\vspace{.5pt}$ be the Koszul functor, and let $\mathcal{F}: C^b({\rm Rep}^-(\tilde{Q}^{\rm \hspace{.5pt}op})) \to C^{-, b}({\rm Proj}\,\tLa)$ be the complex Koszul functor.

\begin{enumerate}[$(1)$]

\item  If $M\in {\rm Rep}^-(\tilde{Q}^{\rm \hspace{.5pt} op})$ and $s\in \Z,$ then
$\mf{t}^{\,s \hspace{.5pt} r_{\hspace{-1.5pt}_Q}}(\rho^s\cdot F(M)^\cdt)= F(\rho^s \hspace{-1.5pt} \cdot \hspace{-1.5pt} M)^\cdt \,[-s \hspace{.5pt} r_{\hspace{-1.5pt}_Q}].$


\vspace{2pt}

\item If $M^\cdt\in C^b({\rm Rep}^-(\tilde{Q}^{\rm \hspace{.5pt}op}))$ and $s\in \Z$, then
$\mk{t}^{\,s \hspace{.5pt} r_{\hspace{-1.5pt}_Q}}(\rho^s\cdot \mathcal{F}(M^\cdt))=\mathcal{F}(\vartheta^s \cdot M^\cdt).$

\end{enumerate}

\end{Lemma}

\noindent{\it Proof.} We need only to consider the case where $s=1.$ For simplicity, we write $r=r_{\hspace{-1.5pt}_Q}.\vspace{.5pt}$
Let $M\in {\rm Rep}^-(\tilde{Q}^{\rm \hspace{.5pt} op})$. Put $L^\cdt=\mk{t}^{\,r}(\rho\cdot F(M)^\cdt)$.
For each $n\in \mathbb{Z}$, since $\rho^{-1}\cdot \tilde{Q}^{\,r-n}=\tilde{Q}^{-n}$, we have
$$
\begin{array}{rcl}
F(\rho\hspace{-1.5pt} \cdot \hspace{-1.5pt} M)^n[-r]
&=& {\oplus}_{x\in \tilde{Q}^{\,r-n}} P[x]\otimes M(\rho^{-1} \cdot x) = {\oplus}_{y\in \tilde{Q}^{-n}} P[\rho \hspace{-1.5pt} \cdot \hspace{-1.5pt} y]\otimes M(y)\\ \vspace{-10pt}\\
&=& \rho \cdot ({\oplus}_{y\in \tilde{Q}^{-n}}  P[y]\otimes M(y) )= L^n, \\ \vspace{-10pt}\\
\end{array}$$
Moreover, $d^n_{F(\rho \cdot M) [-r]} = (-1)^r d_{F(\rho \cdot M)}^{n-r}=(-1)^r (d^{\hspace{.5pt} n-r}_{F(\rho \cdot M)}(y,x))_{(y, x)\in \tilde{Q}^{r-n-1}\times \tilde{Q}^{r-n}}. \vspace{2pt}$
For $(y, x)\in \tilde{Q}^{r-n-1}\times \tilde{Q}^{r-n},\vspace{2pt}$ we have $(y', x')=(\rho^{-1}\cdot y, \, \rho^{-1} \cdot x) \in \tilde{Q}^{n-1}\times \tilde{Q}^{-n}$, and \vspace{-10pt}

$$\hspace{-6pt}\begin{array}{rcl}
d^{\hspace{.5pt}n-r}_{F(\rho \cdot M)}(y,x)
\hspace{-8pt} &=& \hspace{-8pt} {\sum}_{\alpha\in \tilde{Q}_1(y, x)} \, P[\alpha]\otimes (\rho \cdot M)(\alpha^{\rm o})
={\sum}_{\alpha\in \tilde{Q}_1(y, x)}\, P[\alpha]\otimes M ((\rho^{-1} \cdot\alpha)^{\rm o})\\ \vspace{-9pt}\\
\hspace{-8pt} &=& \hspace{-8pt} {\sum}_{\beta\in \tilde{Q}_1(y', x')} \, P[\rho  \cdot \beta]\otimes M(\beta^{\rm o})
= \rho \cdot ({\sum}_{\beta\in \tilde{Q}_1(y', x')} \, P[\beta]\otimes M(\beta^{\rm o}) ) \\ \vspace{-9pt}\\
\hspace{-8pt} &=& \hspace{-8pt} \rho \cdot d_{F(M)}^n(y', x')
= d_{\rho \hspace{.7pt}\cdot F(M)}^n(y', x'),
\end{array}$$
where the third equation follows because $\rho \cdot \tilde{Q}_1(y', x')=\tilde{Q}_1(y, x).\vspace{1pt}$ This shows that $d^n_{F(\rho \cdot M)^\cdot [-r]}=(-1)^rd_{\rho \cdot F(M)}^n=d_L^n.\vspace{1pt}$ Hence, $F(\rho \cdot  M)^\cdt \,[-r]=\mk{t}^{\,r}(\rho\cdot F(M)^\cdt).\vspace{1pt}$ 

Next, we claim that $F(\rho \cdot f)^{n-r}= \rho\cdot F(f)^n$, for any morphism $f: M\to N$ in ${\rm Rep}^-(\tilde{Q}^{\rm \,op})$ and $n\in \Z$. 
Indeed,
$$\begin{array}{rcl}
F(\rho \cdot f)^{n-r} & = & \oplus_{x\in \tilde{Q}^{r-n}}\, \id_{P[x]}\otimes M(\rho^{-1} \cdot x) 
= \oplus_{y\in \tilde{Q}^{-n}}\, \id_{P[\rho\cdot y]}\otimes M(y)\\ \vspace{-8pt}\\
& = & \rho\cdot (\oplus_{y\in \tilde{Q}^{-n}}\, \id_{P[y]}\otimes M(y)) 
= \rho\cdot F(f)^n.
\end{array}$$

%
%
%
%
%
%
%
%
%
%


Now, let $M^\cdt\in C^b({\rm Rep}^-(\tilde{Q}^{\rm \hspace{.5pt}op})).$ For convenience, write $X^\cdt=\mk{t}^{\,r}(\rho\cdot \mathcal{F}(M^\cdt))$ and $Y^\ydt=\mathcal{F}(\vartheta\cdot  M^\cdt)=\mathcal{F}(\rho \cdot M^\cdt)[-r].$ Fix $n\in \Z$. By Statement (1), we obtain $\rho \cdot F(M^i)^{n-i}=F(\rho \cdot M^i)^{n-r-i},$ for $i\in \Z$. This yields
$$X^n=\oplus_{i\in \mathbb{Z}}\, \rho \cdot F(M^i)^{n-i}=\oplus_{i\in \mathbb{Z}}\, F(\rho \cdot M^i)^{n-r-i}=\mathcal{F}(\rho \cdot M^\cdt)^{n-r} = Y^n.$$
Moreover, we have
$$d_X^n=(d^n_X(j,i))_{(j,i)\in \Z\times \mathbb{Z}}: \oplus_{i\in \mathbb{Z}}\, \rho \cdot F(M^i)^{n-i} \to
\oplus_{j\in \mathbb{Z}}\, \rho \cdot F(M^j)^{n+1-j}
\vspace{2pt}$$
where $d^n_X(j,i): \rho \cdot F(M^i)^{n-i}\to \rho \cdot F(M^j)^{n+1-j}$ is given by

$$d^n_X(j,i)=(-1)^r \rho\cdot d^n_{\mathcal{F}(M^\ydt)}(j,i) = \left\{\begin{array}{ll}
(-1)^{r+i} \rho\cdot d_{F(M^i)}^{\hspace{.4pt}n-i}, & j=i;   \\ \vspace{-9pt} \\
(-1)^r \rho\cdot F(d_M^i)^{n-i},             & j=i+1; \\ \vspace{-9pt} \\
0,                          & j\ne i, i+1.
\end{array}\right.
\vspace{5pt}$$
On the other hand, we have
$$d_Y^n=(d_Y^n(j,i))_{(j, i)\in \Z\times \Z}: \oplus_{i\in \mathbb{Z}} \, F(\rho\cdot M^i)^{n-r-i}
\to \oplus_{j\in \mathbb{Z}}\, F(\rho\cdot M^j)^{n-r+1-j},$$ where
$$d_Y^n(j,i)=(-1)^r d^{n-r}_{\cF(\rho\cdot M^\cdt)}(j,i)=\left\{\begin{array}{ll}
(-1)^{r+i} d_{F(\rho \hspace{.4pt} \cdot M^i)}^{\hspace{.4pt}n-r-i}, & j=i;   \\ \vspace{-9pt} \\
(-1)^r F(\rho \cdot d_M^i)^{n-r-i},             & j=i+1; \\ \vspace{-9pt} \\
0,                          & j\ne i, i+1.
\end{array}\right.
\vspace{5pt}$$

For each integer $i$, by Statement (1), we obtain $d_{F(\rho \hspace{.4pt} \cdot M^i)}^{\hspace{.4pt}n-r-i}=\rho\cdot d_{F(M^i)}^{\hspace{.4pt}n-i},\vspace{2pt}$ and by our claim, we have $\rho\cdot F(d_M^i)^{n-i}=F(\rho \cdot d_M^i)^{n-r-i}$. This shows that $d_X^n=d_Y^n$.
The proof of the lemma is completed.

\medskip

The Galois $G$-covering $\pi: \tilde{Q}\to Q$ induces naturally a Galois $G$-covering functor $\pi: \tLa\to \La$ with a trivial $G$-stabilizer; see \cite[(7.8)]{BaL}. With the covering functor $\pi$, one associates its {\it push-down} $\pi_\lambda: {\rm Mod}\hspace{0.8pt}\tLa\to \Mod\La;$ see \cite[(3.2)]{BoG}, with a $G$-stabilizer $\delta=(\delta_{\rho^n})_{n\in \mathbb{Z}};$
see \cite[(6.3)]{BaL}. Note that $\pi_\lambda$ induces a commutative diagram \vspace{-5pt}

$\begin{array}{llllllll}
\hspace{-20pt} \begin{array}{c}  \\ \\ \\
 (4.2) \end{array}  &&&&&&
\xymatrix{
C(\Mod\hspace{-.5pt}\tLa) \ar[d]_{\pi_\lambda^C} \ar[r]^{\mathcal{P}_{\hspace{1pt}\tilde{\hspace{-1pt}\it\Lambda}}}& K(\Mod\hspace{-.5pt}\tLa)\ar[r]^{\mathcal{L}_{\hspace{1pt}\tilde{\hspace{-1pt}\it\Lambda}}}\ar[d]_{\pi^K_\lambda}& D(\Mod\hspace{-.5pt}\tLa)\ar[d]^{\pi^D_\lambda}\\
C(\Mod\hspace{-.5pt}\La)\ar[r]^{\mathcal{P}_{\it\Lambda}}& K(\Mod\hspace{-.5pt}\La)
\ar[r]^{\mathcal{L}_{\it\Lambda}}& D(\Mod\hspace{-.5pt}\La),}
\end{array} $

\medskip

\noindent where $\pi^C_\lambda, \pi^K_\lambda, \pi^D_\lambda$ 
are $G$-stable with $\pi^C_\lambda$ being exact and $\pi^K_\lambda, \pi^D_\lambda$ being triangle-exact; see \cite[(5.1)]{BaL}. Applying
the Koszul functor $F: {\rm Rep}^-(\tilde{Q}^{\rm \hspace{.5pt}op}) \to RC^{-,b}({\rm Proj}\, \tLa)\vspace{1pt}$ followed by $\pi_\lambda^C$, we obtain a composite functor
called {\it Koszul push-down}.

\medskip

\begin{Lemma}\label{F-pi}


The Koszul push-down is a $k$-linear functor \vspace{-2pt}
$$\xymatrix{F_\pi: {\rm Rep}^-(\tilde{Q}^{\rm \hspace{.5pt}op}) \ar[r]& RC^{-,b}({\rm Proj}\hspace{.5pt}\La)
}\vspace{-2pt}$$ such that, given any $M\in {\rm Rep}^-(\tilde{Q}^{\rm \hspace{.5pt}op})$, we have

\begin{enumerate}[$(1)$]

\item $F_\pi(M)^\cdt\in RC^{-, b}({\rm proj}\hspace{.5pt}\La)$ if and only if $M\in {\rm rep}^-(\tilde{Q}^{\hspace{.4pt}\rm op});$

\vspace{2pt}

\item $F_\pi(M)^\cdt\in RC^{\hspace{.4pt}b}({\rm proj}\hspace{.5pt}\La)$ if and only if $M\in {\rm rep}^b(\tilde{Q}^{\hspace{.6pt}\rm op})$.

\end{enumerate}
\end{Lemma}

\noindent{\it Proof.} 
Let $M\in {\rm Rep}^-(\tilde{Q}^{\rm \hspace{.5pt}op})$. Then, $F(M)^\dt\hspace{.4pt}\in RC^{-,b}({\rm Proj}\,\tLa),$ whose image under $\pi_\lambda^C\vspace{1pt}$ is $F_\pi(M)^\dt\,.$ Observe that $\pi_\lambda(P[x])=P[\pi_\lambda(x)]\vspace{1pt}$ for all $x\in \tilde{Q}_0;$ see \cite[(6.3)]{BaL}. For each $n\in \mathbb{Z}$, by the definition of $\pi_\lambda$; see \cite[Section 6]{BaL}, we have $F_\pi(M)^n=\oplus_{x\in \tilde{Q}^{-n}}\, P[\pi(x)]\otimes M(x)\in {\rm Proj}\hspace{.5pt}\La;$ and $d_{F_\pi(M)}^n: F_\pi(M)^n \to F_\pi(M)^{n+1}\vspace{2pt}$
is the $\La$-linear morphism given by the matrix
$(d_{F_\pi(M)}^{\hspace{.6pt}n}(y,x))_{(y,x) \in \tilde{Q}^{-n-1}\times \tilde{Q}^{-n}},$
where
$$
d^n_{F_\pi(M)}(y,x)=\textstyle{\sum}_{\alpha\in \tilde{Q}_1(y,x)}P[\pi(\alpha)]\otimes M(\alpha^{\rm o})\hspace{-1pt}: \hspace{-2pt} P[\pi(x)]\hspace{-2pt} \otimes \hspace{-2pt} M(x)\to P[\pi(y)]\hspace{-2pt} \otimes \hspace{-2pt} M(y),
$$
which is a radical morphism. Moreover, since $\pi_\lambda$ is exact; see \cite[(3.2)]{BoG}, $F_\pi(M)^\dt$ has bounded homology. Hence, $F_\pi(M)^\dt\in RC^{-,b}({\rm Proj}\,\La)$. This establishes the first part of the lemma. Moreover, Statements (1) and (2) follow from Lemma \ref{F-pi} and Proposition \ref{F-bch}.
The proof of the lemma is completed.

\medskip

The following statement describes morphisms between the images (up to shift) of the Koszul push-down.

\medskip

\begin{Lemma}\label{DPD-Morph}

Let $\hspace{-1pt}\xymatrixcolsep{16pt}\xymatrix{F_\pi: {\rm Rep}^-(\tilde{Q}^{\rm \hspace{.5pt}op}) \ar[r]& RC^{-,b}({\rm Proj}\hspace{.5pt}\La)}\hspace{-2pt}\vspace{-2pt}$ be the Koszul push-down, and let $\varphi^\cdt: F_\pi(M)^\cdt\to F_\pi(N)^\cdt[s]$ be a morphism with $M, N\in {\rm Rep}^-(\tilde{Q}^{\rm \hspace{.5pt}op})$ and $s\in \Z$.

\begin{enumerate}[$(1)$]

\item If $\varphi^\cdt$ is non-radical, then $s\equiv 0 \hspace{.5pt} ({\rm mod}\, r_{\hspace{-1.5pt}_Q})$.

\item If $\varphi^\cdt$ is non-zero and radical, then $s\equiv 1 ({\rm mod}\, r_{\hspace{-1.5pt}_Q}).$

\end{enumerate} \end{Lemma}

\noindent{\it Proof.} For each $n\in \Z$, by definition, $F_\pi(M)^n=\oplus_{x\in \tilde{Q}^{-n}}\, P[\pi(x)]\otimes M(x)\vspace{1pt}$ and $F_\pi(N)^{n+s}=\oplus_{y\in \tilde{Q}^{-n-s}}\, P[\pi(y)]\otimes N(y),\vspace{1pt}$ while $\varphi^n: F_\pi(M)^n\to F_\pi(N)^{n+s}$ can be written as a matrix
$\varphi^n=(\varphi^n(y, x))_{(y,x)\in \tilde{Q}^{-n-s} \times \tilde{Q}^{-n}},\vspace{1.5pt}$
with $\La$-linear morphisms $\varphi^n(y,x): P[\pi(x)]\otimes M(x) \to P[\pi(y)]\otimes N(y).$ By Lemma \ref{rqz-pm}, we have
$$ \varphi^n(y,x)={\textstyle\sum}_{\gamma\in Q_{\le 1}(\pi(y), \pi(x))}\, P[\gamma]\otimes f_\gamma, \; \mbox{ where } f_\gamma\in \Hom_k(M(x), N(y)).$$

Suppose that $\varphi^m\ne 0$ for some integer $m$. Then there exists a pair $(y, x)$ in $\tilde{Q}^{-m-s} \times \tilde{Q}^{-m}$ such that $\varphi^m(y,x)\ne 0.$ If $\varphi^\ydt$ is not radical, then we may assume that $\varphi^m(y,x)$ is not radical. In this case, $Q_{\le 1}(\pi(y), \pi(x))$ contains a trivial path, that is, $\pi(y)=\pi(x)$. By Lemma \ref{Quiver-cov}(1), $s\equiv 0 \,({\rm mod}\, r_{\hspace{-1.5pt}_Q})$.

Suppose that $\varphi^\cdt$ is radical. In particular, $\varphi^m(y,x)$ is radical. By Lemma \ref{rqz-pm}, $Q$ has an arrow $\gamma: \pi(y)\to \pi(x).$ Since $\pi: \tilde{Q}\to Q$ is a covering, $\gamma=\pi(\beta)$ for some arrow $\beta: y\to z$ in $\tilde{Q}$. In particular, $z\in \tilde{Q}^{-m-s+1}$ with $\pi(z)=\pi(x)$. By Lemma \ref{Quiver-cov}(1), we obtain $s\equiv 1 ({\rm mod}\, r_{\hspace{-1.5pt}_Q})$. The proof of the lemma is completed.

\medskip

The preceding result allows us to determine when two objects in the images (up to shift) of the Koszul push-down are isomorphic.

\medskip

\begin{Lemma}\label{hat-F-iso}

Let $\hspace{-1pt}\xymatrixcolsep{16pt}\xymatrix{F_\pi: {\rm Rep}^-(\tilde{Q}^{\rm \hspace{.5pt}op}) \ar[r]& RC^{-,b}({\rm Proj}\hspace{.5pt}\La)}\hspace{-1pt}\vspace{-2.5pt}$ be the Koszul push-down. If $M, N\in {\rm Rep}^-(\tilde{Q}^{\hspace{.5pt}\rm op})$ and $s\in \Z$, then

\begin{enumerate}[$(1)$]

\item $F_\pi(M)^\cdt \cong F_\pi(N)^\cdt$ if and only if $M\cong N,$ and more generally,

\vspace{1pt}

\item $F_\pi(M)^\cdt\cong F_\pi(N)^\cdt[s]$ if and only if $s=n r_{\hspace{-1.5pt}_Q}$ and $M\cong \rho^{n} \cdot N$ with $n\in \Z$.

\end{enumerate}

\end{Lemma}

\noindent{\it Proof.} (1) Let $M, N\in {\rm Rep}^-(\tilde{Q}^{\hspace{.5pt}\rm op})$. Assume that $\varphi^\pdt: F_\pi(M)^\cdt \to F_\pi(N)^\cdt$ is an isomorphism in $RC^{-,b}({\rm Proj}\hspace{.5pt}\La)$ with inverse $\psi^\cdt: F_\pi(N)^\cdt \to F_\pi(M)^\cdt$. For each $n\in \Z$, the isomorphism $\varphi^n$ can be written as a matrix as follows:
$$(\varphi^n(y,x))_{(y,x)\in \tilde{Q}^{-n}\times \tilde{Q}^{-n}}: \oplus_{x\in \tilde{Q}^{-n}}\, P[\pi(x)]\otimes M(x)\to \oplus_{y\in \tilde{Q}^{-n}}\, P[\pi(y)]\otimes N(y),$$
where $\varphi^n(y,x): P[\pi(x)]\otimes M(x) \to P[\pi(y)]\otimes N(y)$ is a $\La$-linear morphism. Write $$\varphi^n(x,x)=\id_{P[\pi(x)]}\otimes f_x+ {\textstyle\sum}_{\alpha\in Q_1(\pi(x), \pi(x))}\,  P[\alpha]\otimes f_\alpha,$$ where $f_x, f_\alpha\in \Hom_k(M(x), N(x))$. If $x\ne y$, then $\pi(x)\ne \pi(y)$ by Lemma \ref{Quiver-cov}(2), and hence, $\varphi^n(y,x)$ is a radical morphism. Writing $\psi^n$ in the same way and using the uniqueness stated in Lemma \ref{rqz-pm}, we conclude that $f_x$ is an isomorphism for every $x\in \tilde{Q}_0$. Let $\beta: z\to x$ be an arrow in $\tilde{Q}$ with $x\in \tilde{Q}^{-n}$. Using the equation $\varphi^{n+1} \circ d_{F_\pi(M)}^n = d_{F_\pi(N)}^n\circ \varphi^n\vspace{1pt}$ and the uniqueness stated in Lemma \ref{rqz-pm}, we deduce that $N(\beta^{\rm o})\circ f_x=f_z\circ M(\beta^{\rm o}).\vspace{1pt}$ This yields an isomorphism $f=(f_x)_{x\in \tilde{Q}_0}: M\to N\vspace{1pt}$ in ${\rm Rep}^-(\tilde{Q}^{\hspace{.5pt}\rm op})$.

(2) Suppose that $s=n r_{\hspace{-1.5pt}_Q}$ and $M\cong \rho^{n} \cdot N$ with $n\in \Z$. By Lemma \ref{FG-action}(1),
$F(M)^\cdt\cong \left(\rho^{n}\cdot F(N)^\cdt\right) [n r_{\hspace{-1.5pt}_Q}]=\rho^{n}\cdot F(N)^\cdt [s].\vspace{1pt}$ Applying the $G$-stable functor $\pi_\lambda^C$, we obtain $F_\pi(M)^\cdt \cong F_\pi(N)^\cdt[s]$.
Conversely, if $F_\pi(M)^\cdt\cong F_\pi(N)^\cdt[s]$ then, by Lemma \ref{DPD-Morph}(1), $s=n r_{\hspace{-1.5pt}_Q}$ for some $n\in \Z$. As we have shown, $F_\pi(N)^\cdt[n r_{\hspace{-1.5pt}_Q}] \cong F_\pi(\rho^n\cdot N)^\cdt,$ and hence, $F_\pi(M)^\cdt \cong F_\pi(\rho^n\cdot N)^\cdt.$ By Statement (1), we obtain $M\cong \rho^n\cdot N$. The proof of the lemma is completed.

\medskip

The Koszul push-down $\hspace{-1pt}\xymatrixcolsep{16pt}\xymatrix{F_\pi: {\rm Rep}^-(\tilde{Q}^{\rm \hspace{.5pt}op}) \ar[r]& RC^{-,b}({\rm Proj}\hspace{.5pt}\La)}\hspace{-1pt}\vspace{-2.5pt}$ can be extended to bounded complexes over ${\rm Rep}^-(\tilde{Q}^{\hspace{.5pt}\rm op}).\vspace{1pt}$ Indeed, since $\pi_\lambda^C$ is exact and sends projective modules to projective ones, composing the commutative diagrams (4.1) and (4.2),
we obtain a commutative diagram

\vspace{-5pt}

$\begin{array}{lllll}
\hspace{-20pt} \begin{array}{c}  \\ \\ \\
 (4.3) \end{array}  &&&
\xymatrix{
C^b({\rm Rep}^-(\tilde{Q}^{\rm \hspace{.5pt}op})) \ar[r]^{\mathcal{P}_{\hspace{-1pt}\tilde{Q}^{\rm \hspace{.5pt}op}}} \ar[d]^-{\mathcal{F}_\pi} &  K^b({\rm Rep}^-(\tilde{Q}^{\rm \hspace{.5pt}op})) \ar[r]^{\mathcal{L}_{\tilde{Q}^{\rm \hspace{.5pt}op}}} \ar[d]^-{\mathscr{F}_\pi} &
D^b({\rm Rep}^{-}(\tilde{Q}^{\rm \hspace{.5pt}op})) \ar[d]^-{\mf{F}_\pi} \\
C^{-, b}({\rm Proj}\hspace{.6pt}\La) \ar[r]^{\mathcal{P}_{_{\hspace{-1pt}\it\Lambda}}} & K^{-, b}({\rm Proj}\hspace{.6pt}\La) \ar[r]^{\mathcal{E}_{_{\hspace{-1pt}\it\Lambda}}}& D^b(\Mod^{\hspace{.5pt}b\hspace{-2pt}}\La),\vspace{5pt}
} \end{array}$

\medskip

\noindent where $\mathscr{F}_\pi$ and $\mf{F}_\pi$ are triangle-exact. In the sequel, we shall call $\mathcal{F}_\pi$ the {\it complex Koszul push-down} and $\mf{F}_\pi$ the {\it derived Koszul push-down} associated with the minimal gra\-dable covering $\pi: \tLa\to \La$. In order to collect some properties of $\mf{F}_\pi$, for each $x\in \tilde{Q}_0$, we denote by $I_{x^{\rm \hspace{.4pt}o}}$ the associated indecomposable injective representation of $\tilde{Q}^{\,\rm op}$ and by $S[\pi(x)]$ the simple $\La$-module supported by $\pi(x)$.

\medskip

\begin{Prop}\label{KPD-image}

Let $\xymatrixcolsep{16pt}\xymatrix{\mf{F}_\pi: D^b({\rm Rep}^{-}(\tilde{Q}^{\rm \hspace{.5pt}op})) \ar[r] & D^b(\Mod^{\hspace{.5pt}b\hspace{-2.5pt}}\La)}\vspace{-2.5pt}$ be the derived Koszul push-down, and let $M, N\in {\rm Rep}^-(\tilde{Q}^{\hspace{.5pt}\rm op})\vspace{-1.5pt}$ with $s, t\in \Z$.

\begin{enumerate}[$(1)$]

%

\item If $x\in \tilde{Q}^n$ for some integer $n$, then $\mf{F}_\pi(I_{x^{\rm \hspace{.4pt}o}})\cong S[\pi(x)][n].$

\vspace{1.5pt}

\item $\mf{F}_\pi(M)\cong \mf{F}_\pi(N)[s]\vspace{1pt}$ if and only if $s=n \hspace{.4pt} r_{\hspace{-1.5pt}_Q}$ and $M\cong \rho^n\cdot N$ for some $n\in \Z.$

\vspace{.5pt}

\item If $\Hom_{D^b({\rm Mod}^{\hspace{.4pt}b}\hspace{-1.5pt}\it\Lambda)}(\mf{F}_\pi(M)[t], \mf{F}_\pi(N)[s])\ne 0,\vspace{.5pt}$ then $s-t\equiv 0, 1\, ({\rm mod}\, r_{\hspace{-1.5pt}_Q}).$

\vspace{1.5pt}

\item If $u^\cdt: \mf{F}_\pi(M)[t]\to \mf{F}_\pi(N)[s]$ and $v^\cdt: \mf{F}_\pi(N)[s]\to \mf{F}_\pi(M)[t]$ are such that
$v^\cdt \hspace{.4pt} u^\cdt\ne 0$, then $s \equiv t \, ({\rm mod}\, r_{\hspace{-1.5pt}_Q}).$

\end{enumerate}

\end{Prop}

\noindent{\it Proof.} (1) We shall make use of the commutative diagrams (4.2) and (4.3). By Lemma \ref{hat-F-functor}, we obtain $\mf{F}_\pi(M)= \mathcal{P}_{_{\hspace{-1pt}\it\Lambda}}(\mathcal{E}_{_{\hspace{-1pt}\it\Lambda}}(F_\pi(M)^\cdt))\,.$ In particular, if $x\in \tilde{Q}^n$, by Lemma \ref{F-inj-image}, we have $\mf{F}_\pi(I_{x^{\rm \hspace{.4pt}o}})=\mathcal{P}_{_{\hspace{-1pt}\it\Lambda}}(\mathcal{E}_{_{\hspace{-1pt}\it\Lambda}}(\pi_\lambda^C(F(I_{x^{\rm o}})^\cdt)))
\cong \pi_\lambda^D(S[x][n]) \cong S[\pi(x)][n].$

\vspace{1pt}

(2) By Lemma \ref{F-pi}, $F_\pi(M)^\cdt, F_\pi(N)^\cdt\in RC^{-,b}({\rm Proj}\hspace{.6pt}\La)$. Since the projection functor $\mathcal{P}_{_{\hspace{-1pt}\it\Lambda}}: RC^{-,b}({\rm Proj}\hspace{.6pt}\La)\to K^{-,b}({\rm Proj}\hspace{.6pt}\La)$ reflects isomorphisms; see (\ref{Der-Hmtp}), we see that $\mf{F}_\pi(M)\cong \mf{F}_\pi(N)[s]$ if and only if $F_\pi(M)^\cdt\cong F_\pi(N)^\cdt[s].$  By Proposition \ref{hat-F-iso}(2), this is equivalent to the existence of an integer $n$ such that $s=n r_{\hspace{-1.5pt}_Q}$ and $M\cong \rho^n\cdot N.$

(3) Let $\Hom_{D^b({\rm Mod}^{\hspace{.4pt}b}\hspace{-1.5pt}\it\Lambda)}(\mf{F}_\pi(M)[t], \mf{F}_\pi(N)[s])\ne 0.\vspace{1pt}$
Since $\mathcal{P}_{_{\hspace{-1pt}\it\Lambda}}$ is full; see (\ref{Der-Hmtp}), we obtain $\Hom(F_\pi(M)^\cdt[t], F_\pi(N)^\cdt[s])\ne 0;$ and by Lemma \ref{DPD-Morph}, $s-t\equiv 0, 1\, ({\rm mod}\, r_{\hspace{-1.5pt}_Q})$.

\vspace{1pt}

(4) Let $u^\cdt: \mf{F}_\pi(M)[t]\to \mf{F}_\pi(N)[s]\vspace{1pt}$ and $v^\cdt: \mf{F}_\pi(N)[s]\to \mf{F}_\pi(M)[t]$ be non-zero morphisms in $D^b({\rm Mod}^{\hspace{.4pt}b}\hspace{-1.5pt}\it\Lambda).$ Then, $u^\cdt= \mathcal{E}_{_{\hspace{-1pt}\it\Lambda}}(\mathcal{P}_{_{\hspace{-1pt}\it\Lambda}}(f^\cdt))$ and $v^\cdt= \mathcal{E}_{_{\hspace{-1pt}\it\Lambda}}(\mathcal{P}_{_{\hspace{-1pt}\it\Lambda}}(g^\cdt)),$ for some non-zero morphisms $f^\pdt: F_\pi(M)^\cdt[s]\to F_\pi(N)^\cdt[t]$ and $g^\cdt: F_\pi(N)^\cdt[t]\to F_\pi(M)^\cdt[s]$ in $RC^{-,b}({\rm Proj}\,\La)$.
If $s\not\equiv t \, ({\rm mod}\, r_{\hspace{-1.5pt}_Q})$ then, by Lemma \ref{DPD-Morph}(1), $f^\cdt$ and $g^\cdt$ are radical morphisms. Since ${\rm rad}^2(\La)=0$, we obtain $g^\cdt\hspace{.4pt} f^\cdt=0$, and hence, $v^\cdt \hspace{.4pt} u^\cdt=0$. The proof of the proposition is completed.

\medskip

We shall need the following preparatory lemma.

\medskip

\begin{Lemma}\label{mfG-stable}

The derived Koszul push-down $\mf{F}_\pi: C^b({\rm Rep}^{-}(\tilde{Q}^{\hspace{.5pt}\rm op}))\to D^b(\Mod^{\hspace{.5pt}b\hspace{-2pt}}\La)$ is a $\mf{G}$-precovering, where $\mf{G}$ is the shifted translation group of $C^b({\rm Rep}^{-}(\tilde{Q}^{\hspace{.5pt}\rm op}))$.

%
%
%
%

\end{Lemma}

\noindent{\it Proof.} We shall make use of the commutative diagrams (4.1), (4.2) and (4.3). Write $r=r_{_Q}$ for the sake of simplicity. We first construct a $\mf{G}$-stabiliser for the complex Koszul push-down $\,\mathcal{F}_\pi: C^b({\rm Rep}^{-}(\tilde{Q}^{\hspace{.5pt}\rm op}))\to C^{-, b}({\rm Proj}\hspace{.6pt}\La)$. Recall that the push-down $\pi_\lambda: \Mod^{\hspace{.4pt}b\hspace{-1.5pt}}\tLa\to \ModbLa$ admits a $G$-stabiliser $\delta=(\delta_{\rho^p})_{p\in \mathbb{Z}};\vspace{1pt}$ see \cite[(6.3)]{BaL}, which induces a $G$-stabiliser $\delta^C=(\delta_{\rho^p}^C)_{p\in \mathbb{Z}}\vspace{1pt}$ for $\pi_\lambda^C: C(\Mod^{\hspace{.4pt}b\hspace{-1.5pt}}\tLa)\to C(\ModbLa);$ see \cite[(5.6)]{BaL}. Let $X^\cdt\in C^b(\Mod^{\hspace{.4pt}b\hspace{-1.5pt}}\tLa)$ and $p, q\in \Z$. By the definition of $\delta^C_{\rho^p};$ see \cite[(5.2)]{BaL}, we obtain
$\delta^C_{\rho^p, \mf{t}(X^\cdt)}=(\delta_{\rho^p, X^n})_{n\in \mathbb{Z}}$ and $\delta^C_{\rho^p, X^\cdt}=(\delta_{\rho^p, X^n})_{n\in \mathbb{Z}}.\vspace{1pt}$ Moreover, in view of Lemma \ref{twist}, we see that $\rho^q \cdot \mf{t}^p(X^\cdt)=\mf{t}^p(\rho^q\cdot X^\cdt)$, and $\pi_\lambda^C(\kappa_{p, \, X^\cdt})=\kappa_{p, \, \pi_\lambda^C(X^\cdt)},$ and $\pi_\lambda^C(\kappa_{p, \,\rho^q\cdot X^\cdt})=\kappa_{p, \, \pi_\lambda^C(\rho^q\cdot X^\cdt)}.$ This yields
$$
\begin{array}{rcl}
\pi_\lambda^C(\kappa_{p,\, X^\cdt})\circ \delta^C_{\rho^q, \mf{t}^p(X^\cdt)} & = &
\kappa_{p, \, \pi_\lambda^C(X^\cdt)} \circ \delta^C_{\rho^q, \mf{t}^p(X^\cdt)} 
= ((-\id_{\pi_\lambda(X^n)})^{pn} \circ \delta_{\rho^q, X^n})_{n\in \mathbb{Z}}\\ \vspace{-8pt}\\
& = & (\delta_{\rho^q, X^n}\circ (-\id_{\pi_\lambda(\rho^q \cdot X^n)})^{pn})_{n\in \mathbb{Z}} 
= \delta^C_{\rho^q, X^\cdt}  \circ  \kappa_{p, \, \pi_\lambda^C(\rho^q\cdot X^\cdt)} \\ \vspace{-8pt}\\
& = & \delta^C_{\rho^q, X^\cdt} \circ \pi_\lambda^C(\kappa_{p, \, \rho^q\cdot X^\cdt}).
\end{array} $$

\vspace{1pt}

Fix $M^\ydt\in C^b({\rm Rep}^{-}(\tilde{Q}^{\rm \hspace{.5pt}op}))$ and $p\in \Z$. In view of Lemma \ref{FG-action}(2), we see that $\mathcal{F}(\vartheta^p\cdot M^\cdt)= \mf{t}^{\,pr}(\rho^p\cdot \mathcal{F}( M^\cdt))$. As previously defined, we have a natural isomorphism
$\xymatrixcolsep{16pt}\xymatrix{\kappa_{pr, \, \rho^p\cdot \mathcal{F}(M^\cdt)}: \mathcal{F}(\vartheta^p\cdot M^\cdt)\ar[r] & \rho^p\cdot \mathcal{F}(M^\cdt).}$ Let $\hspace{-2pt}\xymatrixcolsep{16pt}\xymatrix{\zeta_{\vartheta^p, M^\cdt}^C: \mathcal{F}_\pi(\vartheta^p\cdot M^\cdt) \ar[r] & \mathcal{F}_\pi(M^\cdt)}\hspace{-2pt}$ be the composite of the following natural isomorphisms: \vspace{-2pt}
$$\xymatrixcolsep{18pt}\xymatrix{\pi_\lambda^D(\mf{t}^{pr}(\rho^s\cdot \mathcal{F}(M^\cdt)))\ar[rrr]^-{\pi_\lambda^C\left(\kappa_{pr, \,\rho^p\cdot \mathcal{F}(M^\cdt)}\right)} &&&
\pi_\lambda^C(\rho^p\cdot \mathcal{F}(M^\cdt)) \ar[rr]^-{\delta^C_{\rho^p, \mathcal{F}(M^\cdt)}} &&
\mathcal{F}_\pi(M^\cdt).} \vspace{-3pt}$$
This gives rise to a functorial isomorphism $\xymatrixcolsep{18pt}\xymatrix{\zeta_{\vartheta^p}^C: \mathcal{F}_\pi\circ \vartheta^p\ar[r] & \mathcal{F}_\pi.}$ Given an integer $q$, by definition,
$\xymatrixcolsep{18pt}\xymatrix{\zeta_{_{\vartheta^q, \vartheta^p\cdot M^\cdt}}^C: \mathcal{F}_\pi(\vartheta^q\cdot (\vartheta^p\cdot M^\cdt)) \ar[r] & \mathcal{F}_\pi(\vartheta^p\cdot M^\cdt)}$ is the compo\-site of $\hspace{-2pt}\xymatrixcolsep{18pt}\xymatrix{\pi_\lambda^C\left(\kappa_{qr, \, \rho^q\cdot \mathcal{F}(\vartheta^p\cdot M^\cdt)}\right) : \pi_\lambda^C(\mf{t}^{qr}(\rho^q\cdot \mathcal{F}(\vartheta^p\cdot M^\cdt)))\ar[r] &
\pi_\lambda^C(\rho^q\cdot \mathcal{F}(\vartheta^p\cdot M^\cdt))
}$ and $\xymatrixcolsep{18pt}\xymatrix{\delta^C_{\rho^q, \mathcal{F}(\vartheta^p\cdot M^\cdt)}: \pi_\lambda^C(\rho^q\cdot \mathcal{F}(\vartheta^p\cdot M^\cdt)) \ar[r] & \mathcal{F}_\pi(\vartheta^p\cdot M^\cdt).}$ Applying the equations stated-above, we deduce that \vspace{1pt}
$$\begin{array}{rcl}
\pi_\lambda^C\hspace{-3pt}\left(\kappa_{pr, \, \rho^p \cdot \mathcal{F}(M^\cdt)}\right) \circ \delta^C_{\rho^q, \mathcal{F}(\vartheta^p \cdot M^\cdt)}
& = & \pi_\lambda^C\hspace{-3pt}\left(\kappa_{pr, \, \rho^p\cdot \mathcal{F}(M^\cdt)}\right) \circ \delta^C_{\rho^q, \mf{t}^{pr}(\rho^p\cdot \mathcal{F}( M^\cdt))} \\ \vspace{-8pt}\\
& = &  \delta^C_{\rho^q, \rho^p\cdot \mathcal{F}(M^\cdt)} \circ
\pi_\lambda^C\hspace{-3pt}\left(\kappa_{pr, \,\rho^q\cdot(\rho^p\cdot \mathcal{F}(M^\cdt))}\right) \\ \vspace{-8pt}\\
& = & \delta^C_{\rho^q, \rho^p\cdot \mathcal{F}( M^\cdt)} \circ
\pi_\lambda^C\hspace{-3pt}\left(\kappa_{pr, \, \rho^{q+p} \cdot \mathcal{F}(M^\cdt))}\right).
\end{array}$$

This enables us to obtain \vspace{3pt}
$$\begin{array}{rcl}
&&\zeta^C_{\vartheta^p, M^\cdt} \hspace{-1pt}\circ \hspace{-1pt} \zeta^C_{\vartheta^q, \vartheta^p\cdot M^\cdt}\\ \vspace{-8pt}\\
\hspace{-8pt} & =  &
\delta_{\rho^p, \mathcal{F}(M^\cdt)}^C \hspace{-1pt}\circ \hspace{-1pt} \pi_\lambda^C\hspace{-3pt}\left(\kappa_{pr,\, \rho^p \cdot \mathcal{F}(M^\cdt)}\right) \hspace{-1pt}\circ \hspace{-1pt} \delta^C_{\rho^q, \mathcal{F}(\vartheta^p\cdot M^\cdt)} \hspace{-1pt}\circ \hspace{-1pt} \pi_\lambda^C\hspace{-3pt}\left(\kappa_{qr, \, \rho^q \cdot \mathcal{F}(\vartheta^p\cdot M^\cdt)}\right) \\ \vspace{-8pt}\\
& =  &
\delta_{\rho^p, \mathcal{F}(M^\cdt)}^C \circ \delta^C_{\rho^q, \rho^p\cdot \mathcal{F}(M^\cdt)} \circ
\pi_\lambda^C\hspace{-3pt}\left(\kappa_{pr, \,\rho^{q+p} \cdot \mathcal{F}(M^\cdt))}\right) \hspace{-1pt}\circ \hspace{-1pt} \pi_\lambda^C\hspace{-3pt}\left(\kappa_{qr,\, \mf{t}^{pr}(\rho^{q+p} \cdot \mathcal{F}(M^\cdt))}\right) \\ \vspace{-8pt}\\
& = & \delta_{\rho^{q+p}, \mathcal{F}(M^\cdt)}^C\hspace{-1pt}\circ \hspace{-1pt}
\pi_\lambda^C \hspace{-3pt}\left(\kappa_{(q+p)r, \, \mf{t}^{(q+p)r}(\rho^{q+p} \cdot \mathcal{F}(M^\cdt))} \right) = \zeta^C_{\vartheta^{q+p}, M^\cdt}\,.
\end{array}\vspace{2pt}$$
That is, $\zeta^C=(\zeta^C_{\vartheta^p})_{p\in \mathbb{Z}}$ is a $\mf{G}$-stabilizer for $\mathcal{F}_\pi$. Next, given a bounded complex $M^\cdt$ over ${\rm Rep}^{-}(\tilde{Q}^{\rm \hspace{.5pt}op})$ and an integer $p$, we set
$$\xymatrix{\zeta_{\rho^p, M^\cdt}^D=\mathcal{P}_{_{\hspace{-1pt}\it\Lambda}}(\mathcal{E}_{_{\hspace{-1pt}\it\Lambda}}(\zeta_{\rho^p, M^\cdt}^C)): \mf{F}_\pi(\vartheta^p\cdot M^\cdt) \ar[r] & \mf{F}_\pi(M^\cdt)}. \vspace{-1pt}$$
This yields a natural isomorphism $\zeta^D_{\rho^p}=\{\zeta_{\rho^p, M^\cdt}^D \mid M^\cdt\in D^b({\rm Rep}^{-}(\tilde{Q}^{\rm \hspace{.5pt}op}))\}$ from $\mf{F}_\pi\circ \rho^p$ to $\mf{F}_\pi.$ It is now easy to see that  $\zeta^D=(\zeta^D_{\rho^p})_{p\in \mathbb{Z}}$ is a $\mf{G}$-stabiliser for $\mf{F}_\pi.$

To conclude the proof, we fix two complexes $M^\cdt, N^\cdt$ in $D^b({\rm Rep}^-(\tilde{Q}^{\hspace{.5pt}\rm op}))\vspace{1pt}$. The Koszul equivalence $\mf{F}: D^b({\rm Rep}^{-}(\tilde{Q}^{\rm \hspace{.5pt}op})) \to D^b(\Mod^{\hspace{.5pt}b\hspace{-1.5pt}}\tLa)$ stated in Theorem \ref{KzEqv} induces an isomorphism
$$\mf{F}_{M^\cdt, N^\cdt}: \oplus_{p\in \Z}\Hom_{D^b({\rm Rep}^-(\tilde{Q}^{\rm \hspace{0.5pt} op}))}(M^\cdt,
\vartheta^p \hspace{-1pt} \cdot \hspace{-1pt} N^\cdt) \to \Hom_{D^b({\rm Mod}^b\tilde{\hspace{-2pt}\it\Lambda})}(\mf{F}(M^\cdt),
\mf{F}(\vartheta^p \hspace{-1pt} \cdot \hspace{-1pt} N^\cdt)),$$
sending $(f^\ydt_p)_{p\in \mathbb{Z}}$ to $(\mf{F}(f_p^\cdt))_{p\in \mathbb{Z}}$, where $f_p^\cdt\in \Hom_{D^b({\rm Rep}^-(\tilde{Q}^{\rm \hspace{0.5pt} op}))}(M^\cdt, \vartheta^p\hspace{-1pt} \cdot \hspace{-1pt} N^\cdt).\vspace{1.5pt}$

\vspace{1pt}

In view of Lemma \ref{FG-action}(2), we deduce that $\mf{F}(\vartheta^p\cdot N^\cdt)=\mf{t}^{pr}(\rho^p\cdot \mf{F}(N^\cdt)).$ Setting $\tilde{\kappa}_{pr, \,\rho^p\cdot \mf{F}(N^\cdt)}=\mathcal{P}_{_{\hspace{-1pt}\tilde{\it\Lambda\hspace{1pt}}}} (\mathcal{L}_{_{\hspace{-.5pt}\tilde{\it\Lambda\hspace{1pt}}}}(\kappa_{pr, \,\rho^p\cdot \mf{F}(N^\cdt)})$,
we obtain an isomorphism
$$\kappa_{_{\hspace{-1pt}M^\cdt\hspace{-1pt}, \hspace{-1pt} N^\cdt}}: \hspace{-2pt} \oplus_{p\in \Z} \Hom_{D^b({\rm Mod}^b\tilde{\hspace{-2pt}\it\Lambda})\hspace{-1pt}}(\mf{F}(M^\cdt),
\mf{F}(\vartheta^p\hspace{-1pt} \cdot \hspace{-1pt} N^\cdt))\hspace{-2pt}\to \hspace{-1pt}
\oplus_{p\in \Z} \Hom_{D^b({\rm Mod}^b\tilde{\hspace{-2pt}\it\Lambda})\hspace{-1pt}} (\mf{F}(M^\cdt),
\rho^p\cdot \mf{F}(N^\cdt)),$$
sending $(\mf{F}(f_p^\cdt))_{p\in \mathbb{Z}}\vspace{1pt}$ to $(\tilde{\kappa}_{pr, \,\rho^p\cdot \mf{F}(N^\cdt)}\circ \mf{F}(f_p^\cdt))_{p\in \mathbb{Z}}$.\vspace{1.5pt} Finally, considering the $G$-precovering $\pi_\lambda^D: D^b(\Mod\hspace{0.5pt}\tLa) \to D^b(\Mod \La)$; see \cite[(6.7)]{BaL}, we obtain an isomorphism
$$\pi_\lambda^D(M^\cdt, N^\cdt):
\oplus_{p\in \Z} \Hom_{D^b({\rm Mod}^b\tilde{\hspace{-2pt}\it\Lambda})\hspace{-1pt}} (\mf{F}(M^\cdt),
\rho^p\cdot \mf{F}(N^\cdt)) \to \Hom_{D^b({\rm Mod}^b\hspace{-2pt}\it\Lambda)\hspace{-1pt}}(\mf{F}(M^\cdt), \mf{F}(N^\cdt)),$$
sending $(\tilde{\kappa}_{pr, \, \rho^p\cdot \mf{F}(N^\cdt)}\circ \mf{F}(f_p^\cdt))_{p\in \mathbb{Z}}$ to ${\sum}_{p\in \mathbb{Z}}\, \delta^D_{\rho^p, \mf{F}(N^\cdt)}\circ \pi_\lambda^D(\tilde{\kappa}_{pr, \rho^p\cdot \mf{F}(N^\cdt)} \circ \mf{F}(f_p^\cdt)),$
where
$$\begin{array}{rcl}
&& \delta^D_{\rho^p, \mf{F}(N^\cdt)}\circ \pi_\lambda^D(\tilde{\kappa}_{pr, \rho^p\cdot \mf{F}(N^\cdt)} \circ \mf{F}(f_p^\cdt)) \\
\vspace{-8pt} \\
&=&\mathcal{P}_{_{\hspace{-1pt}\it\Lambda}}(\mathcal{E}_{_{\hspace{-1pt}\it\Lambda}}(\delta^C_{\rho^p, \mf{F}(N^\cdt)}))
\circ \pi_\lambda^D(\mathcal{P}_{_{\hspace{-1pt}\tilde{\it\Lambda\hspace{1pt}}}}
(\mathcal{L}_{_{\hspace{-.5pt}\tilde{\it\Lambda\hspace{1pt}}}}(\kappa_{pr, \,\rho^p\cdot \mf{F}(N^\cdt)})))
\circ \pi_\lambda^D(\mf{F}(f_p^\cdt)) \\
\vspace{-8pt} \\
&=& \mathcal{P}_{_{\hspace{-1pt}\it\Lambda}}(\mathcal{E}_{_{\hspace{-1pt}\it\Lambda}}(\delta_{\rho^p, \mf{F}(N^\cdt)}^C))
\circ \mathcal{P}_{_{\hspace{-1pt}\it\Lambda}}(\mathcal{E}_{_{\hspace{-1pt}\it\Lambda}}(\pi_\lambda^C(\kappa_{pr, \rho^p\cdot \mf{F}(N^\cdt)})))\circ \mf{F}_\pi(f_p^\cdt)\\
\vspace{-8pt} \\
&=& \mathcal{P}_{_{\hspace{-1pt}\it\Lambda}}(\mathcal{E}_{_{\hspace{-1pt}\it\Lambda}}(\zeta^C_{\vartheta^p, N^\cdt}))\circ \mf{F}_\pi(f_p^\cdt).
\end{array}$$

Composing the above three isomorphisms, we obtain a desired isomorphism
$$\mf{F}_\pi^{M^\cdt, N^\cdt}: \oplus_{p\in \Z}\Hom_{D^b({\rm Rep}^-(\tilde{Q}^{\rm \hspace{0.5pt} op}))}(M^\cdt,
\vartheta^p\hspace{-1pt} \cdot \hspace{-1pt} N^\cdt) \to \Hom_{D^b({\rm Mod}^b\hspace{-2pt}\it\Lambda)\hspace{-1pt}}(\mf{F}_\pi(M^\cdt), \mf{F}_\pi(N^\cdt)),$$ sending
$(f^\ydt_p)_{p\in \mathbb{Z}}$ to ${\sum}_{p\in \mathbb{Z}}\, \zeta^D_{\vartheta^p, N^\cdt}\circ \mf{F}_\pi(f_p^\ydt).\vspace{1.5pt}$
The proof of the lemma is completed.

\medskip

In order to state the main result of this section, we denote by ${\rm Ind}^{-}(\tilde{Q}^{\hspace{.5pt}\rm op})$ a complete set of isomorphism class representatives of the indecomposable objects of ${\rm Rep}^{-}(\tilde{Q}^{\hspace{.5pt}\rm op})$. Given an integer $r\ge 0$, moreover, we write $\mathbb{Z}\,_r=\mathbb{Z}\,$ if $r=0;$ and otherwise, $\mathbb{Z}\,_r=\{0, 1, \ldots r-1\}.$

\medskip

\begin{Theo} \label{Main-1}

Let $\La=kQ/(kQ^+)^2$ with $Q$ a connected locally finite quiver, and let $\pi: \tilde{Q}\to Q$ be a minimal gradable covering
of $Q$.

\begin{enumerate}[$(1)$]

\item The derived Koszul push-down $\mf{F}_\pi: D^b({\rm Rep}^{-}(\tilde{Q}^{\hspace{.5pt}\rm op}))\to D^b(\ModbLa) \vspace{1pt}$ is a Galois $\mf{G}$-covering, where $\mf{G}$ is the shifted translation group of $C^b({\rm Rep}^{-}(\tilde{Q}^{\hspace{.5pt}\rm op}))$.

\vspace{1pt}

\item The non-isomorphic indecomposable objects of $D^b\hspace{-1pt}(\ModbLa)$ are $\mf{F}_\pi(M)[s]$ with
$M\in {\rm Ind}^-(\tilde{Q}^{\hspace{.5pt}\rm op})$ and $s\in \mathbb{Z}_{\hspace{.4pt}r},$ where $r$ is the grading period of $Q$.

\vspace{1.5pt}

\item Every object $X^\cdt$ of $D^b\hspace{-.5pt}(\ModbLa)$ is uniquely $($up to isomorphism$)$ decomposed as \vspace{-1pt}
$$X^\cdt=\textstyle{\oplus}_{\hspace{.4pt} s\in \mathbb{Z}_{\hspace{.4pt}r}} \hspace{.4pt} \mf{F}_\pi(M_s)[s],\vspace{-1pt}$$ where the $M_s\in {\rm Rep}^-(\tilde{Q}^{\hspace{.5pt}\rm op}),$ all but finitely many of them are zero.

\end{enumerate}

\end{Theo}

\smallskip

\noindent{\it Proof.} First of all, Statement (2) will follow easily from Statements (1) and (3). By Lemma \ref{Der-gr}, the $\mf{G}$-action on $D^b({\rm Rep}^{-}(\tilde{Q}^{\rm \hspace{.5pt}op}))\vspace{1.5pt}$ is free, locally bounded and directed; and by Lemma \ref{mfG-stable}(2), $\mf{F}_\pi$ is a $\mf{G}$-precovering. Since the idempotents in $D^b({\rm Rep}^{-}(\tilde{Q}^{\rm \hspace{.5pt}op}))\vspace{1.5pt}$ split; see \cite{BaS}, it follows from Lemma 2.10 in \cite{BaL} that $\mf{F}_\pi$ satisfies Conditions (2) and (3) of the notion of a Galois covering stated in Subsection 1.3.

For the rest of proof, we shall make use of the commutative diagram (4.3). Let $X^\cdt\in D^b(\Mod^{\hspace{.5pt}b\hspace{-2pt}}\La).$ By
Proposition \ref{Der-Hmtp}(3), $X^\cdt\cong\vspace{1pt} \mathcal{E}_{\hspace{-1pt}_{\it\Lambda}}(\mathcal{P}_{_{\hspace{-1pt}\it\Lambda}}(P^\pdt)),$
for some radical complex $P^\pdt$ in $RC^{-, b}({\rm Proj}\hspace{.5pt}\La).$ Note that $P^\pdt\cong \pi_\lambda^C(L^\cdt)$ for some $L^\cdt\in RC^{-, b}({\rm Proj}\hspace{3pt}\tilde{\hspace{-2pt}\La});$ see \cite[(7.9)]{BaL}. Setting $N^\cdt= \mathcal{E}_{_{\tilde{\hspace{-1pt}\it\Lambda}}}(\mathcal{P}_{_{\tilde{\hspace{-1pt}\it\Lambda}}} (L^\pdt)),$ we obtain
$X^\cdt\cong \pi_\lambda^D(N^\cdt)$. Since the derived Koszul functor $\mf{F}: D^b({\rm Rep}^{-}(\tilde{Q}^{\hspace{.5pt}\rm op}))\to D^b({\rm Mod}\tLa)$ is an equivalence, $N^\cdt\cong \mf{F}(M^\cdt)$ for some complex $M^\cdt\in D^b({\rm Rep}^-(\tilde{Q}^{\hspace{.5pt}\rm op}))$. Thus, $X^\cdt\cong \mf{F}_\pi(M^\cdt)$. This shows that $\mf{F}_\pi$ is dense, and hence, it is a Galois $\mf{G}$-covering.


For proving Statement (3), we consider a complex $X^\cdt\in D^b\hspace{-.5pt}(\ModbLa)$. Suppose first that
$X^\cdt=\mf{F}_\pi(M)[t]$ with $M\in {\rm Rep}^-(\tilde{Q}^{\hspace{.5pt}\rm op})$ and $t\in \Z$. Writing
$t= n\, r + s$ with $n\in \Z$ and $s\in \Z_{\,r}\vspace{1pt}$ and setting $N=\rho^n\cdot M$, by Proposition \ref{KPD-image}(1), we obtain $X^\cdt\cong \mf{F}_\pi(N)[s]$. In general, since $\mf{F}_\pi$ is dense and every complex in $D^b({\rm Rep}^-(\tilde{Q}^{\hspace{.5pt}\rm op}))$ is a finite sum of stalk complexes, $X^\cdt$ admits a desired decomposition. For proving the uniqueness, assume that there exists an isomorphism
$$\xymatrix{f^\cdt: X^\cdt= \textstyle{\oplus}_{\hspace{.4pt} s\in \Z_{\,r}} \hspace{.4pt} \mf{F}_\pi(M_s)[s]\ar[r] & \textstyle{\oplus}_{\hspace{.4pt} t\in \Z_{\,r}} \hspace{.4pt} \mf{F}_\pi(N_t)[t]=Y^\cdt,}\vspace{-2pt}$$
whose inverse is $g^\cdt$. Consider the canonical projections $p_s^\cdt: X^\cdt\to \mf{F}_\pi(M_s)[s]$ and $u_t^\cdt: Y^\cdt\to \mf{F}_\pi(N_t)[t]$ and the canonical injections $q_s^\cdt: \mf{F}_\pi(M_s)[s]\to X^\cdt$ and $v_t^\cdt: \mf{F}_\pi(N_t)[t]\to Y^\cdt$. If $s, t\in \Z_{\hspace{.4pt}r}$ are distincts, then we deduce from Proposition \ref{KPD-image}(3) that $(p_s^\cdt v^\cdt v_t^\cdt)(u_t^\cdt u^\cdt q_s^\cdt)=0$ and $(u_s^\cdt u^\cdt q_t^\cdt)(p_t^\cdt v^\cdt v_s^\cdt)=0.$ As a consequence, $(p_s^\cdt v^\cdt v_s^\cdt)(u_s^\cdt u^\cdt q_s^\cdt)=1_{\mf{F}_\pi(M_s)[s]}$ and
$(u_s^\cdt u^\cdt q_s^\cdt)(p_s^\cdt v^\cdt v_s^\cdt)=1_{\mf{F}_\pi(N_s)[s]}.\vspace{1pt}$ That is, $\mf{F}_\pi(M_s)[s]\cong \mf{F}_\pi(N_s)[s].$ Applying Proposition \ref{hat-F-iso}(1), we deduce that $M_s\cong N_s.$ The proof of the theorem is completed.

\medskip

\noindent{\sc Remark.} If $\La$ is a finite dimensional $k$-algebra, then $D^b(\ModbLa)$ is the bounded derived category of all $\La$-modules.

\medskip

In view of Lemma \ref{F-pi}(1), we see that the same result holds true for $D^b({\rm mod}^{\hspace{.4pt}b\hspace{-2.5pt}}\La)$. To state it explicitly, we denote by ${\rm ind}^{-}(\tilde{Q}^{\hspace{.5pt}\rm op})$ a complete set of isomorphism class
representatives of the indecomposable objects of ${\rm rep}^{-}(\tilde{Q}^{\hspace{.5pt}\rm op})$, containing the indecomposable injective representations $I_{x^{\rm o}}$ with $x\in \tilde{Q}_0$ and the finite dimensional indecomposable projective representations $P_{y^{\rm o}}$ with $y\in \tilde{Q}_0$.

\medskip

\begin{Theo} \label{Main-2}

Let $\La=kQ/(kQ^+)^2$ with $Q$ a connected locally finite quiver, and let $\pi: \tilde{Q}\to Q$ be a minimal gradable covering
of $Q$.

\begin{enumerate}[$(1)$]

\item The derived Koszul push-down $\mf{F}_\pi: D^b({\rm rep}^{-}(\tilde{Q}^{\hspace{.5pt}\rm op}))\to D^b({\rm mod}^{\hspace{.4pt}b\hspace{-3pt}}\La) \vspace{1pt}$ is a Galois $\mf{G}$-covering, where $\mf{G}$ is the shifted translation group of $C^b({\rm Rep}^{-}(\tilde{Q}^{\hspace{.5pt}\rm op}))$.

\vspace{1pt}

\item The non-isomorphic indecomposable objects of $D^b({\rm mod}^{\hspace{.4pt}b\hspace{-3pt}}\La)$ are $\mf{F}_\pi(M)[s]$ with
$M\in {\rm ind}^-(\tilde{Q}^{\hspace{.5pt}\rm op})$ and $s\in \mathbb{Z}_{\hspace{.4pt}r},$ where $r$ is the grading period of $Q$.

\vspace{1pt}

\item Every object $X^\cdt$ of $D^b({\rm mod}^{\hspace{.4pt}b\hspace{-3pt}}\La)$ is uniquely $($up to isomorphism$)$ decomposed as \vspace{-1pt}
$$X^\cdt=\textstyle{\oplus}_{\hspace{.4pt} s\in \mathbb{Z}_{\hspace{.4pt}r}} \hspace{.4pt} \mf{F}_\pi(M_s)[s], \vspace{-2pt}$$ where $M_s\in {\rm rep}^-(\tilde{Q}^{\hspace{.5pt}\rm op}),$ all but finitely many of them are zero.

\end{enumerate}

\end{Theo}

\smallskip

\section{\sc Auslander-Reiten components}

\medskip

Throughout this section, let $\La\vspace{0pt}$ stand for a connected elementary locally bounded $k$-category with radical squared zero. It is well known that $D^b({\rm mod}^{b\hspace{-2.7pt}}\La)$ is a Hom-finite Krull-Schmidt $k$-category; see \cite{BaS} and \cite[(1.9)]{BaL}. Our objective of this section is to study the Auslander-Reiten theory in $D^b({\rm mod}^{b\hspace{-2.7pt}}\La)$, and in particular, to describe the shapes of its Auslander-Reiten components. 

\medskip

We shall keep all the notation introduced in the previous section. As usual, we may assume that $\La=kQ/(kQ^+)^2$, where $Q$ is a connected locally finite quiver. Fix a minimal gradable covering $\pi: \tilde{Q}\to Q$ of $Q$, which is a Galois covering with respect to the group $G$ generated by the translation $\rho$ of $\tilde{Q}$.  By Theorem \ref{Main-2}, the derived Koszul push-down $\mf{F}_\pi: D^b({\rm rep}^-(\tilde{Q}^{\hspace{.5pt}\rm op})) \to D^b({\rm mod}^{\hspace{.5pt}b\hspace{-2pt}}\La)$ is a Galois covering with respect to the shifted translation group $\mf{G}$ of $C^b({\rm rep}^-(\tilde{Q}))$ generated by the shifted translation $\vartheta=\rho\circ [-r_{\hspace{-1.5pt}_Q}]$, where $r_{\hspace{-1.5pt}_Q}$ is the grading period of $Q$. A complex in $D^b({\rm mod}^{b\hspace{-3pt}}\La)$ is called {\it simple} if it is isomorphic to a shift of some simple $\La$-module. The following result lists all the simple complexes up to isomorphism in $D^b({\rm mod}^{b\hspace{-3pt}}\La)$ and describes the irreducible morphisms between them, where $S[a]$ with $a\in Q_0$ denotes by the simple $\La$-module supported by $a.$

\medskip

\begin{Prop}\label{S-Irr}

Let $\La=kQ/(kQ^+)^2\vspace{1pt}$ with $Q$ a connected locally finite quiver, and let $\mf{F}_\pi: D^b({\rm rep}^-(\tilde{Q}^{\hspace{.5pt}\rm op}))\to D^b({\rm mod}^b\hspace{-2.7pt}\La)$ be the Koszul derived push-down associated with a minimal gradable covering $\pi: \tilde{Q}\to Q$ of $Q$.

\vspace{-1pt}

\begin{enumerate}[$(1)$]

\item If $a\in Q_0$ and $n \in \Z$, then $S[a][n]\cong \mf{F}_\pi(I_{x^{\rm o}})[s]$ for some $x\in \tilde{Q}_0$ and $s\in \Z_{\, r_{\hspace{-1pt}_Q}}\hspace{-1pt}.$

\vspace{1pt}

\item If $\alpha_i: a\to a_i$, $i=1, \ldots, s,\vspace{.5pt}$ are the arrows in $Q$ starting at a vertex $a$, then $D^b({\rm mod}^{b\hspace{-2.7pt}}\La)$ has an irreducible morphism $f^\cdt: S[a]\to \oplus_{i=1}^s\, S[a_i][1].$

\end{enumerate}

\end{Prop}

\noindent{\it Proof.} (1) Given a vertex $x=(a, t)\in \tilde{Q}$ with $a\in Q_0$ and $t\in \Z$, by Lemma \ref{KPD-image}(1), $\mf{F}_\pi(I_{x^{\rm o}})\cong S[a][t]$. By our choice of $\tilde{Q}$, there exists some vertex $x_0=(a_0, 0)\in \tilde{Q}$ with $a_0\in Q_0$. Let $a\in Q_0$ and $n \in \Z$. Being connected, $Q$ contains a walk from $a_0$ to $a$, say of degree $d$. By Lemma \ref{Quiver-cov}(3), $y=(a, d)\in \tilde{Q}$. In view of Proposition \ref{KPD-image}(2), we see that $S[a][n]\cong \mf{F}_\pi(I_{y^{\rm o}})[n-d]$. Write $n-d=s+ m r_{\hspace{-1pt}_Q}\hspace{-1pt}$ with $s\in \Z_{\, r_{\hspace{-1pt}_Q}}\vspace{1pt}$ and $m\in \Z$. By Proposition \ref{KPD-image}(3), $S[a][n]\cong \mf{F}_\pi(I_{y^{\rm o}})[mr_{\hspace{-1pt}_Q}][s]\cong (\rho^m \cdot I_{y^{\rm o}})[s]\cong \mf{F}_\pi(I_{(\rho^m \cdot y)^{\rm o}})[s].$

\vspace{1pt}

(2) Let $\alpha_i: a\to a_i$, $i=1, \ldots, s$, be the arrows in $Q$ starting at $a$. Considering the minimal gradable covering $\pi: \tilde{Q}\to Q,\vspace{1pt}$ we see that $a=\pi(x)$ for some vertex $x$ and $\alpha_i=\pi(\beta_i)$ for some arrow $\beta_i: x\to x_i$ in $\tilde{Q}_1$, $i=1, \ldots, s$. Then $\beta_1, \ldots, \beta_s$ are the arrows in $\tilde{Q}$ starting at $x$, and consequently, $\beta_1^{\hspace{.4pt}\rm o}, \ldots, \beta_s^{\hspace{.4pt}\rm o}$ are the arrows in $\tilde{Q}^{\hspace{.5pt}\rm op}$ ending at $x$. Therefore, $I_{x^{\rm o}}/{\rm soc}I_{x^{\rm o}}\cong \oplus_{i=1}^s\, I_{x_i^{\rm o}}$. Observe that the canonical projection $p_{x^{\rm o}}: I_{x^{\rm o}}\to \oplus_{i=1}^s\, I_{x_i^{\rm o}}$ is minimal left almost split; see \cite[(2.3)]{BLP}, and in particular, it is irreducible in ${\rm rep}^-(\tilde{Q}^{\hspace{.5pt}\rm op})$. Since ${\rm rep}^-(\tilde{Q}^{\hspace{.5pt}\rm op})$ is hereditary, it is easy to see that $p_{x^{\rm o}}$ remains irreducible in $D^b({\rm rep}^-(\tilde{Q}^{\hspace{.5pt}\rm op}))$. Suppose that $x\in \tilde{Q}^n$ for some $n\in \Z$. Then $\mf{F}_\pi(I_{x^{\rm o}})\cong S[a][n]$; and since $x_i\in \tilde{Q}^{n+1}$, we have $\mf{F}_\pi(I_{x_i^{\rm o}})\cong S[a_i][n+1],$ for $i=1, \ldots, s$. Since $\mf{F}_\pi$ is a Galois $\mf{G}$-covering, $\mathscr{F}_\pi(p_x): S[a][n]\to \oplus_{i=1}^s\, S[a_i][n+1]$ is an irreducible morphism in $D^b({\rm mod}^{b\hspace{-2.7pt}}\La);$ see \cite[(3.3)]{BaL}, and consequently, we obtain an irreducible morphism $f^\cdt: S[a]\to \oplus_{i=1}^s\, S[a_i][1]$. The proof of the proposition is completed.

\medskip

Next, we shall describe the almost split triangles in $D^b({\rm mod}^b\hspace{-2pt}\La)$. For this purpose, we denote by ${\rm ind}^b(\tilde{Q}^{\rm op})$ the set of finite dimensional representations in ${\rm ind}^-(\tilde{Q}^{\rm op})$.

\medskip

\begin{Theo}\label{ART}

Let $\La=kQ/(kQ^+)^2\vspace{.6pt}$ with $Q$ a connected locally finite quiver, and let $\mf{F}_\pi: D^b({\rm rep}^-(\tilde{Q}^{\hspace{.5pt}\rm op}))\to D^b({\rm mod}^b\hspace{-2.7pt}\La)$ be the derived Koszul push-down associated with a minimal gradable covering $\pi: \tilde{Q}\to Q$ of $Q$.

\begin{enumerate}[$(1)$]

\item If $M\in {\rm ind}^b(\tilde{Q}^{\hspace{.5pt}\rm op})$ is not projective, then $D^b({\rm mod}^b\hspace{-3pt}\La)$ has an almost split triangle $\xymatrixcolsep{16pt}\xymatrix{\mf{F}_\pi(L) \ar[r] & \mf{F}_\pi(N) \ar[r] & \mf{F}_\pi(M) \ar[r] & \mf{F}_\pi(L)[1],}\vspace{-2pt}$ which is induced from an almost split sequence
    $\xymatrixcolsep{18pt}\xymatrix{0\ar[r] & L \ar[r] & N \ar[r] & M \ar[r] & 0}$ in ${\rm rep}^-(\tilde{Q}).$

    \vspace{2pt}

\item If $P_{x^{\rm o}}\in {\rm ind}^b(\tilde{Q}^{\rm op})\vspace{1pt}$ with $x\in \tilde{Q}_0,$ then $D^b({\rm mod}^b\hspace{-2.7pt}\La)$ has an almost split triangle $\xymatrixcolsep{16pt}\xymatrix{\mf{F}_\pi(I_{x^{\rm o}})\ar[r] & \mf{F}_\pi(I_{x^{\rm o}}/{\rm soc}I_{x^{\rm o}}) \oplus \mf{F}_\pi({\rm rad}P_{x^{\rm o}})[1]
\ar[r] & \mf{F}_\pi(P_{x^{\rm o}})[1] \ar[r] & \mf{F}_\pi(I_{x^{\rm o}})[1].}$

\vspace{0pt}

\item Every almost split triangle in $D^b({\rm mod}^b\hspace{-3pt}\La)$ is isomorphic to a
shift of an almost split triangle stated in $(1)$ or $(2)$.

\end{enumerate} \end{Theo}

\noindent{\it Proof.} (1) Let $M$ be a non-projective representation in ${\rm ind}^b(\tilde{Q}^{\hspace{.5pt}\rm op})$.
Then, ${\rm rep}^-(\tilde{Q}^{\hspace{.5pt}\rm op})$ has an almost split sequence $\xymatrixcolsep{18pt}\xymatrix{0\ar[r] & L \ar[r] & N \ar[r] & M \ar[r] & 0,}$ which induces an almost split triangle $\xymatrixcolsep{18pt}\xymatrix{L \ar[r] & N \ar[r] & M \ar[r] & L[1]}\vspace{-1pt}$ in $D^b({\rm rep}^-(\tilde{Q}^{\hspace{.5pt}\rm op}))$; see \cite[(2.8), (7.2)]{BLP}. Applying the Galois $G$-covering $\mf{F}_\pi: D^b({\rm rep}^-(\tilde{Q}^{\hspace{.5pt}\rm op})) \to D^b({\rm mod}^{\hspace{.5pt}b\hspace{-2pt}}\La)$, we obtain an almost split triangle $\xymatrixcolsep{18pt}\xymatrix{\mf{F}_\pi(L) \ar[r] & \mf{F}_\pi(N) \ar[r] & \mf{F}_\pi(M) \ar[r] & \mf{F}_\pi(L)[1]}\vspace{-2pt}$ in $D^b({\rm mod}^b\hspace{-2.7pt}\La)$; see \cite[(3.7)]{BaL}.

\vspace{1pt}

(2) Let $P_{x^{\rm o}}\in $ be finite dimensional for some $x\in \tilde{Q}_0$. By the dual of Theorem 7.5 in \cite{BLP}, there exists
an almost split triangle \vspace{-3pt}
$$\xymatrixcolsep{16pt}\xymatrix{I_{x^{\rm o}}\ar[r] & (I_{x^{\rm o}}/{\rm soc}I_{x^{\rm o}}) \oplus ({\rm rad}P_{x^{\rm o}})[1]
\ar[r] &P_{x^{\rm o}}[1] \ar[r] & I_{x^{\rm o}}[1]}\vspace{-5pt}$$ in
$D^b({\rm rep}^-(\tilde{Q}^{\hspace{.5pt}\rm op}))$. Applying the Galois $G$-covering $\mf{F}_\pi,\vspace{1pt}$ we obtain a desired almost split triangle in $D^b({\rm mod}^b\hspace{-2.7pt}\La)$; see \cite[(3.7)]{BaL} and \cite[(6.1)]{Liu2}.

\vspace{1pt}

(3) Let us consider an almost split triangle $\zeta^\cdt: \xymatrixcolsep{16pt}\xymatrix{X^\cdt \ar[r] & Y^\cdt \ar[r] & Z^\cdt \ar[r] & X^\cdt[1]}\vspace{-2pt}$ in $D^b({\rm mod}^b\hspace{-2.7pt}\La).$ By Theorem \ref{Main-2}(2), we may assume that $Z^\pdt=\mf{F}_\pi(M)$ for some object $M\in {\rm ind}^-(\tilde{Q})$. Combining the results in \cite[(3.7)]{BaL} and \cite[(6.1)]{Liu2}, we see that $\zeta^\cdt=\mf{F}_\pi(\eta^\cdt)$, where $\eta^\cdt$ is an almost split triangle in $D^b({\rm rep}^-(\tilde{Q}^{\hspace{.5pt}\rm op}))$ ending at $M$. By the dual of Theorem 7.5 in \cite{BLP}, in case
$M$ is not projective, $\eta^\cdt$ is an almost split triangle $\xymatrixcolsep{18pt}\xymatrix{L \ar[r] & N \ar[r] & M \ar[r] & L[1],}\vspace{-1.5pt}$ induced from an almost split sequence $\xymatrixcolsep{18pt}\xymatrix{0\ar[r] & L \ar[r] & N \ar[r] & M \ar[r] & 0}\hspace{-2pt}$ in $\rep^-(\tilde{Q}^{\hspace{.5pt}\rm op});$ and otherwise, $\eta$ is an almost split triangle
$\xymatrixcolsep{18pt}\xymatrix{\hspace{-1pt}I_{x^{\rm o}}\ar[r] & (I_{x^{\rm o}}/{\rm soc}I_{x^{\rm o}}) \oplus {\rm rad}P_{x^{\rm o}}[1] \ar[r] &P_{x^{\rm o}}[1] \ar[r] & I_{x^{\rm o}}[1]}\hspace{-2pt}$ with $P_{x^{\rm o}}\in {\rm ind}^b(\tilde{Q}^{\hspace{.5pt}\rm op}).$ The proof of the theorem is completed.

\medskip

By the dual of Corollary 4.15 in \cite{BLP}, ${\rm rep}^-(\tilde{Q}^{\,\rm op})$ has symmetric irr-spaces. We choose the vertex set of the Auslander-Reiten quiver $\Ga_{{\rm rep}^-(\tilde{Q}^{\,\rm op})}$ of ${\rm rep}^-(\tilde{Q}^{\,\rm op})$ in such a way that it contains all the $I_{x^{\rm o}}$ with $x\in \tilde{Q}_0$ and the finite dimensional $P_{y^{\rm o}}$ with $y\in \tilde{Q}_0.$ By the dual of Theorem 4.6 in \cite{BLP}, all the $I_{x^{\rm o}}$ with $x\in \tilde{Q}_0$ lie in a connected component of $\Ga_{{\rm rep}^-(\tilde{Q}^{\,\rm op})},\vspace{1pt}$ called {\it preinjective component}. Moreover, a connected component of $\Ga_{{\rm rep}^-(\tilde{Q}^{\,\rm op})}$ is called {\it preprojective} if it contains some of the finite dimensional $P_{y^{\rm o}}$ with $y\in \tilde{Q}_0$, and {\it regular} if it contains none of the $I_{x^{\rm o}}$ with $x\in Q_0$ and the finite dimensional $P_{y^{\rm o}}$ with $y\in Q_0$. A representation in $\Ga_{{\rm rep}^-(\tilde{Q}^{\,\rm op})}\vspace{.5pt}$ is called {\it preinjective, preprojective} or {\it regular} if it lies in a preinjective, preprojective or regular component.

\medskip

Next, by the dual of Lemma 7.7(3) in \cite{BLP}, $D^b({\rm rep}^-(\tilde{Q}^{\,\rm op}))$ has symmetric irr-spaces. We shall choose the vertices of the Auslander-Reiten quiver $\Ga_{D^b({\rm rep}^-(Q))}\vspace{1pt}$ of $D^b({\rm rep}^-(Q))$ to be the shifts of the vertices of $\Ga_{{\rm rep}^-(Q)}\vspace{.6pt}$. In view of the dual of Theorem 7.5 in \cite{BLP}, we see that a regular component of $\Ga_{{\rm rep}^-(Q)}$ is a connected component of $\Ga_{D^b({\rm rep}^-(Q))}\vspace{.6pt}$. Moreover, by the dual statements of Lemmas 7.7 and 7.8 in \cite{BLP}, the preinjective component of $\Ga_{{\rm rep}^-(Q)}$ is glued together with the shifts by $1$ of all possible preprojective components to form a connected subquiver of $\Ga_{D^b({\rm rep}^-(Q))},\vspace{.6pt}$ which is contained in the so-called  {\it connecting component}.

\medskip

Furthermore, we deduce from Theorem \ref{ART}
that $D^b({\rm mod}^b\hspace{-2pt}\it\Lambda)$ has symmetric irr-spaces. Observe that
the free $\mf{G}$-action on $D^b({\rm rep}^-(\tilde{Q}^{\hspace{.5pt}\rm op}))$ induces a $\mf{G}$-action on $\Ga_{D^b({\rm rep}^-(\tilde{Q}^{\hspace{.5pt}\rm op}))}\vspace{1pt}$.
Fix a complete set $\Oa$ of $\mf{G}$-orbit representatives of the vertices of $\Ga_{D^b({\rm rep}^-(\tilde{Q}^{\hspace{.5pt}\rm op}))},\vspace{1pt}$ which contains all the $I_{x^{\rm o}}$ and the finite dimensional $P_{y^{\rm o}}$. Since $\mf{F}_\pi: D^b({\rm rep}^-(\tilde{Q}^{\,\rm op}))\to D^b({\rm mod}^b\hspace{-2pt}\it\Lambda)$ is a Galois $\mf{G}$-covering, we may choose the vertices of $\Ga_{D^b({\rm mod}^b\hspace{-2pt}\it\Lambda)}\vspace{1pt}$ to be the objects $\mf{F}_\pi(X^\cdt)$ with $X^\cdt\in \Oa$. By Theorem 4.7 in \cite{BaL}, we have an induced Galois $\mf{G}$-covering $\mf{F}_\pi: \Ga_{D^b({\rm rep}^-(\tilde{Q}^{\hspace{.5pt}\rm op}))} \to \Ga_{D^b({\rm mod}^b\hspace{-2pt}\it\Lambda)}.\vspace{1pt}$ If $\Ga$ is a connected component of $\Ga_{D^b({\rm rep}^-(\tilde{Q}^{\hspace{.5pt}\rm op}))},\vspace{1pt}$ then the full subquiver $\mf{F}_\pi(\Ga)$ of $\Ga_{D^b({\rm mod}^b\hspace{-1.7pt}\it\Lambda)}\vspace{1pt}$ generated by the $\mf{F}_\pi(X^\cdt)$ with $X^\cdt\in \Ga$ is a connected component.

\medskip

\begin{Theo}\label{ARQ}

Let $\La=kQ/(kQ^+)^2\vspace{.5pt}$ with $Q$ a connected locally finite quiver, and let $\mf{F}_\pi: D^b({\rm rep}^-(\tilde{Q}^{\hspace{.5pt}\rm op}))\to D^b({\rm mod}^b\hspace{-2.7pt}\La)$ be the derived Koszul push-down associated with a minimal gradable covering $\pi: \tilde{Q}\to Q$ of $Q$.

\vspace{-1pt}

\begin{enumerate}[$(1)$]

\item Let $\mathcal{C}$ be the connecting component of $\Ga_{D^b({\rm rep}^-(\tilde{Q}^{\hspace{.5pt}\rm op}))}\vspace{1pt}$ and $\rho$ the translation of $\tilde{Q}.$ If $n\in \Z$, then $\rho^n\cdot \mathcal{C}=\mathcal{C}$ and $\mf{F}_\pi(\mathcal{C})[n]=\mf{F}_\pi(\mathcal{C})$ whenever $n\equiv \, 0 \, ({\rm mod}\, r_{\hspace{-1.5pt}_Q}).$

\vspace{1.5pt}

\item If $\Ga$ is a connected component of $\Ga_{D^b({\rm rep}^-(\tilde{Q}^{\hspace{.5pt}\rm op}))},\vspace{1pt}$ then $\mf{F}_\pi:  \Ga \to \mf{F}_\pi(\Ga)$ is a translation quiver isomorphism.

\vspace{1.5pt}

\item The connected components of $\Ga_{D^b({\rm mod}^b\hspace{-2pt}\it\Lambda)}\vspace{.4pt}$ are $\mf{F}_\pi(\Ga)[s]$, where $s\in \Z_{\,r_{\hspace{-1.5pt}_Q}}\vspace{.6pt}$ and $\Ga$ is either the connecting component of  $\Ga_{D^b({\rm rep}^-(\tilde{Q}^{\hspace{.5pt}\rm op}))}\vspace{1pt}$ or a regular component of $\Ga_{{\rm rep}^-(\tilde{Q}^{\hspace{.5pt}\rm op})},$ which are pairwise distinct in case $Q$ is not of Dynkin type.

\end{enumerate}\end{Theo}

\vspace{1.5pt}

\noindent{\it Proof.} (1) For simplicity, we write $r=r_{\hspace{-1.3pt}_Q}.$ Let $n\in \Z$. Since $\rho^n\cdot I_{x^{\rm o}}=I_{(\rho^n\cdot x)^{\rm o}}\in \mathcal{C}$ for all $x\in \tilde{Q}_0$, we see that $\rho^n\cdot \mathcal{C}=\mathcal{C}$.
If $n=t r$ for some $t\in \Z$, by Proposition \ref{KPD-image}(2), we obtain $\mf{F}_\pi(I_{(\rho^t\cdot x)^{\rm o}})=\mf{F}_\pi(\rho^t\cdot I_{x^{\rm o}})= \mf{F}_\pi(I_{x^{\rm o}})[t r_{\hspace{-1pt}_Q}]\in \mf{F}_\pi(\mathcal{C})\cap \mf{F}_\pi(\mathcal{C})[n].$ Therefore, $\mf{F}_\pi(\mathcal{C})[n]=\mf{F}_\pi(\mathcal{C})$.

\vspace{1pt}

(2) Let $\Ga$ be a connected component of $\Ga_{D^b({\rm rep}^-(\tilde{Q}^{\hspace{.5pt}\rm op}))}$. Then, $\Ga=\mathcal{C}[m]$ or $\Ga=\mathcal{R}[m]$, where $m\in \Z$ and $\mathcal{R}$ is a regular component of $\Ga_{{\rm rep}^-(\tilde{Q}^{\hspace{.5pt}\rm op})};\vspace{1pt}$ see \cite[(7.9),(7.10)]{BLP}. The Galois $\mf{G}$-covering $\mf{F}_\pi: \Ga_{D^b({\rm rep}^-(\tilde{Q}^{\hspace{.5pt}\rm op}))} \to \Ga_{D^b({\rm mod}^b\hspace{-1.5pt}\it\Lambda)}\vspace{1pt}$ restricts to a Galois $\mf{G}_{\hspace{-1pt}\it\Gamma}$-covering $\mf{F}_\pi: \Ga\to \mf{F}_\pi(\Ga)$, where $\mf{G}_{\hspace{-1pt}\it\Gamma}=\{g\in \mf{G} \mid g(\Ga)=\Ga\};$ see \cite[(4.7)]{BaL}. Let $g$ be a non-identity element of $\mf{G}$. Then $g=\vartheta^n=\rho^n\circ [-nr]$ with $nr\ne 0$. We claim that $g\cdot \Ga\ne \Ga.$ Indeed, since $\rho\cdot \mathcal{C}=\mathcal{C}$ by Statement (1), $\rho\cdot \mathcal{R}$ is a regular component of $\Ga_{{\rm rep}^-(\tilde{Q}^{\hspace{.5pt}\rm op})}.\vspace{1pt}$ If $\Ga=\mathcal{R}[m]$, then $g\cdot \Ga= (\rho^n\cdot \mathcal{R})[m-nr]\ne \Ga$ since $nr\ne 0$. Otherwise, $\Ga=\mathcal{C}[m]$ and $g\cdot \Ga=\mathcal{C}[m-nr].$ Since $r>0$, we see that $\tilde{Q}^{\rm op}$ is not of Dynkin type. Thus, the preinjective component of $\Ga_{{\rm rep}^-(\tilde{Q}^{\rm op})}$ does not intersect any preprojective component; see \cite[(2.14)]{BLP} and \cite{DRi}. In view of the dual of Theorem 7.5 in \cite{BLP}, $\mathcal{C}$ contains only the preinjective representations and the shifts by 1 of the preprojective representations. As a consequence, $\mathcal{C}[m-nr]\ne \mathcal{C}[m],$ and hence, $g\cdot \Ga \ne \Ga.$ This establishes our claim. That is, $\mf{G}_{\hspace{-1pt}\it\Gamma}$ is trivial, and consequently,
$\mf{F}_\pi: \Ga\to \mf{F}_\pi(\Ga)$ is an isomorphism; see the remark following \cite[(4.6)]{BaL}.

\vspace{1pt}

(3) Let $X^\cdt$ be any complex in $\Ga_{D^b({\rm mod}^b\hspace{-1.5pt}\it\Lambda)}\vspace{1pt}$. By Theorem \ref{Main-2}(2), $X^\cdt=\mf{F}_\pi(M)[s]$, for some $M\in \Ga_{{\rm rep}^-(\tilde{Q}^{\rm op})}\vspace{1pt}$ and $s\in \Z_{\,r}$. If $M$ lies in a regular component $\mathcal{R}$, then $X^\cdt=\mf{F}_\pi(M)[s]\in \mf{F}_\pi(\mathcal{R})[s].$
If $M$ is preinjective, then $X^\cdt\in \mf{F}_\pi(\mathcal{C})[s]$. Assume now that $M$ is preprojective. Then $X^\cdt\in \mf{F}_\pi(\mathcal{C})[s-1]$. If $r=0$ or $s\ge 1$, we obtain $s-1\in \Z_r$. Otherwise, $r>0$ and $s=0$. By Statement (1), we obtain $\mf{F}_\pi(\mathcal{C})[s-1]=\mf{F}_\pi(\mathcal{C})[-1]=\mf{F}_\pi(\mathcal{C})[r-1]$ with $r-1\in \Z_r.$

Suppose, moreover, that $X^\cdt\in \mf{F}_\pi(\Ga)[s] \cap \mf{F}_\pi(\Ta)[t]$, where $s, t\in \Z_{\,r}\vspace{1.3pt}$ and $\Ga, \Ta$ are either $\mathcal{C}$ or some regular components of $\Ga_{{\rm rep}^-(\tilde{Q}^{\hspace{.5pt}\rm op})}$. If $\Ga$ or $\Ta$ is a regular component of $\Ga_{{\rm rep}^-(\tilde{Q}^{\hspace{.5pt}\rm op})}$, then we deduce from the uniqueness stated in Theorem \ref{Main-2}(3) that $\Ga=\Ta$ and $s=t$. If $\Ga=\Ta=\mathcal{C}$ and $ s\ne t$, then $\mathcal{C}=\mathcal{C}[s-t]$ with $s-t\ne 0$. As seen above, this occurs only if $\tilde{Q}^{\,\rm op}$ is of Dynkin type, that is, $Q$ is of Dynkin type. The proof of the theorem is completed.

\medskip

We shall see that the simple complexes play an important role in the description of the Auslander-Reiten components of $D^b({\rm mod}^b\hspace{-2.3pt}\it\Lambda).$ For convenience, we shall call $\mathscr{C}=\mf{F}_\pi(\mathcal{C})$ the {\it connecting component} of $\Ga_{D^b({\rm mod}^b\hspace{-1.8pt}\it\Lambda)}\vspace{-1pt}$ determined by the minimal gradable covering $\pi: \tilde{Q}\to Q$. Recall that a complex in $D^b({\rm mod}^b\hspace{-2.3pt}\it\Lambda)$ is called {\it perfect} if it is isomorphic to a bounded complex over ${\rm proj}\hspace{.5pt}\La$. By Proposition \ref{F-bch}(3), the indecomposable perfect complexes in $D^b({\rm mod}^b\hspace{-2.3pt}\it\Lambda)$ are precisely $\mf{F}_\pi(M)[n]$, with $M\in {\rm ind}^b(\tilde{Q})$ and $n\in \Z$. Moreover, an infinite path in a quiver is called {\it right infinite} if it has a starting point; {\it left infinite} if it has an ending point; and {\it double infinite} if it has neither a starting point nor an ending point.

\medskip

\begin{Theo}\label{ARC-1}

Let $\La=kQ/(kQ^+)^2,$ where $Q$ is a connected locally finite non-Dynkin quiver, and let $\mathscr{C}$ be the
connecting component of $\Ga_{D^b({\rm mod}^b\hspace{-1.5pt}\it\Lambda)}\vspace{-1.5pt}$ determined by a minimal gradable covering
$\pi: \tilde{Q}\to Q$ of $Q$.

\vspace{-2pt}

\begin{enumerate}[$(1)$]

\item The component $\mathscr{C}$ contains some $S[a]$ with $a\in Q\hspace{.5pt};$ and in this case, $S[b][n]\in \mathscr{C}$ with $b\in Q$ and $n\in \Z$ if and only if $n\equiv d \hspace{.5pt}({\rm mod}\, r_{\hspace{-1.5pt}_Q})$, where $d$ is the degree of some walk in $Q$ from $a$ to $b$.

\vspace{1pt}

\item The simple complexes in $\mathscr{C}$ form a section of shape $\tilde{Q}$, and consequently, $\mathscr{C}$ embeds in $\Z\hspace{.4pt} \tilde{Q}.$

    \vspace{1pt}

\item The connected components of $\Ga_{D^b({\rm mod}^b\hspace{-1.5pt}\it\Lambda)}$ containing simple complexes are the  $\mathscr{C}[s]$ with $s\in \Z_{\,r_{_{\hspace{-1pt}Q}}}.$

\vspace{1pt}

\item If $Q$ has no right infinite path, then $\mathscr{C}$ is left stable and contains only perfect complexes$\,;$ and if $Q$ has no infinite path, then $\mathscr{C}\cong \Z \,\tilde{Q}$.

\vspace{1pt}

\item If $Q$ has right infinite paths, then $\mathscr{C}$ contains a left-most $\Da$ generated by its non-perfect complexes$\,;$ and if $Q$ in addition has no left infinite path, then $\mathscr{C}\cong \N \Da$.

\end{enumerate}

\end{Theo}

\noindent{\it Proof.} (1) By our choice of $\tilde{Q}$, the connecting component $\mathcal{C}$ of $\Ga_{D^b({\rm rep}^-(\tilde{Q}^{\hspace{.5pt}\rm op}))}\vspace{.5pt}$ contains some $I_{x^{\rm o}}$, where $x=(a, 0)$ with $a\in Q_0$. In view of Proposition \ref{KPD-image}(2), we see that $S[a]=\mf{F}_\pi(I_{x^{\rm o}})\in \mathscr{C}.$ Let $b\in Q_0$ and $n\in \Z$. Suppose first that $n=d + t r_{\hspace{-1.5pt}_Q}$, where $t\in \Z$ and $d$ is the degree of some walk in $Q$ from $a$ to $b$. By Lemma \ref{Quiver-cov}(3), $y=(b, n)\in \tilde{Q}^n$, and by Proposition \ref{KPD-image}(2), $S[b][n]=\mf{F}_\pi(I_{y^{\rm o}})\in \mathscr{C}.\vspace{1pt}$ Conversely, suppose that $S[b][n]\in \mathscr{C}$. Being connected, $Q$ contains a walk from $a$ to $b$, say of degree $d$. By Lemma \ref{Quiver-cov}(3), $z=(b, d)\in \tilde{Q}$, and by Proposition \ref{KPD-image}(2), $S[b][n]=\mf{F}_\pi(I_{z^{\rm o}})[n-d]$. Write $n-d=s+ t r_{\hspace{-1pt}_Q}$ with $t\in \Z$ and $s\in \Z_{\, r_{\hspace{-1pt}_Q}}$. In view of Proposition \ref{KPD-image}(3), we see that $S[b][n] =\mf{F}_\pi(I_{z^{\rm o}})[t r_{\hspace{-1pt}_Q}+s]=\mf{F}_\pi(I_{(\rho^t\cdot z)^{\rm o}})[s]\in \mathscr{C}\cap \mathscr{C}[s]$. As a consequence, $\mathscr{C}=\mathscr{C}[s]$. Since $Q$ is not of Dynkin type, $s=0$ by Theorem \ref{ARQ}(3). That is, $n\equiv d \, (\hspace{-.6pt}{\rm mod}\, r_{\hspace{-1pt}_Q})$.

(2) By the dual of Proposition 7.2 in \cite{BLP}, $\mathcal{C}$ has a section $\Sa$ of shape $\tilde{Q}$ which is generated by the $I_{x^{\rm o}}$ with $x\in \tilde{Q}_0$. Since $\mf{F}_\pi: \mathcal{C}\to \mathscr{C}$ is an isomorphism; see (\ref{ARQ}), $\mf{F}_\pi(\Sa)$ is a section of $\mathscr{C}$ generated by the simple complexes $\mf{F}_\pi(I_{x^{\rm o}})$ with $x\in \tilde{Q}_0$. As seen in the proof of Statement (1), all the simple complexes in $\mathscr{C}$ belong to $\mf{F}_\pi(\Sa)$.

(3) Let $S[b]][n]$ be a simple complex in $\Ga_{D^b({\rm mod}^b\hspace{-2.5pt}\it\Lambda)},$ where $b\in Q_0$ and $n\in \Z$. By Proposition \ref{S-Irr}(1), $S[b][n]=\mf{F}_\pi(I_{x^{\rm o}})[s]$ for some $x\in \tilde{Q}_0$ and $s\in \Z_{\, r_{\hspace{-1pt}_Q}}$, which belongs to $\mathscr{C}[s]$.

(4) Suppose that $Q$ contains no right infinite path. Then $\tilde{Q}$ has no right infinite path, that is, $\tilde{Q}^{\hspace{.5pt}\rm op}$ has no left infinite path. In this case, all the $I_{x^{\rm o}}$ with $x\in \tilde{Q}_0$ are finite dimensional, and consequently, ${\rm rep}^-(\tilde{Q}^{\hspace{.5pt}\rm op})={\rm rep}^b(\tilde{Q}^{\hspace{.5pt}\rm op})$. In particular, all the complexes in $\mathscr{C}$ are perfect. Moreover, if $Q$ has left infinite path, then $\tilde{Q}^{\hspace{.5pt}\rm op}$ has right infinite paths. By the dual of Proposition 7.9(3) in \cite{BLP}, $\mathcal{C}$ is left stable. If $Q$ has no left infinite path, then $\tilde{Q}^{\hspace{.5pt}\rm op}$ has no right infinite path, and hence, it has no infinite path. By the dual of Proposition 7.9(1) in \cite{BLP}, $\mathcal{C}$ is stable of shape $\Z \tilde{Q}$. Now, we deduce from Theorem \ref{ARQ}(2) that $\mathscr{C}$ is always left stable and $\mathscr{C}\cong \Z \tilde{Q}$ if $Q$ has no infinite path.

\vspace{1pt}

(5) Suppose that $Q$ contains right infinite paths. Then $\tilde{Q}^{\hspace{.5pt}\rm op}$ contains left infinite paths. That is, some of the $I_{x^{\rm o}}$ with $x\in \tilde{Q}$ are infinite dimensional. By the dual of Lemma 4.5(1) in \cite{BLP}, the infinite dimensional representations in the preinjective component of $\Ga_{{\rm rep}^-(\tilde{Q}^{\hspace{.5pt}\rm op})}$ form a left-most section $\Ta$, which is also a left-most section of $\mathcal{C}$. Moreover, by the dual of Theorem 4.7 in \cite{BLP}, the preprojective components of $\Ga_{{\rm rep}^-(\tilde{Q}^{\hspace{.5pt}\rm op})}$ contain only finite dimensional representations. As a consequence, $\Da=\mf{F}_\pi(\Ta)$ is a left-most section of $\mathscr{C}$ generated by its non-perfect complexes. Suppose further that $Q$ contains no left infinite path. Then $\tilde{Q}$ contains no right infinite path. By the dual of Proposition 7.9(2) in \cite{BLP}, $\mathcal{C}$ is of shape $\N \Ta$, and consequently, $\mathscr{C}$ is of shape $\N \Da$. The proof the theorem is completed.

\medskip

Next, we shall describe the Auslander-Reiten components of $D^b({\rm mod}^b\hspace{-2.3pt}\it\Lambda)$ without simple complexes. For this purpose, we denote by $\A_\infty^+$ a quiver which is a right infinite path and by $\A_\infty^-$ a quiver which is a left infinite path. Recall that a translation quiver  is called a {\it stable tube} if it is of shape $\Z \A_\infty / \hspace{-3pt} <\hspace{-1.5pt}\tau^n \hspace{-2pt}> $ for some positive integer $n$; and a {\it wing} if it is of the following shape$\,:$
\vspace{0pt}
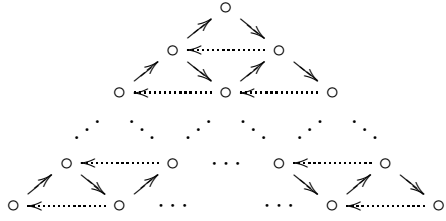
$$\xymatrixrowsep{16pt}\xymatrixcolsep{20pt}
\xymatrix@!=0.1pt{&&&& \circ \ar[dr]&&&&\\
&&&\circ \ar[ur]\ar[dr] \ar@{<.}[rr]&& \circ \ar[dr]&&&\\
&& \circ\ar[ur] \ar@{<.}[rr] && \circ\ar[ur] \ar@{<.}[rr] && \circ
&&\\
}$$ \vspace{-16pt}
$$\iddots \hspace{8pt} \ddots \hspace{8pt} \iddots \hspace{8pt}  \ddots  \hspace{8pt} \iddots \hspace{8pt} \ddots$$
\vspace{-24pt}
$$\xymatrixrowsep{16pt}\xymatrixcolsep{20pt}\xymatrix@!=0.1pt{& \circ\ar[dr] \ar@{<.}[rr] && \circ   &\cdots& \circ\ar[dr]  \ar@{<.}[rr]  && \circ \ar[dr] & \\
\circ \ar[ur]\ar@{<.}[rr] && \circ \ar[ur]  &\cdots&&\cdots&
\circ\ar[ur] \ar@{<.}[rr]
 && \circ}\vspace{6pt}
$$
where the dotted arrows indicate the translation; see \cite[(3.3)]{Rin2}.

\medskip

\begin{Theo}\label{ARC-2}

Let $\La=kQ/(kQ^+)^2$ with $Q$ a connected locally finite quiver, and let $\mathscr{R}$ be a connected component of $\Ga_{D^b({\rm mod}^b\hspace{-2pt}\it\Lambda)}\vspace{1pt}$ without simple complexes.

\vspace{-2pt}

\begin{enumerate}[$(1)$]

\item If $\mathscr{R}$ contains only perfect complexes, then it is a stable tube or of shape $\Z \A_\infty$ or $\N^-\A_\infty^-$, where the first case occurs only if $Q$ is gradable of Euclidean type.

\item If $\mathscr{R}$ has non-perfect complexes, then it is a wing or of shape $\N \A_\infty^+$, where the left-most section is generated by the non-perfect complexes.

\end{enumerate}

\end{Theo}

\noindent{\it Proof.} Let $\pi: \tilde{Q}\to Q$ be a minimal gradable covering of $Q$. By Theorems \ref{ARQ} and \ref{ARC-1}, $\mathscr{R}=\mf{F}_\pi(\mathcal{R})[s]$, where $s\in \Z_{\,r_{_{\hspace{-1pt}Q}}}$ and $\mathcal{R}$ is a regular component of $\Ga_{{\rm rep}^-(\tilde{Q}^{\hspace{.5pt}\rm op})}.$ \vspace{1pt}

(1) Suppose that $\mathscr{R}$ contains only perfect complexes. Then $\mathcal{R}$ contains only finite dimensional representations. Assume that $\tilde{Q}^{\,\rm op}$ is finite. By Lemma \ref{Quiver-cov}, $Q$ is a finite gradable quiver. It is well known that $\mathcal{R}$ is a stable tube or of shape $\Z \A_\infty$, where the first case occurs only if $\tilde{Q}^{\,\rm op}$ is of Euclidean type; see \cite{DRi}, \cite{Rin} and \cite[(3.6)]{Rin2}. Observe that $\tilde{Q}^{\,\rm op}$ is of Euclidean type if and only if $Q$ is gradable of Euclidean type. If $\tilde{Q}^{\,\rm op}$ is infinite, then we deduce from
the dual statements of Theorems 4.14(1) and 4.14(3) in \cite{BLP} that $\mathcal{R}$ is of shape $\Z \A_\infty$ or $\N^-\A_\infty^-$. Statement (1) follows now from Theorem \ref{ARQ}(2).

(2) Suppose that $\mathscr{R}$ contains some non-perfect complexes. Then, $\mathcal{R}$ contains some infinite dimensional representations. Applying the dual statements of Lemma 4.5(1) and Theorems 4.14(2) and 4.14(4) in \cite{BLP}, we conclude that
$\mathcal{R}$ is a wing or of shape $\N \A_\infty^+$, where the left-most section is generated by the infinite dimensional representations. Statement (2) follows from Theorem \ref{ARQ}(2). The proof of the theorem is completed.

\medskip

If $Q$ has no infinite path, such as $\La$ is finite dimensional of finite global dimension, then we have a nicer description of the Auslander-Reiten components of $D^b({\rm mod}^b\hspace{-2pt}\it\Lambda)$.

\medskip

\begin{Theo}\label{ARC-3}

Let $\La=kQ/(kQ^+)^2,$ where $Q$ is a connected locally finite quiver with no infinite path. If $\tilde{Q}$ is a minimal gradable covering of $Q$, then every connected component of $\Ga_{D^b({\rm mod}^b\hspace{-2pt}\it\Lambda)}$ is a stable tube or of shape $\mathbb{Z}\tilde{Q}^{\hspace{.5pt}\rm op}$ or $\mathbb{Z}\mathbb{A}_\infty$, where the first case occurs only if $Q$ is gradable of Euclidean type.

\end{Theo}

\noindent{\it Proof.} Let $\tilde{Q}$ be a minimal gradable covering of $Q$. Then $\tilde{Q}^{\hspace{.5pt}\rm op}$ has no infinite path, and ${\rm rep}^-(\tilde{Q}^{\hspace{.5pt}\rm op})={\rm rep}^+(\tilde{Q}^{\hspace{.5pt}\rm op})={\rm rep}^b(\tilde{Q}^{\hspace{.5pt}\rm op}).$ In this case, the connecting component of $\Ga_{D^b({\rm rep}^-(\tilde{Q}^{\hspace{.5pt}\rm op}))}\vspace{.5pt}$ is of shape $\Z \tilde{Q}$; see \cite[(7.9)]{BLP}. Moreover, every regular component of $\Ga_{{\rm rep}^-(\tilde{Q}^{\hspace{.5pt}\rm op})}$ is stable tube or of shape $\Z\A_\infty$, where the first case occurs only if $Q$ is gradable of Euclidean type; see \cite[(4.16)]{BLP}, \cite{DRi}, \cite{Rin} and \cite[(3.6)]{Rin2}. The result now follows from  Theorem \ref{ARQ}(2). The proof of the theorem is completed.

\medskip

To conclude the paper, we shall determine when $D^b({\rm mod}^b\hspace{-2.5pt}\La)$ has only finite many Auslander-Reiten components. For this purpose, we say that a quiver is of {\it type $\tilde{\A}_n$} with $n\ge 1\vspace{1pt}$ if its underlying graph is a cycle of $n$ edges
; compare \cite[Page 7]{Rin2}.

\smallskip

\begin{Theo}\label{ARC-4}

Let $\La=kQ/(kQ^+)^2,$ where $Q$ is a connected locally finite quiver.

\begin{enumerate}[$(1)$]

\item If $Q$ is of Dynkin type, then $\Ga_{D^b({\rm mod}^b\hspace{-2pt}\it\Lambda)}$ is connected of shape $\Z Q\hspace{.4pt}.$

\vspace{.8pt}

\item If $Q$ is an oriented cycle of $n$ arrows, then $\Ga_{D^b({\rm mod}^b\hspace{-1.5pt}\it\Lambda)}$ consists of $n$ components of shape $\Z\A_\infty$ 
    and $n$ double infinite paths with only simple complexes.

\vspace{.9pt}

\item If $Q$ is of type $\tilde{\A}_n$ with $0<r_{\hspace{-1pt}_Q}<n$, then
$\Ga_{D^b({\rm mod}^b\hspace{-1.8pt}\it\Lambda)}\hspace{-1.5pt}$ 
consists of $2\hspace{.4pt}r_{\hspace{-1pt}_Q}$ components of shape $\mathbb{Z} \mathbb{A}_\infty$ and $r_{\hspace{-1pt}_Q}$ components of shape $\,\mathbb{Z}\tilde{Q}$.

\vspace{.8pt}

\item In all other cases, $\Ga_{D^b({\rm mod}^b\hspace{-1.8pt}\it\Lambda)}$ has infinite many connected components.

\end{enumerate}\end{Theo}

\noindent{\it Proof.} Let $\pi: \tilde{Q}\to Q$ be a minimal gradable covering of $Q$. Suppose first that $Q$ is of Dynkin type. Then $\tilde{Q}^{\hspace{.5pt}\rm op}$ is of Dynkin type with ${\rm rep}^-(\tilde{Q}^{\, \rm op})={\rm rep}^b(\tilde{Q}^{\hspace{.5pt}\rm op})$.
It is known that $\Ga_{D^b({\rm rep}^b(\tilde{Q}^{\hspace{.5pt}\rm op}))}$ is connected of shape $\Z\tilde{Q}$; see \cite[(5.6)]{Ha2}, and by Theorem \ref{ARQ}, $\Ga_{D^b({\rm mod}^b\hspace{-2pt}\it\Lambda)}\vspace{1pt}$ is connected of shape $\Z Q\hspace{.4pt}.$

Suppose that $Q$ is an oriented cycle of $n$ arrows. Then $r_{\hspace{-1pt}_Q}=n$ and $\tilde{Q}^{\hspace{.4pt}\rm op}$ is a double infinite path. By the dual of Theorem 5.17(1) in \cite{BLP}, $\Ga_{{\rm rep}^-(\tilde{Q}^{\hspace{.5pt}\rm op})}\vspace{1pt}$ has exactly one regular component $\mathcal{R}$, which contains only finite dimensional representations and is of shape $\Z \A_\infty$. Since the $P_{x^{\rm o}}$ and the $I_{x^{\rm o}}$ are all infinite dimensional, $\Ga_{{\rm rep}^-(\tilde{Q}^{\hspace{.5pt}\rm op})}\vspace{1pt}$ has no preprojective component, and by the dual of Proposition 3.6 in \cite{BLP}, none of the $I_{x^{\rm o}}$ is the ending term of an almost split sequence in $\rep^-(\tilde{Q}^{\hspace{.4pt}\rm op}).$ This implies that the preinjective component $\mathcal{I}$ of $\Ga_{{\rm rep}^-(\tilde{Q}^{\hspace{.5pt}\rm op})}\vspace{-1pt}$ is a double infinite path containing only the $I_{x^{\rm o}}$ with $x\in \tilde{Q}$, and it is the connecting component of $\Ga_{D^b({\rm rep}^-(\tilde{Q}^{\hspace{.5pt}\rm op}))}\vspace{1pt}$. By Theorem \ref{ARQ}, the connected components of $\Ga_{D^b({\rm mod}^b\hspace{-2pt}\it\Lambda)}$ are the $\mf{F}_\pi(\mathcal{R})[s]$ with $s\in \Z_n$, which are of shape $\Z \A_\infty$ without non-perfect complexes, and the $\mf{F}_\pi(\mathcal{I})[s]$ with $s\in \Z_n$, which are double infinite paths with only simple complexes.

Let now $Q$ be of type $\tilde{\A}_n$ with $0<r_{\hspace{-1pt}_Q}<n$. In particular, $Q$ is not an oriented cycle, and hence, $\tilde{Q}$ is of type $\A_\infty^\infty$ with no infinite path. In this case, ${\rm rep}^-(\tilde{Q})={\rm rep}^b(\tilde{Q})$, and consequently, $\Ga_{D^b({\rm mod}^b\hspace{-1.8pt}\it\Lambda)}\hspace{-1.5pt}$ contains only perfect complexes. By the dual of Theorem 5.17(2) in \cite{BLP}, $\Ga_{{\rm rep}^-(\tilde{Q}^{\hspace{.5pt}\rm op})}\vspace{1pt}$ has exactly two regular components $\mathcal{R}$ and $\mathcal{L}$, both are of shape $\Z\A_\infty$. Moreover, by the dual of Proposition 7.9(1), the connecting component $\mathcal{C}$ of $\Ga_{D^b({\rm rep}^-(\tilde{Q}^{\hspace{.5pt}\rm op}))}\vspace{1pt}$ is of shape $\Z \tilde{Q}$. By Theorem \ref{ARQ}, the components of $\Ga_{D^b({\rm mod}^b\hspace{-2pt}\it\Lambda)}$ are the $\mf{F}_\pi(\mathcal{R})[s]$ and the
$\mf{F}_\pi(\mathcal{L})[s]$ with $s\in \Z_{\,r_{\hspace{-1pt}_Q}}$, which are of shape $\Z \A_\infty$, and the $\mf{F}_\pi(\mathcal{C})[s]$ with $s\in \Z_{\,r_{\hspace{-1pt}_Q}}$, which are of shape $\Z \tilde{Q}$.

Suppose finally that $Q$ is neither of Dynkin type nor of type $\tilde{\A}_n$ with $r_{\hspace{-1pt}_Q}>0$. If $Q$ is gradable, then $\tilde{Q}^{\hspace{.4pt}\rm op}$ is finite but not of Dynkin type. Otherwise, $Q$ is not of type $\tilde{\A}_n$ for any $n\ge 1$, and consequently, $\tilde{Q}$ is infinite but not of infinite Dynkin types $\A_\infty$, $\A_\infty^\infty,$ $\mathbb{D}_\infty$. In either case, by the dual of Theorem 6.4 in \cite{BLP}, $\Ga_{{\rm rep}^-(\tilde{Q}^{\hspace{.4pt}\rm op})}\vspace{.5pt}$ contains infinitely many regular components. By Theorem \ref{ARQ}(2), $\Ga_{D^b({\rm mod}^b\hspace{-2pt}\it\Lambda)}\vspace{.5pt}$ has infinitely many connected components. The proof of the theorem is completed.

\medskip

\noindent{\sc Remark.} The second part of Theorem \ref{ARC-4} (2) was partially obtained in \cite[(5.7)]{HKR}.

\end{document}